\documentclass[11pt]{article}
\usepackage{amssymb}
\usepackage{amsfonts}
\usepackage{graphicx,color,epsfig}
\usepackage{amsmath}
\usepackage[utf8]{inputenc}

 \usepackage[left=2cm,right=2cm,top=2cm,bottom=3cm]{geometry}
 \usepackage{makeidx}
 \usepackage{verbatim}
 \usepackage{latexsym}

 \usepackage{setspace}
 \usepackage{bm}

\newcommand{\bs}{\boldsymbol}

\makeatletter \@addtoreset{equation}{section} \makeatother

\newtheorem{theorem}{Theorem}
\newtheorem{lemma}{Lemma}
\newtheorem{corollary}{Corollary}
\newtheorem{remark}{Remark}
\newtheorem{definition}{Definition}
\newtheorem{proposition}{Proposition}
\newtheorem{assumption}{Assumption}

\def\tr{\mathrm{Tr}}

\def\var{\mathbf{Var}}

\def\ex{\mathbb{E}}
\def\im{\mathrm{Im}}

\begin{document}
	\pagestyle{plain}
	\title{On the behaviour of large empirical autocovariance matrices between the past and the future.}
	\author{P. Loubaton $^{1}$, D. Tieplova $^{1,2}$ \\
$^{1}$ Laboratoire d'Informatique Gaspard Monge, UMR 8049, Universit\'e Paris-Est Marne la Vall\'ee, France \\
$^{2}$ Dept. of Theoretical Physics, Institute for Low Temperature Physics and Engineering, Ukraine}

	\date{This work is partially funded by the B\'ezout Labex, funded by ANR, reference ANR-10-LABX-58, 
and by the ANR Project HIDITSA, reference ANR-17-CE40-0003. }
	
	\setcounter{page}{1} \maketitle

        {\bf Abstract} The asymptotic behaviour of the distribution of the squared singular values of the sample autocovariance matrix between the past and the future of a high-dimensional complex Gaussian uncorrelated sequence is studied. Using Gaussian tools, it is established the distribution behaves as a deterministic probability measure 
whose support ${\cal S}$ is characterized. It is also established that the singular values to the square are almost surely located in a neighbourhood of ${\cal S}$. \\

        \section{Introduction.}
        \subsection{The addressed problem and the results.}
        In this paper, we consider a sequence of integer $(M(N))_{N \geq 1}$, and positive definite $M(N) \times M(N)$ hermitian matrices $(R_N)_{N \geq 1}$.
        For each $N$, we define an independent identically distributed sequence $(y_n)_{n \geq 1}$ (depending on $N$) of zero mean complex Gaussian $M(N)$--dimensional random vectors such that $y_n = R_N^{1/2} \xi_n$ where the  components
        of the $M$--dimensional vector $\xi_n$ are complex Gaussian standard i.i.d. random variables (i.e.
        their real and imaginary parts are i.i.d. and $\mathcal{N}(0,1/2)$ distributed).
        If $L$ is a fixed integer, we consider the 2 block-Hankel $ML \times N$ matrices $W_{p,N}$ and $W_{f,N}$ defined
        by
        \begin{equation}
\label{eq:def-Wp}
W_{p,N} = \frac{1}{\sqrt{N}} \, Y_{p,N} =  \frac{1}{\sqrt{N}} \left( \begin{array}{ccccc} y_1 & y_2 & \ldots & y_{N-1} & y_N \\
                                            y_2 & y_3 & \ldots & y_N & y_{N+1} \\
                                             \vdots & \vdots & \vdots & \vdots & \vdots \\
                                             \vdots & \vdots & \vdots & \vdots & \vdots \\
                                             y_L & y_{L+1} & \ldots & y_{N+L-2} & y_{N+L-1} \\
\end{array} \right)
\end{equation}
and 
\begin{equation}
\label{eq:def-YLf}
W_{f,N} = \frac{1}{\sqrt{N}}  \, Y_{f,N} =  \frac{1}{\sqrt{N}} \left( \begin{array}{ccccc} y_{L+1} & y_{L+2} & \ldots & y_{N-1+L} & y_{N+L} \\
                                            y_{L+2} & y_{L+3} & \ldots & y_{N+L} & y_{N+L+1} \\
                                             \vdots & \vdots & \vdots & \vdots & \vdots \\
                                             \vdots & \vdots & \vdots & \vdots & \vdots \\
                                             y_{2L} & y_{2L+1} & \ldots & y_{N+2L-2} & y_{N+2L-1} \\
\end{array} \right)
\end{equation}
and study the behaviour of the empirical eigenvalue distribution $\hat{\nu}_N$ of the $ML \times ML$ matrix
$W_{f,N} W_{p,N}^{*} W_{p,N} W_{f,N}^{*}$ in the asymptotic regime where $M$ and $N$ converge towards $+\infty$
in such a way that
\begin{equation}
\label{eq:asymptotic-regime}
c_N = \frac{ML}{N} \rightarrow c_*, \, c_* > 0
\end{equation}
Using Gaussian tools, we evaluate the asymptotic behaviour of the resolvent $Q_N(z) = (W_{f,N} W_{p,N}^{*} W_{p,N} W_{f,N}^{*} - z I)^{-1}$, and establish that the sequence $(\hat{\nu}_N)_{N \geq 1}$ has the same almost sure asymptotic behaviour than
a sequence $(\nu_N)_{N \geq 1}$ of deterministic probability measures. In the following, $\nu_N$ will be referred
to as the deterministic equivalent of $\hat{\nu}_N$. We evaluate the Stieltjes transform of $\nu_N$, characterize its support, study the
properties of its density, and eventually establish that almost surely, for $N$ large enough, all the eigenvalues of $W_{f,N} W_{p,N}^{*} W_{p,N} W_{f,N}^{*}$
are located in a neighbourhood of the support of $\nu_N$.         
\subsection{Motivation}
Matrix $W_{f,N} W_{p,N}^{*} = \frac{Y_{f,N}  Y_{p,N}^{*}}{N}$ represents the traditional empirical estimate of the
autocovariance matrix $R^{L}_{f|p,y}$ between the past and the future of $y$ defined as
$$
R^{L}_{f|p,y} = \mathbb{E} \left[ \left( \begin{array}{c} y_{n+L} \\  y_{n+L+1} \\ \vdots \\ y_{n+2L-1} \end{array} \right) \left(
y_n^{*}, y_{n+1}^{*}, \ldots, y_{n+L-1}^{*} \right) \right]
$$
This matrix plays a key role in statistical inference problems related to multivariate time series with rational spectrum. In order to explain
this, we consider a $M$--dimensional multivariate time series $(v_n)_{n \in \mathbb{Z}}$ generated as 
\begin{equation}
\label{eq:dynamic-factor-model}
v_n = u_n  + y_n
\end{equation}
where $(y_n)_{n \in \mathbb{Z}}$ is as above  a Gaussian "noise" term such that $\mathbb{E}(y_{n+k} y_n^{*}) = R \delta_k$ for some 
unknown positive definite matrix $R$, and where $(u_n)_{n \in \mathbb{Z}}$ is a "useful" non observable Gaussian signal with rational spectrum. $u_n$ can thus
be represented as
\begin{equation}
\label{eq:state-space}
x_{n+1}  =  A x_n + B \omega_n, \;  u_n  = C x_n + D \omega_n
\end{equation}
where $(\omega_n)_{n \in \mathbb{Z}}$ is a $K \leq M$--dimensional white noise sequence ($\mathbb{E}(\omega_{n+k} \omega_n^{*}) = I_K \, \delta_k$),
$A$ is a deterministic $P \times P$ matrix whose spectral radius $\rho(A)$ is strictly less than 1, and where $B,C,D$ are deterministic matrices.
The $P$-dimensional Markovian sequence $(x_n)_{n \in \mathbb{Z}}$ is called the state-space sequence associated to (\ref{eq:state-space}).
The state space representation (\ref{eq:state-space}) is said to be minimal if the dimension $P$ of the state space sequence is minimal.
Given the autocovariance sequence  $(R_{u,n})_{n \in \mathbb{Z}}$ of $u$ (i.e. $R_{u,n} = \mathbb{E}(u_{k+n}u_k^*)$ for each $n$),
the so-called stochastic realization problem of $(u_n)_{n \in \mathbb{Z}}$ consists in characterizing all the minimal state space representations (\ref{eq:state-space}) of $u$, or equivalently in identifying all the minimum Mac-Millan degree \footnote{The Mac-Millan degree of a rational matrix-valued
  function $\Phi$ is defined as the minimal dimension of the matrices $A$ for which $\Phi(z)$ can be represented as $ D + C(zI - A)^{-1} B$} matrix-valued function $\Phi(z) = D + C(zI - A)^{-1} B$ such that
$\rho(A)  < 1$ and
\begin{equation}
  \label{eq:spectral-factorization}
S_u(e^{2i\pi f}) = \sum_{n \in \mathbb{Z}} R_{u,n} e^{-2i\pi n f} = \Phi(e^{2 i \pi f}) \Phi(e^{2 i \pi f})^{*}
\end{equation}
for each $\nu$. Such a function $\Phi$ is called a minimal degree causal spectral factorization of $S_u$. We refer the reader to \cite{lindquist-picci-book-2015} or \cite{vanoverschee-demoor} for more details. \\

The identification of $P$ and of matrices $C$ and $A$ is based on the observation that the autocovariance sequence
of $u$ can be represented as
\begin{equation}
  \label{eq:expre-autocovariance}
R_{u,n}  = \mathbb{E}(u_{n+k} u_n^*) = C A^{n-1} G
\end{equation}
for each $n \geq 1$, where the 3 matrices $(A,C,G)$ are unique up to similarity transforms, thus showing that the matrices $C$ and $A$
associated to a minimal realization are uniquely defined (up to a similarity). Moreover, the autocovariance matrix $R^{L}_{f|p,u}$
between the past and the future of $u$ can be written as
\begin{equation}
\label{eq:factorisation-hankel}
R^{(L)}_{f|p,u} = \mathcal{O}^{(L)} \, \mathcal{C}^{(L)}
\end{equation}
where matrix $\mathcal{O}^{(L)}$ is the $ML \times P$ "observability" matrix 
\begin{equation}
\label{eq:def-OL}
\mathcal{O}^{(L)} = \left( \begin{array}{c} C \\ C A \\ \vdots \\  C A^{L-1} \end{array} \right)
\end{equation}
and matrix $\mathcal{C}^{(L)}$ is the $P \times ML$ "controllability" matrix  
\begin{equation}
\label{eq:def-CL}
\mathcal{C}^{(L)} = \left(  A^{L-1} G,  A^{L-2} G, \ldots,  G \right)
\end{equation}
For each $L \geq P$, the rank of $R^{(L)}_{f|p,u}$ remains equal to $P$, and each minimal rank factorization of $R^{(L)}_{f|p,u}$
can be written as (\ref{eq:factorisation-hankel}) for some particular triple $(A,C,G)$. In particular, if
$R^{(L)}_{f|p,u} = \Theta \Gamma \tilde{\Theta}^{*}$ is the singular value decomposition of $R^{(L)}_{f|p,u}$, matrix
$\Theta \Gamma^{1/2}$ coincides with the observability matrix $\mathcal{O}^{(L)}$ of a pair $(C,A)$. $C$ and $A$ are immediately
obtained from the knowledge of the structured matrix $\mathcal{O}^{(L)}$. This discussion shows that the evaluation of $P$, $C$ and $A$ from
the autocovariance sequence of $u$ is an easy problem. We mention that, while $C$ and $A$ are essentially unique, there exist in general more than
1 pair $(B,D)$ for which (\ref{eq:state-space}) holds because the minimal degree spectral factorization problem (\ref{eq:spectral-factorization}) has
more than 1 solution. We refer the reader to \cite{lindquist-picci-book-2015} or \cite{vanoverschee-demoor}. \\

We notice that as $(y_n)_{n \in \mathbb{Z}}$ in (\ref{eq:dynamic-factor-model}) is an uncorrelated sequence, it holds that
$R_{v,n} = \mathbb{E}(v_{n+k} v_k^{*})$ coincides with $R_{u,n}$ for each $n \geq 1$. Therefore, $P$ and matrices
$C$ and $A$ can still be identified from the autocovariance sequence of the noisy version $v$ of $u$. In practice, however,
the exact autocovariance sequence $(R_{v,n})_{n \geq 1}$ is in general unknown, and it is necessary to estimate $P$ and $(C,A)$ from
the sole knowledge of $N$ samples $v_1=u_1+y_1, v_2=u_2+y_2, \ldots, v_{N}=u_N+y_N$. For this, $P$ is first estimated as the number of
significant singular values of the empirical estimate  $\hat{R}^{L}_{f|p,v}$ of the true matrix $R^{L}_{f|p,v} = R^{L}_{f|p,u}$ defined by
$$
\hat{R}^{L}_{f|p,v} = \frac{V_{f,N} V_{p,N}^{*}}{N}
$$
where $V_{f,N}$ and $V_{p,N}$ are defined in the same way than $Y_{f,N}$ and $Y_{p,N}$. If $(\hat{\gamma}_p)_{p=1, \ldots, P}$ and $\hat{\Theta} = (\hat{\theta}_1, \ldots, \hat{\theta}_P)$ are the $P$ largest singular values and corresponding 
left singular vectors of matrix $\hat{R}^{(L)}_{f|p,v}$, and if $\hat{\Gamma}$ is the $P \times P$ diagonal matrix 
with diagonal entries  $(\hat{\gamma}_p)_{p=1, \ldots, P}$, $ML \times P$ matrix 
$\hat{\mathcal{O}}^{(L)} = \hat{\Theta} \hat{\Gamma}^{1/2}$ is an estimator of an observability matrix $\mathcal{O}^{(L)}$. $\hat{\mathcal{O}}^{(L)}$
has not necessarily the structure of an observability matrix, but it is easy to estimate $A$ by finding the minimum of the quadratic 
fuction
$$
\left\| \hat{\mathcal{O}}^{(L)}_{\mathrm{down}}  A  - \hat{\mathcal{O}}^{(L)}_{\mathrm{up}} \right\|
$$
where the operator "down" (resp. "up") suppresses the last (resp. the first) $M$ rows from $ML \times P$ matrix 
$\hat{\mathcal{O}}^{(L)}$. This approach provides a consistent estimate of $P,C,A$ when $N \rightarrow +\infty$ while $M$, $L$ and $P$ are fixed
parameters. We refer the reader to \cite{ciuso-picci-2004} for a detailed  analysis of this statistical inference scheme.  \\

If $M$ is large and that the sample size $N$ cannot be arbitrarily larger than $M$, the ratio $ML/N$ may not be small enough to make reliable
the above statistical analysis. It is thus relevant to study the behaviour of the above estimators in asymptotic regimes where $M$ and $N$ both converge towards $+\infty$ in such a way that $\frac{ML}{N}$ converges towards a non zero constant. In this context, the truncated singular value decomposition of $\hat{R}^{(L)}_{f|p,v}$ does not provide a consistent estimate of an observability matrix $\mathcal{O}^{(L)}$, and it appears relevant to study the largest singular values and corresponding singular vectors of $\hat{R}^{(L)}_{f|p,v}$ when $M$ and $N$ both converge towards $+\infty$, and to precise how they are related to an observability matrix $\mathcal{O}^{(L)}$. \\

Without formulating specific assumptions on $u$, this problem seems very complicated. In the past, a number of works addressed
high-dimensional inference schemes based on the eigenvalues and eigenvectors of the empirical covariance matrix of the observation
(see e.g. \cite{nadakuditi-edelman-2008}, \cite{mestre2008modified}, \cite{nadakuditi-silverstein-2010},
  \cite{hachem-loubaton-et-al-jmva-2013}, \cite{vinogradova-couillet-hachem-1}, \cite{vinogradova-couillet-hachem-2},
  \cite{couillet-jmva-2015}, \cite{vallet-mestre-loubaton-2016}) when the useful signal lives in a low-dimensional
  deterministic subspace. Using results related to spiked large random matrix models (see e.g. \cite{benaych-nadakuditi-ann-math}
  \cite{benaych-rao-2}, \cite{paul-2007}), based on perturbation technics, a number of important statistical problems could be addressed
  using large random matrix theory technics. Our ambition is to follow the same kind of approach to address the estimation problem
  of $P,A,C$ when $u$ satisfies some low rank assumptions. The first part of this program is to study the asymptotic behaviour of the
  singular values of the
  empirical autocovariance matrix in the absence of signal $W_{f,N} W_{p,N}^{*} = \frac{Y_{f,N}  Y_{p,N}^{*}}{N}$. As the singular values
  of $W_{f,N} W_{p,N}^{*}$ are the square roots of the eigenvalues of $W_{f,N} W_{p,N}^{*} W_{p,N} W_{f,N}^{*}$, this is precisely the topic of the
  present paper. Using the obtained results, it should be
  possible to use a perturbation approach in order to evaluate the behaviour of the largest singular values and corresponding left singular vectors in the
  presence of a useful signal, and to deduce from this some improved performance scheme for estimating $P,C,A$.

  \subsection{On the literature.}
  The large sample behaviour of high-dimensional autocovariance matrices was comparatively less studied than the
  high-dimensional covariance matrices. We first mention \cite{jin-bai-el-al-2014} which studied the asymptotic behaviour of
  the eigenvalue distribution of the hermitian matrix $\hat{R}_{\tau} + \hat{R}_{\tau}^{*}$ where $\hat{R}_{\tau}$ is defined as
  $\hat{R}_{\tau} = \frac{1}{N} \sum_{n=1}^{N} x_{n+\tau} x_{n}^{*}$ where $(x_n)_{n \in \mathbb{Z}}$ represents
  a $M$ dimensional non Gaussian i.i.d. sequence, the components of each vector $x_n$ being morever i.i.d. In particular,
  $\mathbb{E}(x_n x_n^{*}) = I$. It is proved that the empirical eigenvalue distribution of  $\hat{R}_{\tau} + \hat{R}_{\tau}^{*}$ converges towards a
  limit distribution independent from $\tau \geq 1$. Using finite rank perturbation technics of the resolvent of the matrix under consideration, the Stieltjes transform of this distribution was shown 
to satisfy a polynomial degree 3 equation. Solving this equation led to an explicit
expression of the  probability density of the limit distribution. \cite{liu-aue-paul-2015} extended these
  results to the case where $(x_n)_{n \in \mathbb{Z}}$ is a non Gaussian linear process $x_n = \sum_{l=0}^{+\infty} A_l z_{n-l}$
  where $(z_n)_{n \in \mathbb{Z}}$ is i.i.d., and where matrices $(A_l)_{l \geq 0}$ are simultaneously diagonalizable. The limit
  eigenvalue distribution was characterized through its Stieltjes transform that is obtained by integration of a certain kernel, itself solution of an integral equation. The proof was based on the observation that in the Gaussian case, 
the correlated vectors $(x_n)_{n \in \mathbb{Z}}$ can be replaced by independent ones using a classical 
frequency domain decorrelation procedure. The results were generalized in the non Gaussian case using the generalized Lindeberg principle. We also mention \cite{batta-bose-2016} (see also the book 
\cite{bose-batta-book}) where the existence of a limit distribution of any symmetric polynomial 
of $(\hat{R}_{\tau}, \hat{R}_{\tau}^{*})_{\tau \in T}$ for some finite set $T$ was proved using the moment method when $x$ is a linear non Gaussian process. \cite{li-pan-yao-jmva-2015} studied the asymptotic behaviour of 
matrix $\hat{R}_{\tau} \hat{R}_{\tau}^{*}$ when $(x_n)_{n \in \mathbb{Z}}$ represents
  a $M$ dimensional non Gaussian i.i.d. sequence, the components of each vector $x_n$ being morever i.i.d. 
Using finite rank perturbation technics, it was shown that the empirical eigenvalue distribution converges 
towards a limit distribution whose Stieltjes transform is solution of a degree 3 polynomial equation. As 
in \cite{jin-bai-el-al-2014}, this allowed to obtain the expression of the corresponding probability 
density function. Using combinatorial technics, \cite{li-pan-yao-jmva-2015} also established that almost surely, 
for large enough dimensions, all the eigenvalues of $\hat{R}_{\tau} \hat{R}_{\tau}^{*}$ are located 
in a neighbourhood of the support of the limit eigenvalue distribution. We finally mention that \cite{li-wang-yao-ann-stat-2016} used the results in \cite{li-pan-yao-jmva-2015} in order to study the largest eigenvalues and corresponding eigenvectors of $\hat{R}_{\tau} \hat{R}_{\tau}^{*}$ when the observation contains a certain spiked useful signal that is more specific than the signals
signals $(u_n)_{n \in \mathbb{Z}}$ that motivated the present paper. \\

We now compare the results of the present paper with the content of the above previous works. We first study 
a matrix that is more general than $\hat{R}_{\tau} \hat{R}_{\tau}^{*}$. While we do not consider linear processes
here, we do not assume that the covariance matrix of the i.i.d. sequence $(y_n)_{n \in \mathbb{Z}}$ is reduced 
to $I$ as in \cite{li-pan-yao-jmva-2015}. This in particular implies that the Stieltjes transform of the 
deterministic equivalent $\nu_N$ of $\hat{\nu}_N$ cannot be evaluated in closed from. Therefore, a dedicated 
analysis of the support and of the properties of $\nu_N$ is provided here. We also mention 
that in contrast with the above papers, we characterize the asymptotic behaviour of the resolvent of 
matrix $W_{f,N} W_{p,N}^{*} W_{p,N} W_{f,N}^{*}$ while the mentionned previous works only studied 
the normalized trace of the resolvent of the matrices under consideration. Studying the full resolvent 
matrix is necessary to address the case where a useful spiked signal $u$ is added to the noise $y$. 
We notice that the above papers addressed the non Gaussian case while we consider the case where 
$y$ is a complex Gaussian i.i.d. sequence. This situation is of course simpler in that various Gaussian 
tools are available, but appears to be 
relevant because in the context of the present paper, $y$ is indeed supposed to represent some additive 
noise, which, in a number of contexts, is Gaussian. In any case, it should be possible to 
extend the present results to the non Gaussian case by using the Lindeberg principle or some 
interpolation scheme. \\

We finally mention that some of the results of this paper may be obtained by adapting general recent results
devoted to the study of the spectrum of hermitian polynomials of GUE matrices and deterministic matrices
(see \cite{belinschi-capitaine-2017} and \cite{male-2012}). If we denote by
$Z_N$ the $M \times (N+2L-1)$ matrix $Z_N = (y_1, \ldots, y_{N+2L-1})$, then
$Z_N$ can be written as $Z_N = R_N^{1/2} X_N$ where the entries of $X_N$ are i.i.d. complex Gaussian
standard variables. Each $M \times M$ block  $\Sigma_{N,k,l}$ ($1 \leq k,l \leq L$) of $\Sigma_{N} = W_{f,N} W_{p,N}^{*} W_{p,N} W_{f,N}^{*}$ is clearly a
polynomial of $X_N, X_N^{*}$ and various $M \times M$ and $M \times (N+2L -1)$ deterministic matrices. Assume
that $M < N+2L-1$. In order to be back to
a polynomial of GUE matrices, it possible to consider the $L(N+2L-1) \times L(N+2L-1)$ matrix $\tilde{\Sigma}_{N}$ whose
$(N+2L -1) \times (N+2L -1)$ blocks are defined by
$$
\tilde{\Sigma}_{N,k,l} = \left( \begin{array}{cc} \Sigma_{N,k,l} & 0 \\ 0 & 0  \end{array} \right)
$$
It is clear that apart $0$, the eigenvalues of $\tilde{\Sigma}_N$ coincide with those of $\Sigma_N$. If $\tilde{X}_N$ is any
$(N+2L-1) \times (N+2L -1)$ matrix with i.i.d. complex Gaussian standard entries whose $M$ first rows
coincide with $X_N$, then, it is easily seen that each block of $\tilde{\Sigma}_N$ coincides with
a hermitian polynomial of $\tilde{X}_N, \tilde{X}_N^{*}$ and deterministic $(N+2L-1) \times (N+2L-1)$ matrices such as 
$$
\tilde{R}_N =  \left( \begin{array}{cc} R_N  & 0 \\ 0 & 0  \end{array} \right)
$$
Expressing $\tilde{X}_N$ as the sum of its hermitian and anti-hermitian
parts, we are back to study the behaviour of the eigenvalues of a matrix whose blocks are hermitian polynomials
of 2 independent GUE matrices and of $(N+2L-1) \times (N+2L-1)$ deterministics matrices. Extending
Proposition 2.2 and Theorem 1.1 in \cite{belinschi-capitaine-2017} to block matrices (as in
Corollary 2.3 in \cite{male-2012}) would lead to the conclusion that $\hat{\nu}_N$ has a deterministic
equivalent $\nu_N$ and that the eigenvalues of $ W_{f,N} W_{p,N}^{*} W_{p,N} W_{f,N}^{*}$ are located in the
neighbourhood of the support of $\nu_N$. While this last consequence would avoid the use of the specific approach
used in section \ref{sec:localization-eigenvalues} of the present paper, the existence of $\nu_N$ is not a sufficient information.
$\nu_N$ should of course be characterized through its Stieltjes transform, and we believe that the adaptation of
Proposition 2.2 and Theorem 1.1 in \cite{belinschi-capitaine-2017} is not the most efficient approach. 
\subsection{Overview of the paper.}
As the entries of matrices $W_{p,N}$ and $W_{f,N}$ are correlated, approaches based on finite rank perturbation of the resolvent $Q_N(z)$ of matrix $W_{f,N}W_{p,N}^{*}W_{p,N}W_{f,N}^{*}$, usually used when independence assumptions hold, are not the most efficient in our context. We rather propose to use Gaussian tools, i.e. integration by parts formula in conjunction with the Poincar\'e-Nash inequality 
(see e.g. \cite{pastur-shcherbina-book}), because they are robust to correlation of the matrix entries. 
Moreover, as the entries of $W_{f,N}W_{p,N}^{*}W_{p,N}W_{f,N}^{*}$ are biquadratic functions of 
$y_1, \ldots, y_{N+2L-1}$, we rather use the well-known linearization trick that consists in 
studying the resolvent ${\bf Q}_N(z)$ of the $2ML \times ML$ hermitized version 
$$
\begin{pmatrix}
0 & W_{f,N} W_{p,N}^*\\
W_{p,N} W_{f,N}^*&0
\end{pmatrix}
$$
of matrix $W_{f,N} W_{p,N}^*$. As is well known, the first $ML \times ML$ diagonal block of
${\bf Q}_N(z)$ coincides with $z Q_N(z^{2})$. Therefore, we characterize the asymptotic behaviour of 
${\bf Q}_N(z)$, and deduce from this the results concerning $Q_N(z)$. The  hermitized version 
is this time a quadratic function of $y_1, \ldots, y_{N+2L-1}$, and the Gaussian calculus that is needed in order 
to study ${\bf Q}_N(z)$ appears much simpler than if $Q_N(z)$ was evaluated directly. \\

In section \ref{sec:poincare-nash}, we evaluate the variance of useful functionals for ${\bf Q}_N(z)$ using
the Poincaré-Nash inequality. In section \ref{sec:stieljes-lemmas}, we establish some useful lemmas 
related to certain Stieltjes transforms. In section \ref{subsec:ipp}, we use the integration by parts formula 
to establish that $\mathbb{E}({\bf Q}_N(z))$ behaves as $I_{2L} \otimes {\bf S}_N(z)$ where ${\bf S}_N(z)$ is defined by 
$$
 \mathbf{S}_N(z)=-\left(\dfrac{c_N\bm{\alpha}_N(z)}{1-c_N^2\bm{\alpha}(z)^2}R_N+z I_M\right)^{-1} 
$$
where $\bm{\alpha}_N(z)$ is defined by $\bm{\alpha}_N(z) = \frac{1}{ML} \mathrm{Tr}\mathbb{E}({\bf Q}_{N,pp}(z))(I_L \otimes R_N)$ where ${\bf Q}_{N,pp}(z)$ represents the first $ML \times ML$ diagonal 
block of ${\bf Q}_N(z)$. We deduce from this that 
$$
E(Q_N(z)) = S_N(z) + \Delta_N(z)
$$
where $S_N(z) = -\left(z I_M + \dfrac{c_N z \alpha_N(z)}{1-c_N^2\alpha_N(z)^2}R_N \right)^{-1}$,
$\alpha_N(z) =  \frac{1}{ML} \mathrm{Tr}\mathbb{E}(Q_{N}(z))(I_L \otimes R_N)$, and where 
$\Delta_N(z)$ is an error term such that
$$
\left| \frac{1}{ML} \mathrm{Tr} \, \Delta_N(z) \right| \leq \frac{1}{N^{2}} P_1(|z|) P_2(\frac{1}{\im(z)})
$$
for each $z \in \mathbb{C}^{+}$, where $P_1$ and $P_2$ are 2 polynomials whose degrees and coefficients do not depend on 
$N$. Using this, we prove in section \ref{sec:deterministic-equivalent} that for each $z \in \mathbb{C}^{+}$, 
$$
\frac{1}{ML} \mathrm{Tr}\mathbb{E}\left[Q_N(z) - I_L \otimes T_N(z)\right]F_N  \rightarrow 0
$$
where $(F_N)_{N \geq 1}$ is any deterministic sequence of matrices such that 
$\sup_{N} \| F_N \| < +\infty$, and where $T_N(z)$ is defined by 
$$
T_N(z) = - \left(zI_M+ \dfrac{zc_Nt_N(z)}{1- zc_N^2t_N^2(z)}R_N\right)^{-1},
$$
$t_N(z)$ being the unique solution of the equation
\begin{equation}
\label{eq:equation-t}
t_N(z)=\dfrac{1}{M}\tr R_N\left(-zI_M-\dfrac{zc_Nt_N(z)}{1- zc_N^2t_N^2(z)}R_N\right)^{-1} 
\end{equation}
such that $t_N(z)$ and $z t_N(z)$ belong to $\mathbb{C}^{+}$ when $z \in \mathbb{C}^{+}$. $t_N(z)$ and $T_N(z)$ are shown to 
coincide with the Stieltjes transforms of a scalar measure $\mu_N$ and of a $M \times M$ positive matrix valued measure $\nu_N^{T}$ respectively, 
and it is proved that $\nu_N = \frac{1}{M} \mathrm{Tr}(\nu_N^{T})$ is a probability measure
such that $\hat{\nu}_N - \nu_N \rightarrow 0$ weakly almost surely. $\nu_N$ is referred to as the 
deterministic equivalent of $\hat{\nu}_N$. In section \ref{sec:study-nuN},
we study the properties and the support of $\nu_N$, or equivalently of $\mu_N$ because the 2 measures are absolutely continuous one with respect to each other. For this, we study the behaviour of 
$t_N(z)$ when $z$ converges towards the real axis. For each $x > 0$, 
the limit of $t_N(z)$ when $z \in \mathbb{C}^{+}$ converges towards $x$ exists and is finite. If $c_N \leq 1$, we deduce from this that $\nu_N$ is absolutely continuous w.r.t. the Lebesgue measure. The corresponding density $g_N(x)$ is real analytic on $\mathbb{R}^{+}$, and converges towards $+\infty$ when $x \rightarrow 0, x > 0$. If $c_N < 1$, it holds that 
$g_N(x) = \mathcal{O}(\frac{1}{\sqrt{x}})$ while $g_N(x) = \mathcal{O}(\frac{1}{x^{2/3}})$ if $c_N = 1$. If $c_N > 1$, $\nu_N$ contains a Dirac mass at $0$ with weight $1 - \frac{1}{c_N}$ and an absolutely continuous component. In order to analyse the support of $\mu_N$ and $\nu_N$, we establish that the function $w_N(z)$ defined by 
$$
w_N(z) = z c_N t_N(z) - \frac{1}{c_N t_N(z)} 
$$
is solution of the equation $\phi_N(w_N(z)) = z$ for each $z \in \mathbb{C} - \mathbb{R}^{+}$
where $\phi_N(w)$ is the function defined by 
$$
\phi_N(w) = c_N w^{2} \, \frac{1}{M} \mathrm{Tr} R_N \left( R_N - w I \right)^{-1} \, \left( c_N \, \frac{1}{M} \mathrm{Tr} R_N \left( R_N - w I \right)^{-1} - 1 \right)
$$
Moreover, if we define $t_N(x)$ for $x > 0$ by the limit of $t_N(z)$ when $z \rightarrow x, z \in \mathbb{C}^{+}$, the equality $\phi_N(w_N(z)) = z$ is also valid on $\mathbb{R}^{+}$. We establish that if $x$ is outside the  support of $\mu_N$, then, it holds that 
$$
 \phi_N(w_N(x))  = x, \; \phi^{'}(w_N(x)) >  0, \; w_N(x) \, \frac{1}{M} \mathrm{Tr} R_N \left( R_N - w(x) I \right)^{-1} < 0
$$
This property allows to prove that apart $\{0 \}$ when $c_N > 1$, the support of $\mu_N$ is a union of intervals whose end points are the extrema of $\phi_N$ whose arguments verify $\frac{1}{M} \mathrm{Tr} R \left( R - w I \right)^{-1} < 0$. A sufficient condition on the eigenvalues of $R_N$ ensuring that the support of $\mu_N$ is reduced to a single interval is formulated. Using the 
Haagerup-Thornbjornsen approach (\cite{HT:05}), it is moreover proved in section \ref{sec:localization-eigenvalues} that for each $N$ large enough, all the eigenvalues of 
$W_{f,N} W_{p,N}^{*} W_{p,N} W_{f,N}^{*}$  lie in a neighbourhood the support of the deterministic equivalent $\nu_N$. The above results do not imply that $\hat{\nu}_N$ converges towards a limit distribution. In order to obtain this kind of result, some extra assumptions have to be 
formulated, such as the existence of a limit empirical eigenvalue distribution for $R_N$ when 
$N \rightarrow +\infty$. If the relevant conditions are met, $\nu_N$, and therefore 
$\hat{\nu}_N$, will converge towards a limit distribution whose Stieltjes transform can be 
obtained by replacing in the above results the empirical eigenvalue of $R_N$ by its limit. 
We do not present the corresponding results here because we believe that results that 
characterize the behaviour of $\nu_N$ for each $N$ large enough are more informative 
than the convergence towards a limit. \\

In section \ref{sec:free-probability}, we finally indicate that the use of free probability tools is an alternative
approach to characterize the asymptotic behaviour of $\hat{\nu}_N$. The results of section \ref{sec:free-probability} are based  on the following observations:
        \begin{itemize}
        \item Up to the zero eigenvalue, the eigenvalues of $W_{f,N}W_{p,N}^{*}W_{p,N} W_{f,N}^{*}$ coincide with the eigenvalues of $ W_{f,N}^{*}W_{f,N} W_{p,N}^{*}W_{p,N}$
        \item While the matrices $W_{f,N}^{*}W_{f,N}$ and $W_{p,N}^{*}W_{p,N}$ do not satisfy the conditions of the usual asymptotic freeness results, it turns out that they are almost surely asymptotically free. Therefore, the eigenvalue distribution
          of $ W_{f,N}^{*}W_{f,N} W_{p,N}^{*}W_{p,N}$ converges towards the free multiplicative convolution product of the limit distributions of  $W_{f,N}^{*}W_{f,N}$ and $W_{p,N}^{*}W_{p,N}$. These two distributions appear to coincide both with the limit distribution of the well known random matrix model
          $\frac{1}{N} X_N^{*} (I_L \times R_N) X_N$ where $X_N$ is a $ML \times N$ complex Gaussian random matrix with standard i.i.d. entries.  
        \end{itemize}
        The asymptotic freeness of  $W_{f,N}^{*}W_{f,N}$ and $W_{p,N}^{*}W_{p,N}$ appear to be a consequence of
        Lemma 6 in \cite{Evans_Tse_2000}. While this approach seems to be simpler than the use of the Gaussian
        tools proposed in the present paper, we mention that the above free probability theory arguments do not allow to
        study the asymptotic behaviour of the resolvent of $W_{f,N}W_{p,N}^{*}W_{p,N} W_{f,N}^{*}$. We recall that in order to evaluate the largest eigenvalues and corresponding eigenvectors of $W_{f,N}W_{p,N}^{*}W_{p,N} W_{f,N}^{*}$ in the presence
        of a useful signal, the asymptotic behaviour of the full resolvent in the absence of signal has to be available.

\section{Some notations, assumptions, and useful results.}
In the following, it is assumed that $L$ is a fixed parameter, and that $M$ and $N$ 
converge towards $+\infty$ in such a way that 
\begin{equation}
\label{eq:asymptotic-regime}
c_N = \frac{ML}{N} \rightarrow c_*, \, c_* > 0
\end{equation}
This regime will be referred to as $N \rightarrow +\infty$ in the following. In the regime 
(\ref{eq:asymptotic-regime}), $M$ should be interpreted as an integer $M=M(N)$ depending on $N$. 
The various matrices we have introduced above thus depend on $N$ and will be 
denoted $R_N, Y_{f,N},Y_{p,N}, \ldots$. In order to simplify the notations, the dependency 
w.r.t. $N$ will sometimes be omitted.

We recall that the resolvent $Q_N(z)$ of 
$W_{f,N} W_{p,N}^{*} W_{p,N} W_{f,N}^{*}$ is defined by 
\begin{equation}
\label{eq:def-QN}
Q_N(z) = \left( W_{f,N} W_{p,N}^{*} W_{p,N} W_{f,N}^{*} - z I \right)^{-1}
\end{equation}
As the direct study of $Q_N(z)$ is not obvious, we rather introduce the resolvent $\mathbf{Q}_N(z)$ of the $2ML\times2ML$ block matrix 
\begin{align*}
\mathbf{M}_N=	\begin{pmatrix}
0 & W_{f,N} W_{p,N}^*\\
W_{p,N} W_{f,N}^*&0
\end{pmatrix}.
\end{align*} 
It is well known that $\mathbf{Q}_N(z)$ can be expressed as 
\begin{equation}
\label{eq:expre-BQ}
\mathbf{Q}_N(z) = \left( \begin{array}{cc} z Q_N(z^{2}) & Q_N(z^{2}) W_{f,N} W_{p,N}^* \\
                                       W_{p,N} W_{f,N}^{*} Q_N(z^{2}) &  z \tilde{Q}_N(z^{2}) \end{array} \right)
\end{equation}
where $\tilde{Q}_N(z)$ is the resolvent of matrix $W_{p,N} W_{f,N}^* W_{f,N} W_{p,N}^{*}$. As shown below, it is rather easy to 
evaluate the asymptotic behaviour of $\mathbf{Q}_N(z)$ using the Poincaré-Nash inequality and the integration 
by part formula (see Propositions \ref{prop:poincare} and \ref{prop:integration-by-parts} below). Formula (\ref{eq:expre-BQ}) will then provide all the necessary information on
$Q_N(z)$. \\

In the following, every $2ML\times 2ML$ matrix $\mathbf{G}$ will be written as 
\begin{align*}
\mathbf{G}=\begin{pmatrix}
\mathbf{G_{pp}}&\mathbf{G_{pf}}\\
\mathbf{G_{fp}}&\mathbf{G_{ff}}
\end{pmatrix},
\end{align*}  
where the 4 matrices $(\mathbf{G}_{i,j})_{i,j \in {p,f}}$ are $ML \times ML$. Sometimes, the blocks will be denoted $\mathbf{G}(pp)$, $\mathbf{G}(pf)$, 
....\\

We denote by $W_N$ the $2ML \times N$ matrix defined by 
\begin{align}
\label{eq:def-WN}
W_N=	\begin{pmatrix}
W_{p,N}\\
W_{f,N}
\end{pmatrix},
\end{align} 
Its elements $(W_{i,j}^{m})_{i \leq 2L, j \leq N, m \leq M}$ satisfy 
\begin{align*}
\mathbb{E}\{W_{i,j}^{m}W_{i^{\prime},j^{\prime}}^{m^{\prime}}\}=\dfrac{1}{N}R_{mm^{\prime},N}\delta_{i+j,i^{\prime}+j^{\prime}}.
\end{align*} 
where $W_{i,j}^m$ represents the element which  lies on the $(m+M(i-1))$-th line and $j$-th column for $1\le m\le M$, $1\le i\le 2L$ and $1\le j\le N$.  Similarly, $\mathbf{Q}_{i_1i_2}^{m_1m_2}$, where $1\le m_1,m_2\le M$ and $1\le i_1,i_2\le 2L$, represents the entry $(m_1+M(i_1-1)),(m_2+M(i_2-1))$ of $\mathbf{Q}$. For each $j=1, \ldots, N$,$\{w_{j}\}_{j=1}^N, \{w_{p,j}\}_{j=1}^N$ and $\{w_{f,j}\}_{j=1}^N$ are the column of matrices $W, W_p$ and $W_f$ respectively. For each $1\le i\le 2L$ and $1\le m\le M$, $\mathbf{f}_i^m$ represents the vector of the canonical basis of $\mathbb{C}^{2ML}$ with 1 at the index $m+(i-1)M$ and zeros elsewhere. In order to simplify the notations, we mention that if 
$i \leq L$, vector $\mathbf{f}_i^m$ may also represents the vector of the canonical basis of $\mathbb{C}^{ML}$ with 1 at the index $m+(i-1)M$ and zeros elsewhere. Vector $\mathbf{e}_j$ with $1\le j\le N$ represents the $j$~--th vector of the canonical basis of $\mathbb{C}^N$. 
Also for any integer $k$, $J_k$ is the $k \times k$ "shift" matrix defined by
\begin{equation}
  \label{eq:def-J}
  (J_k)_{ij} = \delta_{j-i, 1}
\end{equation}
In order to short the
notations, matrix $J^*_k$ is denoted $J^{-1}_k$ , although $J_k$ is of course not invertible.\\

By a nice constant, we mean a positive deterministic constant which does not depend on the dimensions $M$ and $N$ nor of the complex variable $z$. In the following, $\kappa$ will represent a generic nice constant whose value may
change from one line to the other. A nice polynomial $P(z)$ is a polynomial whose degree and coefficients
are nice constants. Finally, we will say that function $f_N(z)=\mathcal{O}_z(\alpha_N)$ if $z$ belongs to a domain $\Omega \in \mathbb{C}$ and there exist two nice polynomials $P_1$ and $P_2$ such that $f_N(z)\le \alpha_N P_1(|z|)P_2(\frac{1}{|\im z|})$ for each $z \in \Omega$. If $\Omega = \mathbb{C}^{+}$, we will just write
$f_N(z)=\mathcal{O}_z(\alpha_N)$ without mentioning the domain. We notice that 
if $P_1$, $P_2$ and $Q_1$, $Q_2$ are nice polynomials, then $P_1(|z|)P_2(\frac{1}{|\im z|})+Q_1(|z|)Q_2(\frac{1}{|\im z|})\le (P_1+Q_1)(|z|)(P_2+Q_2)(\frac{1}{|\im z|})$, from which we conclude that if functions $f_1$ and $f_2$ are $\mathcal{O}_z(\alpha)$ then also $f_1(z)+f_2(z)=\mathcal{O}_z(\alpha)$. \\

The sequence of covariance matrices $(R_N)_{N \geq 1}$ of $M$--dimensional vectors $(y_n)_{n=1, \ldots, N}$ 
is supposed to verify
\begin{equation}
\label{eq:hypothesis-R-bis}
a \, I \leq R_N < b \, I
\end{equation}
for each $N$, where $a>0$ and $b>0$ are 2 nice constants. $\lambda_{1,N} \geq \lambda_{2,N} \geq \ldots \geq \lambda_{M,N}$ represent the eigenvalues of $R_N$ arranged in the decreasing order and $f_{1,N}, \ldots, f_{M,N}$ denote the corresponding eigenvectors. Hypothesis 
(\ref{eq:hypothesis-R-bis}) is obviously equivalent to $\lambda_{M,N} \geq a$ and $\lambda_{1,N} \leq b$ for 
each $N$.  \\

The eigenvalues and eigenvectors of matrix $W_{f,N} W_{p,N}^{*} W_{p,N} W_{f,N}^{*}$ are denoted $\hat{\lambda}_{1,N} \geq \ldots \geq \hat{\lambda}_{M,N}$
and $\hat{f}_{1,N}, \ldots, \hat{f}_{M,N}$ respectively. \\

$\mathcal{C}_{c}^{\infty}(\mathbb{R},\mathbb{R})$ represents the set of all $\mathcal{C}^{\infty}$ real valued compactly supported functions defined on $\mathbb{R}$. \\ 

If $\xi$ is a random variable, we denote by $\xi^{\circ}$ the zero mean random variable
defined by 
\begin{equation}
\label{eq:def-chirond}
\xi^{\circ} = \xi - \mathbb{E} \xi
\end{equation}

We finally recall the 2 Gaussian tools that will be used in the sequel in order to 
evaluate the asymptotic behaviour of $\mathbf{Q}_N(z)$ and $Q_N(z)$. 
\begin{proposition}
  \label{prop:integration-by-parts}
	\textbf{(Integration by parts formula.)}
	Let $\xi=[\xi_1,\ldots,\xi_K]^T$ be a complex Gaussian random vector such that $\ex\{\xi\}=0$, $\ex\{\xi\xi^T\}=0$ and $\ex\{\xi\xi^*\}=\Omega$. If $\Gamma:(\xi)\mapsto\Gamma(\xi,\bar{\xi})$ is a $\mathcal{C}^1$ complex function polynomially bounded together with its derivatives, then 
	\begin{align}\label{integr}
	\ex\{\xi_i\Gamma(\xi)\}=\sum_{k=1}^{K}\Omega_{ik}\ex\left\{\dfrac{\partial\Gamma(\xi)}{\partial\bar{\xi}_k}\right\}.
	\end{align}
\end{proposition}

\begin{proposition}
  \label{prop:poincare}
	\textbf{(Poincaré-Nash inequality.)}
	Let $\xi=[\xi_1,\ldots,\xi_K]^T$ be a complex Gaussian random vector such that $\ex\{\xi\}=0$, $\ex\{\xi\xi^T\}=0$ and $\ex\{\xi\xi^*\}=\Omega$. If $\Gamma:(\xi)\mapsto\Gamma(\xi,\bar{\xi})$ is a $\mathcal{C}^1$ complex function polynomially bounded together with its derivatives, then, noting $\nabla_{\xi}\Gamma=[\frac{\partial\Gamma}{\partial\xi_1},\ldots,\frac{\partial\Gamma}{\partial\xi_K}]^T$ and $\nabla_{\bar{\xi}}\Gamma=[\frac{\partial\Gamma}{\partial\bar{\xi}_1},\ldots,\frac{\partial\Gamma}{\partial\bar{\xi}_K}]^T$
	\begin{align}\label{p-n}
	\var\{\Gamma(\xi)\}\le\ex\left\{\nabla_{\xi}\Gamma(\xi)^T\Omega\overline{\nabla_{\xi}\Gamma(\xi)}\right\}+\ex\left\{\nabla_{\bar{\xi}}\Gamma(\xi)^*\Omega\nabla_{\bar{\xi}}\Gamma(\xi)\right\}.
	\end{align}
\end{proposition}

\section{Use of the Poincaré-Nash inequality.}
\label{sec:poincare-nash}
In this paragraph, we control the variance of various functionals of $\mathbf{Q}_N(z)$ using the 
Poincaré-Nash inequality. For this, it appears useful to evaluate the moments of $\|W_N\|$.  The following
result holds. 
\begin{lemma}\label{moment}
	For any $l\in \mathbb{N}$, it holds that $\sup_{N \geq 1} \ex\{\|W_N\|^{2l}\}<+\infty.$
\end{lemma}

\textbf{Proof.}  We first remark that it is possible to be back to the case 
where matrix $R_N = I_M$. Due to the Gaussianity of the i.i.d. vectors $(y_n)_{n \geq 1}$, it exists 
i.i.d.  $\mathcal{N}_c(0,I_M)$ distributed vectors $(y_{iid,n})_{n \geq 1}$ such that $\mathbb{E}(y_{iid,n} y_{iid,n}^{*}) = I_M$ verifying $y_n = R_N^{1/2} y_{iid,n}$. From this, we obtain immediately that the $2ML \times N$ block Hankel matrix $W_{iid,N}$ built from $(y_{n,iid})_{n=1, \ldots, N}$ satisfies 
\begin{align}
\label{eq:def-Wiid}
W_N=\begin{pmatrix}
R_N^{1/2}&\\
&\ddots\\
&&R_N^{1/2}
\end{pmatrix}W_{iid,N}
\end{align}
As the spectral norm of $R_N$ is assumed uniformly bounded when $N$ increases, 
the statement of the lemma is equivalent to $\sup_{N} \ex\{\|W_{iid}\|^{2l}\} < +\infty$. It is shown in 
 \cite{L:15} that the empirical eigenvalue distribution of $W_{iid,N}W_{iid,N}^{*}$ converges towards 
the Marcenko-Pastur distribution, and that its smallest non zero eigenvalue and its largest eigenvalue (which coincides with $\| W_{iid,N} \|^{2}$)  converge
almost surely towards $(1 - \sqrt{c_*})^{2}$ and $(1 + \sqrt{c_*})^{2}$ respectively.  We express $\ex\{\|W_{iid}\|^{2l}\}$ as 
\begin{align*}
&\ex\{\|W_{iid}\|^{2l}\}=\ex\{\|W_{iid}\|^{2l}\mathbf{1}_{\|W_{iid}\|^{2}\le(1+\sqrt{c_*})^2+\delta}\}+\ex\{\|W_{iid}\|^{2l}\mathbf{1}_{\|W_{iid}\|^{2}>(1+\sqrt{c_*})^2+\delta}\}\\
&\le \kappa+\ex\{\|W_{iid}\|_F^{2l}\mathbf{1}_{\|W_{iid}\|^{2}>(1+\sqrt{c_*})^2+\delta}\}\le \kappa+\ex\{\|W_{iid}\|^{4l}_F\}^{1/2}\ex\{\mathbf{1}_{\|W_{iid}\|^{2}>(1+\sqrt{c_*})^2+\delta}\}^{1/2}
\end{align*}
where $\delta > 0$ is a nice constant. As $\ex\{\|W_{i.i.d.}\|^{4l}_F\} = \mathcal{O}(N^{2l})$, it is sufficient to prove that $\ex\{\mathbf{1}_{\|W_{iid}\|^{2}>(1+\sqrt{c_*})^2+\delta}\}$ is less than any power of $N^{-1}$.
We introduce a smooth function $\phi_0$ defined on $\mathbb{R}$ by
\begin{align*}
\phi_0(\lambda)=\begin{cases}
1,\, \text{for } \lambda\in[-\infty,\,-\delta]\cup[(1+\sqrt{c_*})^2+\delta,\,+\infty],\\
0,\, \text{for } \lambda\in[-\delta/2,\, (1+\sqrt{c_*})^2+\delta/2]
\end{cases}
\end{align*}
and $\phi_0(\lambda)\in(0,\,1)$ elsewhere. Then, it holds that 
\begin{align*}
&\ex\{\mathbf{1}_{\|W_{iid}\|^2>(1+\sqrt{c_*})^2+\delta}\}=\ex\{\mathbf{1}_{\lambda_{max}(W_{iid}W_{iid}^*)>(1+\sqrt{c_*})^2+\delta}\}\le \mathbf{P}[\tr\phi_0(W_{iid}W_{iid}^*)\ge1]\\
&\le \ex\{\tr\phi_0(W_{iid}W_{iid}^*)^{2k}\}
\end{align*}
for any $k\in\mathbb{N}$. Lemma \ref{moment} thus appears as an immediate consequence of the following lemma. 

\begin{lemma}
\label{lem:eval-trace-phi}
For each smooth function $\phi$ such that $\phi(\lambda) = 0$ 
if $\lambda \in  [-\delta/2,\, (1+\sqrt{c_*})^2+\delta/2]$ and $\phi(\lambda)$  
constant on $[-\infty,\,-\delta]\cup[(1+\sqrt{c_*})^2+\delta,\,+\infty]$, it holds that 
 $\forall k\in\mathbb{N}$, $\ex\left\{\left(\tr\phi (W_{iid}W^*_{iid})\right)^{2k}\right\}\le\dfrac{\kappa}{N^{2k}}$.
\end{lemma}
\textbf{Proof.} We prove the Lemma by induction.  We first consider the case $k=1$. 
For more convenience we will write $W$ instead of $W_{iid}$ in the course of the proof. Here and below we take sum for all possible values of indexes, if not specified.  From (\ref{p-n}) we have
\begin{align}\label{var_phi}
&\var\{\tr\phi(WW^*)\}
\le \sum \ex\left\{\left(\dfrac{\partial\tr\phi(WW^*)}{\partial \overline{W}^{m_1}_{i_1,j_1}}\right)^*\ex\{W^{m_1}_{i_1,j_1}\overline{W}^{m_2}_{i_2,j_2}\}\dfrac{\partial\tr\phi(WW^*)}{\partial \overline{W}^{m_2}_{i_2,j_2}}\right\}\notag\\
&+\sum \ex\left\{\dfrac{\partial\tr\phi(WW^*)}{\partial W^{m_1}_{i_1,j_1}}\ex\{W^{m_1}_{i_1,j_1}\overline{W}^{m_2}_{i_2,j_2}\}\left(\dfrac{\partial\tr\phi(WW^*)}{\partial W^{m_2}_{i_2,j_2}}\right)^*\right\}
\end{align}
We only evaluate the first term,  denoted by $\psi$, of the right handside of (\ref{var_phi}), because the second one can be addressed similarly. For this, we first remark that 
\begin{align*}
\dfrac{\partial\tr \phi(WW^*)}{\partial \overline{W}^{m_1}_{i_1,j_1}}
=\tr \left(\phi^\prime(WW^*)\dfrac{\partial WW^*}{\partial \overline{W}^{m_1}_{i_1,j_1}}\right)
= \left(\phi^\prime(WW^*)W\right)^{m_1}_{i_1,j_1}.
\end{align*}
Plugging this into (\ref{var_phi}) we obtain
\begin{align*}
\psi=\sum\dfrac{1}{N}\ex\left\{\left(\phi^\prime(WW^*)W\right)^{*m_1}_{j_1,i_1}
\delta_{m_1,m_2}\delta_{i_1+j_1,i_2+j_2}
\left(\phi^\prime(WW^*)W\right)^{m_2}_{i_2,j_2}\right\}.
\end{align*}
Denoting $l=i_1-i_2$, it is easy to verify that $\psi$ can be written as 
\begin{align}\label{eq:psy}
\psi=\dfrac{1}{N}\sum_{l=-(L-1)}^{L-1}\ex\{\tr\left(\phi^\prime(WW^*)W\right)^*(J^l_L\otimes I_M)\left(\phi^\prime(WW^*)W\right)J^l_N\}.
\end{align}
where we recall that matrix $J_L$ is defined by (\ref{eq:def-J}). For each $ML \times N$ matrices $A$ and $B$, the Schwartz inequality and the inequality between arithmetic and geometric means lead to
\begin{align*}
\left| \dfrac{1}{N}\tr A^{*}(J^{*u}_{L}\otimes I_M)BJ^{*u}_N \right| \le \dfrac{1}{2N}\tr A^*(J^{*u}_LJ^{u}_L\otimes I_M)A+\dfrac{1}{2N}\tr B^*J^{*u}_NJ^u_NB.
\end{align*}
Therefore, since $J^{*u}_LJ^u_L\otimes I_M\le I_{ML}$ and $J^{*u}_NJ^u_N\le I_N$
\begin{align}\label{AB}
\left| \dfrac{1}{ML}\tr A^{*}(J^{*u}_{L}\otimes I_M)BJ^{*u}_N \right| \le\dfrac{\kappa}{N}(\tr A^*A+\tr B^*B).
\end{align}
By taking here $A=B=\phi^\prime(WW^*)W$, we obtain from (\ref{var_phi}) and (\ref{eq:psy}) 
\begin{align}\label{eq:var-phy-2-bis}
\var\{\tr\phi(WW^*)\}\le\dfrac{\kappa}{N}\ex\left\{\tr\left(\phi^\prime(WW^*)\right)^2WW^*\right\}.
\end{align}
Consider the function $\eta(\lambda)=(\phi^\prime(\lambda))^2\lambda$. It is clear that $\eta(\lambda)$ is a compactly supported
smooth function. Therefore (see e.g. \cite{L:15}), it holds that
\begin{align*}
\mathbb{E}\left\{\dfrac{1}{ML}\tr\left((\phi^\prime(WW^*))^2WW^*\right)\right\}=\int_{\mathcal{S}_N}\eta(\lambda)d\mu_{MP,N}(\lambda)+\mathcal{O}\left(\dfrac{1}{N^2}\right),
\end{align*}
where $\mu_{MP,N}$ is the measure associated to Marcenko-Pastur distribution with parameters (1, $c_N$) and where
$\mathcal{S}_{MP,N} \subset [0, (1 + \sqrt{c_N})^{2}]$ represents the support of $\mu_{MP,N}$. 
It is clear that for $N$ large enough, the support of $\phi^\prime$ and $\mathcal{S}_{MP,N}$ do not intersect, so that $\int_{\mathcal{S}_N}\eta(\lambda)d\mu_{MP,N}(\lambda)=0$. Therefore, we obtain that
\begin{align*}
\mathbb{E}\left\{\dfrac{1}{ML}\tr\left((\phi^\prime(WW^*))^2WW^*\right)\right\}=\mathcal{O}\left(\dfrac{1}{N^2}\right).
\end{align*}
This and (\ref{eq:var-phy-2-bis}) lead to the conclusion that $\var\{\tr\phi(WW^*)\}=\mathcal{O}\left(N^{-2}\right)$. To finalize
the case $k=1$, we express $\mathbb{E}\{(\tr\phi(WW^*))^2\}$ as $\mathbb{E}\{(\tr\phi(WW^*))^2\}=\var\{\tr\phi(WW^*)\}+\mathbb{E}\{\tr\phi(WW^*)\}^2$. \cite[Lemma~10.1]{L:15} implies that $\mathbb{E}\{\tr\phi(WW^*)\}=\mathcal{O}(N^{-1})$, which completes the proof for $k=1$.

Now we suppose that for any $n\le k$ we have $\mathbb{E}\{(\tr\phi(WW^*))^{2n}\}=\mathcal{O}(N^{-2n})$ and are about to prove that it holds for $n=k+1$.
As in the previous case we write
\begin{align}\label{eq:induc_pass}
\mathbb{E}\{(\tr\phi(WW^*))^{2(k+1)}\}=\var\{(\tr\phi(WW^*))^{k+1}\}+\Big(\mathbb{E}\{(\tr\phi(WW^*))^{k+1}\}\Big)^2
\end{align}
To evaluate the second term of the r.h.s. of (\ref{eq:induc_pass}), we use the Schwartz inequality and the induction assumption
\begin{align}\label{ineq:tr_phi}
\mathbb{E}\{(\tr\phi(WW^*))^{k+1}\}\le \Big(\mathbb{E}\{(\tr\phi(WW^*))^{2k}\}\mathbb{E}\{(\tr\phi(WW^*))^2\}\Big)^{1/2}=\mathcal{O}\left(\dfrac{1}{N^{k+1}}\right)
\end{align}  
We follows the same steps as in the case $k=1$ to study the first term of the r.h.s. of (\ref{eq:induc_pass}). Using again the
Poincaré-Nash inequality, we obtain that 
 \begin{align}
 \var\{(\tr\phi(WW^*))^{k+1}\}
 \le \dfrac{\kappa}{N}\ex\left\{\left(\tr\phi(WW^*)\right)^{2k}\tr\left(\phi^\prime(WW^*)^2WW^*\right)\right\}.
 \end{align}
 Using Holder's inequality, we obtain 
 \begin{align}
   \label{eq:intermediate-inequality}
\var\{(\tr\phi(WW^*))^{k+1}\}
\le \dfrac{\kappa}{N}\ex\left\{\left(\tr\phi(WW^*)\right)^{2k+2}\right\}^{\frac{k}{k+1}}\ex\left\{\left(\tr(\phi^\prime(WW^*)^2WW^*)\right)^{k+1}\right\}^{\frac{1}{k+1}}.
\end{align} 
The properties of function $\eta(\lambda)=\phi^\prime(\lambda)^2\lambda$ imply that it satisfies the induction hypothesis and that it verifies (\ref{ineq:tr_phi}), i.e.  $\ex\left\{(\tr(\phi^\prime(WW^*)^2WW^*))^{k+1}\right\}=\mathcal{O}(N^{k+1})$. Plugging this into (\ref{eq:intermediate-inequality}), we get that
\begin{align}
\var\{(\tr\phi(WW^*))^{k+1}\}
\le \dfrac{\kappa}{N^2}\ex\left\{\left(\tr\phi(WW^*)\right)^{2k+2}\right\}^{\frac{k}{k+1}}.
\end{align} 
From this, (\ref{ineq:tr_phi}) and (\ref{eq:induc_pass}), we immediately obtain
\begin{equation}
  \label{eq:inequality-trace-phi-2k+2}
\mathbb{E}\{(\tr\phi(WW^*))^{2k+2}\}\le\frac{\kappa_1}{N^2}\ex\{(\tr\phi(WW^*))^{2k+2}\}^{\frac{k}{k+1}}+\frac{\kappa_2}{N^{2k+2}}
\end{equation}
We denote by $z_{k,N}$ the term $z_{k,N}=N^{2k+2} \, \mathbb{E}\{(\tr\phi(WW^*))^{2k+2}\}$. Then, (\ref{eq:inequality-trace-phi-2k+2}) implies that
$$
z_{k,N} \leq \kappa_1 \, (z_{k,N})^{k/(k+1)} + \kappa_2
$$
This inequality leads to the conclusion that sequence $(z_{k,N})_{N \geq 1}$ is bounded, or equivalently that \\
$\mathbb{E}\{(\tr\phi(WW^*))^{2k+2}\} \leq \frac{\kappa}{N^{2k+2}}$ as expected. This completes the proof of Lemmas \ref{lem:eval-trace-phi}
and \ref{moment}. $\blacksquare$ \\

We now evaluate the variance of useful functionals of the resolvent $\mathbf{Q}_N(z)$.
\begin{lemma}\label{var}
	Let $(F_N)_{N \geq 1}$, $(G_N)_{N \geq 1}$ be sequences of deterministic $2ML \times 2ML$ matrices and 
$(H_N)_{N \geq 1}$ a sequence of deterministic $N\times N$ matrices such that $\max\{\sup_N \|F_N\|,\,\sup_N\|G_N\|,\,\sup_N\|H_N\|\} \le \kappa$, and consider sequences of deterministic $2 ML$--dimensional vectors $(a_{1,N})_{N \geq 1}$, $(a_{2,N})_{N \geq 1}$ such that $sup_N \|a_{i,N}\| \le\kappa$ for $i=1,2$. Then, for each $z \in \mathbb{C}^+$, it holds that
	\begin{align}\label{var_1}
	&\var\left\{\dfrac{1}{ML}\tr F\mathbf{Q}\right\}\le\dfrac{C(z)\kappa^2}{N^{2}},\\
	&\var\left\{\dfrac{1}{ML}\tr F\mathbf{Q}GWHW^*\right\}\le\dfrac{C(z)\kappa^6}{N^{2}}.\label{var_2}\\
	&\var\left\{a_1^*\mathbf{Q}a_2\right\}\le \dfrac{C(z)\kappa^4}{N}\label{var_3}
	\end{align}
	where $C(z)$ can be written as $C(z) = P_1(|z|)P_2\left(\frac{1}{\im z}\right)$
	for some nice polynomials $P_1$ and $P_2$.
\end{lemma}
\textbf{Proof.} We first prove (\ref{var_1}) and denote by $\xi$ the term $\xi=\frac{1}{ML}\tr F\mathbf{Q}$. The Poincare-Nash inequality leads to
\begin{align*}
&\var\{\xi\}\le \sum_{\substack{i_1,j_1,m_1\\i_2,j_2,m_2}}\ex\left\{\left(\dfrac{\partial\xi}{\partial \overline{W}_{i_1,j_1}^{m_1}}\right)^*\ex\{W_{i_1,j_1}^{m_1}\overline{W}_{i_2,j_2}^{m_2}\}\dfrac{\partial\xi}{\partial\overline{W}_{i_2,j_2}^{m_2}}\right\}\\
&+\sum_{\substack{i_1,j_1,m_1\\i_2,j_2,m_2}}\ex\left\{\dfrac{\partial\xi}{\partial W_{i_1,j_1}^{m_1}}\ex\{W_{i_1,j_1}^{m_1}\overline{W}_{i_2,j_2}^{m_2}\}\left(\dfrac{\partial\xi}{\partial W_{i_2,j_2}^{m_2}}\right)^*\right\}.
\end{align*}
We just evaluate the first term of r.h.s., denoted by $\phi$.
For this, we  need the expression of the  derivative of $\mathbf{Q}$ with respect to the complex conjugates of the entries of $W$.
We denote by $\Pi_{pf}$ and $\Pi_{fp}$ as $2ML\times 2ML$ matrices defined by  $\Pi_{pf}=\left(\begin{smallmatrix}
0&I_{ML}\\
0&0
\end{smallmatrix}\right)$ and $\Pi_{fp}=\left(\begin{smallmatrix}
0&0\\
I_{ML}&0
\end{smallmatrix}\right)$. Then, after some algebra, we obtain that
\begin{align}
  \label{deriv}
  &\dfrac{\partial \mathbf{Q}}{\partial\overline{W}_{i,j}^{m}} = - {\bf Q} \left( \begin{smallmatrix} w_{j,f} \\ 0 \end{smallmatrix} \right)
  ({\bf f}_{i+L}^{m})^{T} {\bf Q} \, \mathbf{1}_{i \le L} - {\bf Q} \left( \begin{smallmatrix} 0 \\ w_{j,p}  \end{smallmatrix} \right)
  ({\bf f}_{i-L}^{m})^{T} {\bf Q} \, \mathbf{1}_{i > L} \notag\\
  &= - {\bf Q} \Pi_{pf} W {\bf e}_j \, ({\bf f}_i^{m})^{T} \Pi_{pf} {\bf Q}  - {\bf Q} \Pi_{fp} W {\bf e}_j \, ({\bf f}_i^{m})^{T} \Pi_{fp} {\bf Q}
\end{align}
From this, we deduce immediately that 
$$
\dfrac{\partial\xi}{\partial \overline{W}_{i_1,j_1}^{m_1}}= -\frac{1}{ML} \, \Big(\Pi_{pf}\mathbf{Q}F\mathbf{Q}\Pi_{pf}W+\Pi_{fp}\mathbf{Q}F\mathbf{Q}\Pi_{fp}W\Big)_{i_1,j_1}^{m_1}
$$
Using that $\ex\{W_{i_1,j_1}^{m_1}\overline{W}_{i_2,j_2}^{m_2}\}=\frac{1}{N}R_{m_1m_2}\delta_{i_1+j_1,i_2+j_2}$, we obtain that $\phi$ is given by
\begin{align*}
&\phi
=\frac{1}{N (ML)^{2}} \sum_{\substack{i_1,j_1,m_1\\i_2,j_2,m_2}}(\mathbf{e}_{j_1})^T(\Pi_{pf}\mathbf{Q}F\mathbf{Q}\Pi_{pf}W+\Pi_{fp}\mathbf{Q}F\mathbf{Q}\Pi_{fp}W)^*\mathbf{f}_{i_1}^{m_1}R_{m_1m_2}\\
&\times\delta_{i_1+j_1,i_2+j_2}(\mathbf{f}_{i_2}^{m_2})^T
(\Pi_{pf}\mathbf{Q}F\mathbf{Q}\Pi_{pf}W+\Pi_{fp}\mathbf{Q}F\mathbf{Q}\Pi_{fp}W)\mathbf{e}_{j_2}
\end{align*}
We put $u=i_1-i_2$ and remark that $\sum_{m_1,m_2,i_1-i_2=u}\mathbf{f}_{i_1}^{m_1}R_{m_1m_2}(\mathbf{f}_{i_2}^{m_2})^T= J^{*u}_L\otimes R$ and that $\sum_{j_2-j_1=u}\mathbf{e}_{j_2}\mathbf{e}_{j_1}^T=J_{N}^{*u}$. Therefore, $\phi$ can be written as
\begin{align} \label{phi}
&\phi=\dfrac{1}{MLN}\ex\Big\{\sum_{u=-(L-1)}^{L-1}\dfrac{1}{ML}\tr(\Pi_{pf}\mathbf{Q}F\mathbf{Q}\Pi_{pf}W+\Pi_{fp}\mathbf{Q}F\mathbf{Q}\Pi_{fp}W)^*(J^{*u}_L\otimes R)\notag\\
&\times(\Pi_{pf}\mathbf{Q}F\mathbf{Q}\Pi_{pf}W+\Pi_{fp}\mathbf{Q}F\mathbf{Q}\Pi_{fp}W)J_N^{*u}\Big\}
\end{align}
Each term inside the sum over $u$ can be written as $\dfrac{1}{ML}\tr A^{*}(I_L \otimes R^{1/2}) (J^{*u}_{L}\otimes I)(I_L \otimes R^{1/2}) AJ^{*u}_N$, where the expression of
$ML \times N$ matrix $A$ is omitted. 
As $\| R \|$ is bounded by the nice constant $b$ (see (\ref{eq:hypothesis-R-bis})), (\ref{AB}) and (\ref{phi}) lead to the conclusion
that we just need to evaluate $\frac{1}{ML} \ex\{\tr A^*A\}$. Using the Schwartz inequality, we obtain immediately
that
\begin{align}
  \label{eq:inequality-trA*A}
&\ex\{\tr A^*A\} \leq 2 \ex\{\tr \left((\Pi_{pf}\mathbf{Q}F\mathbf{Q}\Pi_{pf}W)^{*}\Pi_{pf}\mathbf{Q}F\mathbf{Q}\Pi_{pf}W\right)\}\\
&\hskip5cm+ 2 \ex\{\tr \left((\Pi_{fp}\mathbf{Q}F\mathbf{Q}\Pi_{fp}W)^{*}\Pi_{fp}\mathbf{Q}F\mathbf{Q}\Pi_{fp}W\right)\}\notag
\end{align}
 Since $\left(\Pi_{pf}\mathbf{Q}F\mathbf{Q}\Pi_{pf}\right)^{*} \Pi_{pf}\mathbf{Q}F\mathbf{Q}\Pi_{pf} \leq \| {\bf Q}\|^{4} \| F \|^{2} \, I$
 and $\|\mathbf{Q}\|\le \frac{1}{\im z}$, we get that
 \begin{align*}
\frac{1}{ML}  \ex\{\tr \left((\Pi_{pf}\mathbf{Q}F\mathbf{Q}\Pi_{pf}W)^{*}\Pi_{pf}\mathbf{Q}F\mathbf{Q}\Pi_{pf}W\right)\} &\leq \frac{1}{(\im z)^{4}} \, \|F \|^{2} \,
 \frac{1}{ML} \ex\{\tr W^*W \}\\
 &\leq \frac{1}{(\im z)^{4}} \, \|F \|^{2} \, \mathbb{E}(\| W \|^{2})
  \end{align*}
 Lemma \ref{moment} thus implies that
 $$
 \frac{1}{ML} \ex\{\tr \left((\Pi_{pf}\mathbf{Q}F\mathbf{Q}\Pi_{pf}W)^{*}\Pi_{pf}\mathbf{Q}F\mathbf{Q}\Pi_{pf}W\right)\} \leq \kappa^{2} P\left(\frac{1}{\im z}\right) 
   $$
\mbox{for some nice polynomial $P$. The term $ \frac{1}{ML} \ex\{\tr \left(\Pi_{fp}\mathbf{Q}F\mathbf{Q}\Pi_{fp}W)^{*}\Pi_{fp}\mathbf{Q}F\mathbf{Q}\Pi_{fp}W\right)\}$}  can be handled similarly. Therefore, (\ref{phi}) leads to $\phi \leq \kappa^{2} \, \frac{1}{N^{2}}  P\left(\frac{1}{\im z}\right)$.  
 This establishes  (\ref{var_1}). 
  
 To prove (\ref{var_2}) one can also use  Poincaré-Nash inequality for $\xi=\frac{1}{ML}\tr F\mathbf{Q}GWHW^*$. After some calculations, we get that
 the variance of $\xi$ is upperbounded by a term given by
 \begin{equation}
   \label{eq:poincare-xi}
   \frac{\kappa_1}{N^{2}} \, \left(\frac{1}{ML}\tr (F\mathbf{Q}GWH)^*(F\mathbf{Q}GWH) + \frac{1}{ML}\tr (F\mathbf{Q}WH)^*(F\mathbf{Q}WH) + \eta_1
   + \eta_2 \right)
   \end{equation}
 where $\kappa_1$ is some nice constant, and where $\eta_1$ and $\eta_2$ are defined by
 \begin{eqnarray}
   \label{eq:def-eta1}
   \eta_1 & = & \frac{1}{ML}\tr (\Pi_{pf}\mathbf{Q}GWHW^*F\mathbf{Q}\Pi_{pf}W)^*(\Pi_{pf}\mathbf{Q}GWHW^*F\mathbf{Q}\Pi_{pf}W) \\
   \eta_2 & = & \frac{1}{ML}\tr (\Pi_{fp}\mathbf{Q}GWHW^*F\mathbf{Q}\Pi_{fp}W)^*(\Pi_{fp}\mathbf{Q}GWHW^*F\mathbf{Q}\Pi_{fp}W)
 \end{eqnarray}
 Using Lemma \ref{moment} as well as the inequality $ \mathbf{Q} \mathbf{Q}^{*} \leq \frac{1}{\im^2z} I$, we obtain
 immediately (\ref{var_2}). \\

(\ref{var_3}) is the consequence of (\ref{var_1}) since 
 $a_1^*\mathbf{Q}a_2=\tr \mathbf{Q}a_2a_1^* = \tr \mathbf{Q} F$ for $F=a_2a_1^*$. This completes the proof of Lemma \ref{var}. $\blacksquare$ \\

 In the following, we also need to evaluate the variance of more specific terms. The following result appears to be a
 consequence of Lemma \ref{var} and of the particular structure (\ref{eq:expre-BQ}) of matrix $\mathbf{Q}(z)$.
 \begin{corollary}
   \label{eq:consequences-var}
   Let $(F_{1,N})_{N \geq 1}$ be a sequence of deterministic $ML \times ML$ matrices such that $\sup_{N} \| F_{1,N} \| \leq \kappa$,
   and $(H_N)_{N \geq 1}$ a sequence of deterministic $N \times N$ matrices satisfying $\sup_{N} \| H_N \| \leq 1$. Then, if $\im z^{2} > 0$,
   the following evaluations hold:
   \begin{equation}
     \label{eq:var-Qij}
    \var \left\{ \frac{1}{ML} \mathrm{Tr} F_1 {\bf Q}_{ij}(z) \right\} \leq \kappa^{2} \frac{1}{N^{2}} P_1(|z^{2}|) P_2(\frac{1}{\im z^{2}})
   \end{equation}
   where $i$ and $j$ belong to $\{ p,f \}$;
   \begin{equation}
     \label{eq:var-QijWW*}
     \var \left\{ \frac{1}{ML} \mathrm{Tr} \left[ H W^{*} \Pi_{i_1j_1} \left( \begin{array}{cc} F_1 & 0 \\ 0 & 0 \end{array} \right) \mathbf{Q}(z) \Pi_{i_2j_2} W
       \right] \right \}
     \leq \kappa^{2} \frac{1}{N^{2}} P_1(|z^{2}|) P_2(\frac{1}{\im z^{2}})
   \end{equation}
   where $i_1,j_1,i_2,j_2$ still belong to  $\{ p,f \}$, but verify $i_1 \neq j_1$ and $i_2 \neq j_2$. 
 \end{corollary}
{\bf Proof.}  We first prove (\ref{eq:var-Qij}), and first consider the case where $i=j=p$. We 
define the $2ML \times 2ML$ matrix $F$ by 
for $F =  \left( \begin{array}{cc} F_1 & 0 \\ 0 & 0 \end{array} \right)$, and remark that 
$\frac{1}{ML} \mathrm{Tr} F_1 {\bf Q}_{pp}(z)$ coincides with $\xi = \frac{1}{ML} \mathrm{Tr} F {\bf Q}(z)$. 
We follow the proof of (\ref{var_1}), and evaluate the right hand side of (\ref{eq:inequality-trA*A}) in a more 
accurate manner by taking into account the particular structure of the present matrix $F$. It is easy to check that
\begin{align*}
&\frac{1}{ML}  \ex\{\tr \left(\Pi_{pf}\mathbf{Q}F\mathbf{Q}\Pi_{pf}W)^{*}\Pi_{pf}\mathbf{Q}F\mathbf{Q}\Pi_{pf}W\right)\} \\
&\hskip5cm= \frac{1}{ML}  \ex\{\tr \left( W_f^{*} {\bf Q}_{pp}^{*}F_1^{*}{\bf Q}_{fp}^{*}{\bf Q}_{fp} F_1 {\bf Q}_{pp}  W_f
\right)\}
\end{align*}
As ${\bf Q}_{fp}(z) = W_p W_f^{*} Q(z^{2})$, we obtain that 
$$
{\bf Q}_{fp}^{*}(z) {\bf Q}_{fp}(z) = (Q(z^{2}))^{*} W_f W_p^{*} W_p W_f^{*} Q(z^{2}) \leq \| W \|^{4} 
\frac{1}{(\im z^{2})^{2}} \, I
$$
if $\im (z^{2}) > 0$. Therefore, it holds that 
$$
F_1^{*}{\bf Q}_{fp}^{*}{\bf Q}_{fp} F_1  \leq  \kappa^{2} \| W \|^{4} 
\frac{1}{(\im z^{2})^{2}} \, I
$$
From this, using the expression of ${\bf Q}_{pp}=zQ(z^2)$, we obtain similarly that
$$
 W_f^{*} {\bf Q}_{pp}^{*}F_1^{*}{\bf Q}_{fp}^{*}{\bf Q}_{fp} F_1 {\bf Q}_{pp} W_f 
\leq \kappa^{2} \| W \|^{6} \frac{|z|^2}{(\im z^{2})^{4}}
$$
Lemma \ref{moment} thus leads to the conclusion that 
$$
 \frac{1}{ML}  \ex\{\tr \left( W_f^{*} {\bf Q}_{pp}^{*}F_1^{*}{\bf Q}_{fp}^{*}{\bf Q}_{fp} F_1 {\bf Q}_{pp} { W}_f
\right)\} \leq \kappa^{2} \frac{\kappa_1|z|^2}{(\im z^{2})^{4}} 
$$
where $\kappa_1$ is a nice constant such that $\mathbb{E}( \| W_N \|^{6}) \leq \kappa_1$ for each $N$. Using similar 
arguments, we obtain that 
$$
 \frac{1}{ML} \ex\{\tr \left(\Pi_{fp}\mathbf{Q}F\mathbf{Q}\Pi_{fp}W)^{*}\Pi_{fp}\mathbf{Q}F\mathbf{Q}\Pi_{fp}W\right)\} \leq  \kappa^{2} \frac{\kappa_1 |z^{2}|^{2}}{(\im z^{2})^{4}}
$$ 
This, in turn, implies (\ref{eq:var-Qij}) for $i=j=p$. As the arguments are essentially the same for 
the other values of $i$ and $j$, we do not provide the corresponding proofs. \\

In order to establish (\ref{eq:var-QijWW*}), we follow the proof (\ref{var_2}) for 
$F = \Pi_{i1j_1}  \left( \begin{array}{cc} F_1 & 0 \\ 0 & 0 \end{array} \right)$, 
$G = \Pi_{i_2j_2}$. It is necessary to check that the 4 terms inside the bracket of 
(\ref{eq:poincare-xi}) can be upperbounded by $\kappa^{2} P_1(|z^{2}|) P_2(\frac{1}{\im z^{2}})$ 
for nice polynomials $P_1$ and $P_2$. As above, the use of the particular expression of 
matrices $(Q_{i,j})_{i,j \in \{f,p \}}$ allows to establish this property. The corresponding 
easy calculations are omitted.  $\blacksquare$

\section{Various lemmas on Stieltjes transform}
\label{sec:stieljes-lemmas}

In this paragraph, we provide a number of useful results on certain Stieltjès transforms. In the following, if $A$ is a Borel set of $\mathbb{R}$, we denote by $\mathcal{S}_M(A)$ the set of all Stieltjes transforms of 
$M \times M$ matrix valued positive finite measures carried by $A$. $\mathcal{S}_1(A)$ is denoted  $\mathcal{S}(A)$. 
We first begin by stating well known properties of Stieltjès transforms.  
\begin{proposition}\label{stil_tran}
	The following properties hold true:
	
	1. Let $f$ be the Stieltjes transform of a positive finite measure $\mu$, then
	
	– the function $f$ is analytic over $\mathbb{C}^+$,
	
	– if $z\in \mathbb{C}^+$ then $f(z) \in \mathbb{C}^+$,
	
	– the function $f$ satisfies: $|f(z)|\le \frac{\mu(\mathbb{R})}{\im z}$, for $z\in\mathbb{C}^+$
	
	– if $\mu(-\infty, 0) = 0$ then its Stieltjes transform $f$ is analytic over $\mathbb{C}/\mathbb{R}^+$. Moreover, $z \in \mathbb{C}^+$ implies $zf(z) \in \mathbb{C}^+$.
	
	– for all $\phi\in\mathcal{C}^{\infty}_c(\mathbb{R},\mathbb{R})$ we have
	\begin{align}
	\int_{\mathbb{R}}\phi(\lambda)d\mu(\lambda)=\dfrac{1}{\pi}\lim\limits_{y\downarrow0}\im \left\{\int_{\mathbb{R}}\phi(x)f(x+iy)dx\right\}
	\end{align}
	
	2. Conversely, let $f$ be a function analytic over $\mathbb{C}^+$ such that $f(z) \in\mathbb{C}^+$ if $z \in \mathbb{C}^+$ and for which $\sup_{y \geq \epsilon} |iy f(iy)| < +\infty$ for some $\epsilon > 0$. Then, $f$ is the Stieltjès transform 
of a unique positive finite measure $\mu$ such that $\mu(\mathbb{R}) = \lim_{y \rightarrow +\infty} -iy f(iy)$. 
Moreover, the following inversion formula holds:
	\begin{align}
	\mu([a,b])=\lim\limits_{\nu\rightarrow0+}\dfrac{1}{\pi}\int_{a}^{b}\im  f(\xi+i\nu)d\xi,
	\end{align}
	whenever $a$ and $b$ are continuity points of $\mu$. If moreover $zf(z) \in\mathbb{C}^+$ for $z$	in $\mathbb{C}^+$ then, $\mu(\mathbb{R}^-) = 0$. In particular, $f$ is given by
	\begin{align*}
f(z)=\int_{0}^{+\infty}\dfrac{\mu(d\lambda)}{\lambda-z}
	\end{align*}
	and has an analytic continuation on $\mathbb{C}/\mathbb{R}^+$.
	
	3. Let $F$ be an $P \times P$ matrix-valued function analytic on $\mathbb{C}^+$ verifying

	– $\im (F(z)) > 0$ if $z \in \mathbb{C}^{+}$

	– $\sup_{y > \epsilon} \| iy F(iy) \| < +\infty$ for some $\epsilon > 0$.
	
	Then, $F \in \mathcal{S}_P(\mathbb{R})$, and if $\mu^{F}$ is the corresponding $P \times P$ associated positive
        measure, it holds that
        \begin{equation}
          \label{eq:mass-omega}
          \mu^{F}(\mathbb{R}) = \lim_{y \rightarrow +\infty} -iy F(iy)
          \end{equation}
          If moreover  $\im  (zF(z)) > 0$, then, $F \in 
	\mathcal{S}_P(\mathbb{R}^{+})$.

%
\end{proposition}

We now state a quite useful Lemma.

\begin{lemma}
\label{beta}
Let $\beta(z)\in\mathcal{S}(\mathbb{R}^+)$, and consider function $\bm{\beta}(z)$ defined by $\bm{\beta}(z)=z\beta(z^2)$. Then $\bm{\beta}\in \mathcal{S}(\mathbb{R})$. Moroever,
it holds that
\begin{align}\label{stieltjes_bf}
&\mathbf{G}(z)=\left(-zI_M-\dfrac{c\bm{\beta}(z)}{1-c^2\bm{\beta}^2(z)}R\right)^{-1}\in \mathcal{S}_M(\mathbb{R})\\
&G(z)=\left(-zI_M-\dfrac{cz\beta(z)}{1-zc^2\beta^2(z)}R\right)^{-1}\in \mathcal{S}_M(\mathbb{R}^+).\label{stieltjes_nonbf}
\end{align}
and that 
\begin{equation}
  \label{eq:bound-GG*}
        {\bf G}(z) \left({\bf G}(z)\right)^{*} \leq \frac{I_M}{(\im z)^{2}}, \;  G(z) \left(G(z) \right)^{*} \leq \frac{I_M}{(\im z)^{2}}
        \end{equation}
Finally, matrices ${\bf G}(z)$ and $G(z)$ are linked by the relation 
\begin{equation}
\label{eq:relation-bfG-G}
{\bf G}(z) = z G(z^2)
\end{equation}
for each $z \in \mathbb{C}^{+}$. 
\end{lemma}
\textbf{Proof.} Let $\tau$ be the measure carried by $\mathbb{R}^{+}$ corresponding to the Stieltjes transform $\beta(z)$. We first prove that $\bm{\beta}(z)$ is a Stieltjes transform. We first remark that if 
$z \in \mathbb{C}^{+}$, then $z^{2} \in \mathbb{C} - \mathbb{R}^{+}$. $\beta$ analytic on 
$\mathbb{C} - \mathbb{R}^{+}$ thus implies that $\bm{\beta}(z)$ is analytic on $\mathbb{C}^{+}$ 
Moreover, it is clear that 
\begin{align*}
\im \bm{\beta}(z)=\im \int_{\mathbb{R}^+}\dfrac{z d \, \tau(\lambda)}{\lambda-z^2}
=\int_{\mathbb{R}^+}\dfrac{\im  z(\lambda+|z|^2) d \, \tau(\lambda)}{|\lambda-z^2|^2}>0,\, \text{when}\, \im  z>0.
\end{align*}
To evaluate $\bm{\beta}(z)$ for $z \in \mathbb{C}^{+}$, we write
\begin{align*}
\left|\int_{\mathbb{R}^+}\dfrac{z d \, \tau(\lambda)}{\lambda-z^2}\right|
\le\int_{\mathbb{R}^+}\dfrac{d \, \tau(\lambda)}{\left|\frac{\lambda}{z}-z\right|}
\end{align*}
Using that $\left|\frac{\lambda}{z}-z\right| \geq \left| \im  (\frac{\lambda}{z}-z) \right|  \geq \im  z$ for $z \in \mathbb{C}^{+}$ and $\lambda \geq 0$, we get that 
$$
|\bm{\beta}(z)| \le \int_{\mathbb{R}^+}\dfrac{d \, \tau(\lambda)}{\im  z}=\dfrac{ \tau(\mathbb{R}^+)}{\im  z}.
$$
From this and Proposition~\ref{stil_tran}, we obtain that $\bm{\beta}(z)\in\mathcal{S}(\mathbb{R})$. 

To prove (\ref{stieltjes_bf}), it is first necessary to show that $\mathbf{G}$ is analytic on 
$\mathbb{C}^{+}$. For this, we first check that $m(z)=1-c^2\bm{\beta}^2(z)\ne 0$ for $z \in \mathbb{C}^{+}$.
Indeed, write $\bm{\beta}(z)=x+iy$ with $y>0$, then
$m(z)=1-c^{2} x^2+c^{2} y^2-2cxyi$. Hence, if $x=0$ we have $m(z)=1+c^{2}y^2>0$, and if $x\ne0$ then $2xy\ne 0$ since $y>0$. In order to establish that matrix $\left( -zI_M-\dfrac{c\bm{\beta}(z)}{1-c^2\bm{\beta}^2(z)}R \right)$ 
is invertible on $\mathbb{C}^{+}$, we verify that 
\begin{equation}
\label{eq:denominator-beta-negative}
\im  \left( -zI_M-\dfrac{c\bm{\beta}(z)}{1-c^2\bm{\beta}^2(z)}R \right) < 0
\end{equation}
on $\mathbb{C}^{+}$. It is easy to check that
$$
\im  \left( -zI_M-\dfrac{c\bm{\beta}(z)}{1-c^2\bm{\beta}^2(z)}R \right) = - \im  z \, I_M 
- \dfrac{c\im \bm{\beta}(z)(1+c^2|\bm{\beta}(z)|^2)}{|1-c^2\bm{\beta}^2(z)|^2} \, R < - \im  z \, I_M 
$$
Therefore, $\im  z > 0$ and $\im \bm{\beta}(z) > 0$ imply (\ref{eq:denominator-beta-negative}).  
The imaginary part of $\mathbf{G}(z)$ is given by
\begin{equation}
  \label{eq:imboldG>0}
  \im  (\mathbf{G}(z))
  =- \mathbf{G}(z) \im  \left( -zI_M-\dfrac{c\bm{\beta}(z)}{1-c^2\bm{\beta}^2(z)}R \right) \left(\mathbf{G}(z)\right)^* > \im  z \,
   \left(\mathbf{G}(z)  \left( \mathbf{G}(z) \right)^{*}\right) > 0
\end{equation}
Therefore,  $\im  \mathbf{G}(z) > 0$ if $z \in \mathbb{C}^{+}$. We finally remark that 
$\lim_{y \rightarrow +\infty} -iy \mathbf{G}(iy) = I$, which implies that 
$\sup_{y > \epsilon} \| iy {\bf G}(iy) \| < +\infty$ for each $\epsilon > 0$. Proposition 
\ref{stil_tran} eventually implies that $\mathbf{G} \in \mathcal{S}_M(\mathbb{R})$. Moreover, if ${\bs \tau}^{{\bf G}}$ is the underlying
$M \times M$ positive matrix valued measure, (\ref{eq:mass-omega}) leads to ${\bs \tau}^{{\bf G}}(\mathbb{R}) = I$. \\

We prove similarly the analyticity of $G(z)$ on $\mathbb{C}^{+}$. We first check that $1-zc^2\beta^2(z)\ne0$ if $z \in \mathbb{C}^{+}$, or
equivalently that $|1-zc^2\beta^2(z)|\ne0$ if $z \in \mathbb{C}^{+}$. We remark that 
\begin{equation}
  \label{eq:factorizarion-1-z(cbeta)2}
  |1-zc^2\beta^2(z)| = |z \beta(z)| |c^{2} \beta(z) - \frac{1}{z \beta(z)}| > \im  z \, \im  \beta(z) \,  \im  \left(c^{2} \beta(z) - \frac{1}{z \beta(z)} \right)
\end{equation}
As $\beta \in \mathcal{S}(\mathbb{R}^{+})$, it holds that $ \im  \left(c^{2} \beta(z) - \frac{1}{z \beta(z)} \right) > 0$ if
$z \in \mathbb{C}^{+}$. Therefore, $1-zc^2\beta^2(z)\ne0$ if $z \in \mathbb{C}^{+}$. As above, we verify that
\begin{equation}
\label{eq:denominator-beta-negative-bis}
\im  \left( -zI_M-\dfrac{c z \beta(z)}{1-z (c \beta(z))^2} R \right) = - \im  z\,I_M - \im  \left(\dfrac{c z \beta(z)}{1-z (c \beta(z))^2} \right) \, R <
-\im  z \, I_M
\end{equation}
It is easy to check that
$$
\im  \left(\dfrac{c z \beta(z)}{1-z (c \beta(z))^2} \right)  = \frac{c}{|1-z (c \beta(z))^2|^{2}} \, \left( \im  (z \beta(z)) + |zc \beta(z)|^{2} \im \beta(z) \right) > 0
$$
if $z \in \mathbb{C}^{+}$, which, of course, leads to (\ref{eq:denominator-beta-negative-bis}). Therefore, matrix
$\left( -zI_M-\dfrac{c z \beta(z)}{1-z (c \beta(z))^2} R \right)$ is invertible if $z \in \mathbb{C}^{+}$, and $G$ is analytic on $\mathbb{C}^{+}$. Moroever, we obtain
immediately that 
\begin{align}
  \label{eq:imG>0}
&\im (G(z))= G(z) \left(\im  z\,I_M + \im  \left(\dfrac{c z \beta(z)}{1-z (c \beta(z))^2} \right) \, R \right) (G(z))^{*} > \im  z \, \left( G(z) G(z)^{*}\right) > 0\\
&\im(zG(z))=G(z)  \im  \left(\dfrac{c z \beta(z)}{1-z (c \beta(z))^2} \right) \, R  (G(z))^{*} >  0 \notag
\end{align} 
for $z \in \mathbb{C}^{+}$. As above, it holds that $\lim_{y \rightarrow +\infty} -iy G(iy) = I$ and that
$\sup_{y > \epsilon} \|iy G(iy)\| < +\infty$ for each $\epsilon > 0$. This implies that $G \in \mathcal{S}_M(\mathbb{R}^{+})$, and that if
$\tau^{G}$ represents the associated $M \times M$ matrix-valued measure, then $\tau^{G}(\mathbb{R}^{+})= I$. 

In order to establish (\ref{eq:bound-GG*}), we follow \cite[Lemma 3.1]{HT:05}. More precisely, we remark that
$$
\im  G(z) = \im  z \, \int_{\mathbb{R}^{+}} \frac{d \tau^{G}(\lambda)}{|\lambda - z|^{2}} < \frac{\tau^{G}(\mathbb{R}^{+})}{\im  z} = \frac{I}{\im  z}
  $$
  Therefore, (\ref{eq:imG>0}) leads to $ \left( G(z) G(z)^{*}\right) \leq \frac{I}{\left(\im  z\right)^{2}}$. The other statement of
  (\ref{eq:bound-GG*}) is proved similarly and this completes the proof.  $\blacksquare$

  \begin{lemma}
  \label{le:gestion-(1-z(cbeta)2}
  We consider a sequence $(\beta_N)_{N \geq 1}$ of elements of $\mathcal{S}(\mathbb{R}^{+})$ whose associated positive measures
  $(\tau_N)_{N \geq 1}$ satisfy for each $N \geq 1$
  \begin{equation}
    \label{eq:mass-tau}
  \tau_N(\mathbb{R}^{+}) = \frac{1}{M} \tr R_N
  \end{equation}
  as well as
  \begin{equation}
    \label{eq:first-moment-tau}
 \int_{\mathbb{R}^{+}} \lambda \, d \, \tau_N(\lambda) = c_N \, \frac{1}{M} \tr R_N \, \frac{1}{M} \tr R_N^{2}
  \end{equation}
  Then, it exist nice constants $\omega, \kappa$ such that
  \begin{equation}
    \label{eq:lower-bound-imbeta}
    \im  \beta_N(z) \geq \frac{ \kappa \, \im  z}{(\omega^{2} + |z|^{2})}
  \end{equation}
  and
  \begin{equation}
    \label{eq:lowerbound-(1-z(cbeta)2-1}
  \left| 1 - z \left( c_N \beta_N(z) \right)^{2} \right| \geq \frac{ \kappa \, (\im  z)^{3}}{(\omega^{2} + |z|^{2})^{2}}
  \end{equation}
  for each $z \in \mathbb{C}^{+}$ and for each $N \geq 1$. Moreover, if $\bs{\beta}_N(z)$ is defined by
  ${\bs \beta}_N(z) = z \, \beta_N(z^{2})$, then, we also have
 \begin{equation}
    \label{eq:lower-bound-imboldbeta}
    \im  \bs{\beta}_N(z) \geq \frac{ \kappa \, \left(\im  z\right)^{3}}{(\omega^{2} + |z|^{4})}
  \end{equation}
  and
  \begin{equation}
    \label{eq:lowerbound-(1-z(cboldbeta)2-1}
  \left| 1 - \left( c_N \bs{\beta}_N(z) \right)^{2} \right| \geq \frac{ \kappa \, (\im  z)^{6}}{(\omega^{2} + |z|^{4})^{2}}
  \end{equation}
  for each $z \in \mathbb{C}^{+}$ and for each $N \geq 1$.  
\end{lemma}
    {\bf Proof.} We first establish (\ref{eq:lower-bound-imbeta}). $\im  \beta_N(z)$ is given by
    $$
    \im  \beta_N(z) = \im  z \, \int_{\mathbb{R}^{+}} \frac{ d \, \tau_N(\lambda)}{|\lambda - z|^{2}}
      $$
      For each $\omega > 0$, it is clear that
      $$
      \int_{\mathbb{R}^{+}} \frac{ d \, \tau_N(\lambda)}{|\lambda - z|^{2}} \geq  \int_{0}^{\omega} \frac{ d \, \tau_N(\lambda)}{|\lambda - z|^{2}} \geq \frac{\tau_N([0, \omega])}{2( \lambda^{2} + |z|^{2})}
          $$
Assumption (\ref{eq:hypothesis-R-bis}) and (\ref{eq:first-moment-tau}) imply that the sequence $(\tau_N)_{N \geq 1}$ is tight. For each $\epsilon > 0$, it thus exists
$\omega > 0$ for which $\tau_N(]\omega, +\infty[) < \epsilon$ for each $N$, or equivalently,
    $\tau_N([0,\omega]) > \tau_N(\mathbb{R}^{+}) - \epsilon$.
    As $\tau_N(\mathbb{R}^{+}) = \frac{1}{M} \mathrm{Tr}(R_N) > a$, we choose $\epsilon = a/2$,
    and obtain that the corresponding $\omega$ verifies  $\tau_N([0,\omega]) > a/2$ for
    each $N$. This completes the proof of (\ref{eq:lower-bound-imbeta}). We now verify (\ref{eq:lowerbound-(1-z(cbeta)2-1}). For this,
    we use (\ref{eq:factorizarion-1-z(cbeta)2}).
    As $\im  \left(\frac{1}{z \beta_N(z)} \right) < 0$, it holds that $\im  \left(c_N^{2} \beta_N(z) - \frac{1}{z \beta_N(z)} \right) \geq c_N^{2} \im  \beta_N(z)$. Therefore,
    we obtain that
    \begin{equation}
      \label{eq:lowerbound-(1-z(cbeta)2-2}
      \left| 1 - z \left( c_N \beta_N(z) \right)^{2} \right|  \geq c_N^{2} \, \im  z \, \left( \im  \beta_N(z) \right)^{2}
    \end{equation}
    which implies  (\ref{eq:lowerbound-(1-z(cbeta)2-1}). \\

    We finally verify (\ref{eq:lower-bound-imboldbeta}) and (\ref{eq:lowerbound-(1-z(cboldbeta)2-1}).
 For this, we first  express
$\bm{\beta}_N(z)$ as
$$
\bm{\beta}_N(z) = z \beta_N(z^{2}) = \int_{\mathbb{R}^{+}} \frac{z}{\lambda - z^{2}} d \, \tau_N(\lambda)
  $$
which leads immediately to

\begin{align*}
\im  \bm{\beta}_N(z) = \im  z \,  \int_{\mathbb{R}^{+}} \frac{\lambda + |z|^{2}}{|\lambda - z^{2}|^{2}} d \, \tau_N(\lambda) \geq \im  z \, |z|^{2} \int_{\mathbb{R}^{+}} \frac{1}{|\lambda - z^{2}|^{2}} d \, \tau_N(\lambda)\\
\geq (\im  z)^{3} \, \int_{\mathbb{R}^{+}} \frac{1}{|\lambda - z^{2}|^{2}} d \, \tau_N(\lambda)
\end{align*}

We observe that for $\omega > 0$, then,
$$
\int_{\mathbb{R}^{+}} \frac{1}{|\lambda - z^{2}|^{2}} d \, \tau_N(\lambda) \geq
\int_{0}^{\omega} \frac{1}{|\lambda - z^{2}|^{2}} d \, \tau_N(\lambda) \geq
\frac{1}{2(\omega^{2} + |z|^{4})} \, \tau_N([0,\omega])
$$
As justified above, it is possible to choose $\omega$ for which $\tau_N([0,\omega]) \geq \frac{a}{2}$ for each $N$. This leads to (\ref{eq:lower-bound-imboldbeta}).

We now remark that   $|1-c_N^2\bm{\beta}_N^2| =|\bm{\beta}_N||\frac{1}{\bm{\beta}_N}-c_N^2\bm{\beta}_N|$. As $\im  {\bs \beta}_N > 0$ on $\mathbb{C}^{+}$, it holds that
$$
\left|\frac{1}{\bm{\beta}_N}-c_N^2\bm{\beta}_N\right| \geq \left| \im  \left( \frac{1}{\bm{\beta}_N}-c_N^2\bm{\beta}_N \right) \right| \geq c_N^{2} \im  \bm{\beta}_N
$$
Using that $|\bs{\beta}_N| \geq \im  \bs{\beta}_N$, we eventually obtain that
$$
|1-c_N^2\bm{\beta}_N^2| \geq c_N^{2} \,  \left( \im  \bm{\beta}_N \right)^{2}
$$
which, in turn, implies (\ref{eq:lowerbound-(1-z(cboldbeta)2-1}).  $\blacksquare$

\section{Expression of matrix $\ex\{\mathbf{Q}\}$ obtained using the integration by parts formula}
\label{subsec:ipp}
We now express $\ex\{\mathbf{Q}\}$ using the integration by parts formula. For this, we have first to establish
some useful properties of $\ex\{\mathbf{Q}(z)\}$ that follow from the invariance properties of the probability
distribution  of the observations $(y_n)_{n=1, \ldots, N}$. In the following, for $k,l \in \{1, 2, \ldots, L \}$, we denote by ${\bf Q}_{pp}^{k,l}$ and ${\bf Q}_{ff}^{k,l}$ the $M \times M$ matrices whose entries are given by
$\left({\bf Q}_{pp}^{k,l}\right)_{m,n} = \left({\bf Q}_{pp}\right)_{(k-1)M+m, (l-1)M+n}$ and $\left({\bf Q}_{ff}^{k,l}\right)_{m,n} = \left({\bf Q}_{ff}\right)_{(k-1)M+m, (l-1)M+n}$for each $m,n \in \{1, 2, \ldots, M \}$.

\begin{lemma}
  \label{le:symetries-E(Q)}
  The matrices $\ex\{\mathbf{Q_{pp}}\}$ and $\ex\{\mathbf{Q_{ff}}\}$ are block diagonal, 
i.e. $\mathbb{E}\left({\bf Q}_{pp}^{k,l}\right) =  \mathbb{E}\left({\bf Q}_{ff}^{k,l}\right) = 0$ if $k \neq l$, and 
	\begin{align}\label{q-pp}
	&\tr\ex\{\mathbf{Q_{pp}}\}(I_L\otimes R)=\tr\ex\{\mathbf{Q_{ff}}\}(I_L\otimes R),\\
	&\ex\{\mathbf{Q_{pf}}\}=\ex\{\mathbf{Q_{fp}}\}=0.\label{q-pf}
	\end{align}
\end{lemma}
\textbf{Proof.} To prove (\ref{q-pf}) we consider the new set of vectors $z_k=e^{-ik\theta}y_k$ and construct the matrices $Z_p$, $Z_f$ in the same way as $Y_p$ and $Y_f$. It is clear 
that sequence $(z_n)_{n \in \mathbb{Z}}$ has the same probability distribution 
that $(y_n)_{n \in \mathbb{Z}}$. $Z_p$ and $Z_f$ can be expressed as
\begin{align*}
&Z_p=\begin{pmatrix}
e^{-i\theta}I_M&\ldots&0\\
\vdots&\ddots&\vdots\\
0&\ldots&e^{-Li\theta}I_M
\end{pmatrix}Y_p\begin{pmatrix}
1&\ldots&0\\
\vdots&\ddots&\vdots\\
0&\ldots&e^{-(N-1)i\theta}
\end{pmatrix},\\
&Z_f=e^{-Li\theta}\begin{pmatrix}
e^{-i\theta}I_M&\ldots&0\\
\vdots&\ddots&\vdots\\
0&\ldots&e^{-Li\theta}I_M
\end{pmatrix}Y_f\begin{pmatrix}
1&\ldots&0\\
\vdots&\ddots&\vdots\\
0&\ldots&e^{-(N-1)i\theta}
\end{pmatrix}.
\end{align*}
Therefore, it holds that
\begin{align*}
Z_fZ_p^*Z_pZ_f^*=\begin{pmatrix}
e^{-i\theta}I_M&\ldots&0\\
\vdots&\ddots&\vdots\\
0&\ldots&e^{-Li\theta}I_M
\end{pmatrix} Y_fY_p^*Y_pY_f^*\begin{pmatrix}
e^{i\theta}I_M&\ldots&0\\
\vdots&\ddots&\vdots\\
0&\ldots&e^{Li\theta}I_M
\end{pmatrix}.
\end{align*}
Similarly to $\mathbf{Q}$ we define matrix $
\mathbf{Q^Z}=\left(\begin{smallmatrix}
-zI_{ML}&\frac{1}{N}Z_fZ_p^*\\
\frac{1}{N}Z_pZ_f^*&-zI_{ML}
\end{smallmatrix}\right)^{-1}
$
and obtain immediately that 
\begin{align*}
\ex\{\mathbf{Q^Z_{pp}}\}=\begin{pmatrix}
e^{-i\theta}I_M&\ldots&0\\
\vdots&\ddots&\vdots\\
0&\ldots&e^{-Li\theta}I_M
\end{pmatrix} \ex\{\mathbf{Q_{pp}}\}\begin{pmatrix}
e^{i\theta}I_M&\ldots&0\\
\vdots&\ddots&\vdots\\
0&\ldots&e^{Li\theta}I_M
\end{pmatrix}. 
\end{align*}
Since $\ex\{\mathbf{Q^Z_{pp}}\}=\ex\{\mathbf{Q_{pp}}\}$ then for any $M\times M$ block $\ex\{\mathbf{Q_{pp}}^{j,k}\}$ we have
\begin{align*}
\ex\{\mathbf{Q_{pp}}^{j,k}\}=e^{-ji\theta}\ex\{\mathbf{Q_{pp}}^{j,k}\}e^{ki\theta}=e^{(k-j)i\theta}\ex\{\mathbf{Q_{pp}}^{j,k}\}.
\end{align*}
This proves that $\ex\{\mathbf{Q_{pp}}^{j,k}\} = 0$ if $k \neq j$ as expected. A similar proof leads to the conclusion that $\ex\{\mathbf{Q_{ff}}\}$
is block diagonal. Moroever, the equality $\ex\{\mathbf{Q^Z_{fp}}\}=\ex\{\mathbf{Q_{fp}}\}$ implies that 
\begin{align*}
\ex\{\mathbf{Q^Z_{fp}}\}=e^{-Li\theta}\begin{pmatrix}
e^{-i\theta}I_M&\ldots&0\\
\vdots&\ddots&\vdots\\
0&\ldots&e^{-Li\theta}I_M
\end{pmatrix} \ex\{\mathbf{Q_{fp}}\}\begin{pmatrix}
e^{i\theta}I_M&\ldots&0\\
\vdots&\ddots&\vdots\\
0&\ldots&e^{Li\theta}I_M
\end{pmatrix} .
\end{align*}
Therefore, each $M \times M$ block $\mathbf{Q_{fp}}^{j,k}$ of $\mathbf{Q_{fp}}$ verifies  $\ex\{\mathbf{Q_{fp}}^{j,k}\}=e^{-(L+j-k)i\theta}\ex\{\mathbf{Q_{fp}}^{j,k}\}$. As $j-k \in \{-(L-1), \ldots, L-1 \}$, this implies that $\ex\{\mathbf{Q_{fp}}^{j,k}\}=0$. This leads immediately to $\ex\{\mathbf{Q_{fp}}\}=0$. We obtain similarly that $\ex\{\mathbf{Q_{pf}}\}=0$. 

To prove (\ref{q-pp}) let us consider sequence $z$ defined by $z_n=y_{-n+N+2L}$ 
for each $n$. Again, the distribution of $z_n$ will remain the same and it is easy to see that $Z_p$ and $Z_f$ are given by 
\begin{align*}
Z_f=\begin{pmatrix}
0&\ldots&I_M\\
\vdots&&\vdots\\
I_M&\ldots&0
\end{pmatrix}Y_p\begin{pmatrix}
0&\ldots&1\\
\vdots&&\vdots\\
1&\ldots&0
\end{pmatrix},\\
Z_p=\begin{pmatrix}
0&\ldots&I_M\\
\vdots&&\vdots\\
I_M&\ldots&0
\end{pmatrix}Y_f\begin{pmatrix}
0&\ldots&1\\
\vdots&&\vdots\\
1&\ldots&0
\end{pmatrix}.	
\end{align*}
From this, we obtain that  
\begin{align*}
\ex\{\mathbf{Q^Z_{pp}}\}=\begin{pmatrix}
0&\ldots&I_M\\
\vdots&&\vdots\\
I_M&\ldots&0
\end{pmatrix}\ex\{\mathbf{Q_{ff}}\}\begin{pmatrix}
0&\ldots&I_M\\
\vdots&&\vdots\\
I_M&\ldots&0
\end{pmatrix}.
\end{align*}	
As $\ex\{\mathbf{Q^Z_{pp}}\} = \ex\{\mathbf{Q_{pp}}\}$, this immediately implies that $\ex\{\mathbf{Q_{ff}}^{j,j}\}=\ex\{\mathbf{Q_{pp}}^{L-j,L-j}\}$, and, as a consequence, that $\ex\{\tr\mathbf{Q_{pp}}(I_L\otimes R)\}=\ex\{\tr\mathbf{Q_{ff}}(I_L\otimes R)\}$ as expected.  $\blacksquare$ \\

Now we return to the expression for ${\bf Q}(z)$.
Using the resolvent identity we get
\begin{align}\label{res id}
z\mathbf{Q}(z)=-I_{2ML}+\mathbf{Q}(z)\mathbf{M}=-I_{2ML}+\sum_{j=1}^{N}\mathbf{Q}(z)
\begin{pmatrix}
0 & w_{f,j}w_{p,j}^*\\
w_{p,j}w_{f,j}^*&0
\end{pmatrix}.
\end{align} 
For every $m_1,m_2=1,\ldots, M$, $i_1=1,\ldots,2L$ and $i_2=1,\ldots,L$ we denote by  
$\mathbf{\hat{A}}_{i_1i_2}^{m_1m_2}$ the $2N\times 2N$ matrix defined by
\begin{align}\label{matr_A}
\begin{split}
&(\mathbf{\hat{A}}_{i_1i_2}^{m_1m_2}(pf))_{jk}=\left(\mathbf{Q}
\left(\begin{smallmatrix}
0 \\
w_{p,j}
\end{smallmatrix}\right)\right)_{i_1}^{m_1}(w_{f,k}^*)_{i_2}^{m_2}, \\
&(\mathbf{\hat{A}}_{i_1i_2}^{m_1m_2}(pp))_{jk}=\left(\mathbf{Q}
\left(\begin{smallmatrix}
0 \\
w_{p,j}
\end{smallmatrix}\right)\right)_{i_1}^{m_1}(w_{p,k}^*)_{i_2}^{m_2}, \\
&(\mathbf{\hat{A}}_{i_1i_2}^{m_1m_2}(ff))_{jk}=\left(\mathbf{Q}
\left(\begin{smallmatrix}
w_{f,j} \\
0
\end{smallmatrix}\right)\right)_{i_1}^{m_1}(w_{f,k}^*)_{i_2}^{m_2},\\
&(\mathbf{\hat{A}}_{i_1i_2}^{m_1m_2}(fp))_{jk}=\left(\mathbf{Q}
\left(\begin{smallmatrix}
w_{f,j} \\
0
\end{smallmatrix}\right)\right)_{i_1}^{m_1}(w_{p,k}^*)_{i_2}^{m_2},\\
\end{split}
\end{align}
We also define matrix $\mathbf{A}_{i_1i_2}^{m_1m_2}$ by $\mathbf{A}_{i_1i_2}^{m_1m_2}=\ex\{\mathbf{\hat{A}}_{i_1i_2}^{m_1m_2}\}$. (\ref{res id}) implies that 
\begin{align}\label{res_id}
z \mathbb{E} \{\mathbf{Q}_{i_1i_2}^{m_1m_2}(z) \} =-\delta_{i_1,i_2}\delta_{m_1,m_2}+ \tr\mathbf{A}_{i_1i_2}^{m_1m_2}(pf)+\tr\mathbf{A}_{i_1i_2}^{m_1m_2}(fp).
\end{align} 
In the reminder of this paragraph, we evaluate for each $i_1,i_2,m_1,m_2$ the elements of matrix 
$\mathbf{A}_{i_1i_2}^{m_1m_2}$ using the Gaussian tools (\ref{integr}) and (\ref{deriv}). As we shall see, 
each element of $\mathbf{A}_{i_1i_2}^{m_1m_2}$ can be written as a functional of matrix $\mathbb{E}({\bf Q})$ 
plus an error term whose contribution vanishes when $N \rightarrow +\infty$. Plugging these expressions 
of $\mathbf{A}_{i_1i_2}^{m_1m_2}$ into (\ref{res_id}) will establish an approximate expression of 
$\mathbb{E}({\bf Q})$. As the calculations are very tedious, we just indicate how each element 
$(\mathbf{A}_{i_1i_2}^{m_1m_2}(ff))_{j,k}$ of matrix $\mathbf{A}_{i_1i_2}^{m_1m_2}(ff)$ can be evaluated. 
By using integration by parts formula and (\ref{deriv}) we obtain 
\begin{multline*}
\ex\left\{\left(\mathbf{Q}
\begin{pmatrix}
w_{f,j} \\
0
\end{pmatrix}\right)_{i_1}^{m_1}(w_{f,k}^*)_{i_2}^{m_2}\right\}
=\sum_{i_3=1}^{L}\sum_{m_3}\ex\{\mathbf{Q}_{i_1i_3}^{m_1m_3}W_{i_3+L,j}^{m_3}\overline{W}_{i_2+L,k}^{m_2}\}\\
=\sum_{i_3=1}^{L}\sum_{\substack{i^{\prime},j^{\prime}\\m^{\prime},m_3}}\ex\{W_{i_3+L,j}^{m_3}\overline{W}_{i^{\prime},j^{\prime}}^{m^{\prime}}\}
\times\ex\left\{\dfrac{\partial \left(\mathbf{Q}_{i_1i_3}^{m_1m_3}\overline{W}_{i_2+L,k}^{m_2}\right)}{\partial \overline{W}_{i^{\prime},j^{\prime}}^{m^{\prime}}}\right\}
=\dfrac{1}{N}\sum_{i_3=1}^{L}\sum_{\substack{i^{\prime},j^{\prime}\\m^{\prime},m_3}}R_{m_3m^{\prime}}\\
\times\delta_{i_3+L+j,i^{\prime}+j^{\prime}}\ex\left\{\mathbf{Q}_{i_1i_3}^{m_1m_3}\delta_{m_2,m^{\prime}}\delta_{i_2+L,i^{\prime}}\delta_{k,j^{\prime}}+\overline{W}_{i_2+L,k}^{m_2}\dfrac{\partial \mathbf{Q}_{i_1i_3}^{m_1m_3}}{\partial\overline{W}_{i^{\prime},j^{\prime}}^{m^{\prime}}}\right\}\\
=\dfrac{1}{N}\sum_{i_3=1}^{L}\sum_{m_3=1}^{M}\ex\left\{\mathbf{Q}_{i_1i_3}^{m_1m_3}R_{m_3m_2}\delta_{i_3,i_2-(j-k)}\right\}
-\dfrac{1}{N}\sum_{\substack{i_3,j^{\prime}\\m_3,m^{\prime}}}\sum_{i^{\prime}=1}^{L}R_{m_3m^{\prime}}\delta_{i_3+L+j,i^{\prime}+j^{\prime}}\\
\times\ex\left\{\overline{W}^{(f)m_2}_{i_2,k}\left(\mathbf{Q}\left(\begin{smallmatrix}
w_{f,j^{\prime}}\\
0
\end{smallmatrix}\right)\right)^{m_1}_{i_1}\mathbf{Q}_{i^{\prime}+Li_3}^{m^{\prime}m_3}\right\}
-\dfrac{1}{N}\sum_{\substack{i_3,j^{\prime}\\m_3, m^{\prime}}}\sum_{i^{\prime}
=L+1}^{2L}R_{m_3m^{\prime}}\delta_{i_3+L+j,i^{\prime}+j^{\prime}}\\
\times\ex\left\{\overline{W}^{(f)m_2}_{i_2,k}\left(\mathbf{Q}\left(\begin{smallmatrix}
0\\
w_{p,j^{\prime}}
\end{smallmatrix}\right)\right)^{m_1}_{i_1}\mathbf{Q}_{i^{\prime}-Li_3}^{m^{\prime}m_3}\right\}
=\dfrac{1}{N}\sum_{i_3=1}^{L}\ex\Big\{\left(\left(\begin{smallmatrix}
\mathbf{Q_{pp}}\\
\mathbf{Q_{fp}}
\end{smallmatrix}\right)(I_L\otimes R)\right)_{i_1i_3}^{m_1m_2}\\
\times\delta_{i_3,i_2-(j-k)}\Big\}
-\dfrac{1}{N}\sum_{m^{\prime},j^{\prime}}\sum_{i_3,i^{\prime}=1}^{L}\delta_{i_3+L+j,i^{\prime}+j^{\prime}}
\ex\left\{\left(\hat{A}_{i_1i_2}^{m_1m_2}(ff)\right)_{j^{\prime},k}(\mathbf{Q_{fp}}(I_L\otimes R))_{i^{\prime}i_3}^{m^{\prime}m^{\prime}}\right\}\\
-\dfrac{1}{N}\sum_{m^{\prime},j^{\prime}}\sum_{i_3,i^{\prime}=1}^{L}\delta_{i_3+j,i^{\prime}+j^{\prime}}\ex\left\{\left(\mathbf{\hat{A}}_{i_1i_2}^{m_1m_2}(pf)\right)_{j^{\prime},k}(\mathbf{Q_{pp}}(I_L\otimes R))_{i^{\prime}i_3}^{m^{\prime}m^{\prime}}\right\}
 \end{multline*}
 Now we define for every $i_1=1,\ldots, 2L$, $i_2=1,\ldots, L$ and $m_1,m_2=1,\ldots, M$ $2N\times 2N$ matrix $\mathbf{B}_{i_1i_2}^{m_1m_2}$ with blocks
 
 \begin{multline*}
\Big(\mathbf{B}_{i_1i_2}^{m_1m_2}(fp)\Big)_{j,k}=\dfrac{1}{N}\ex\left\{\left(\begin{smallmatrix}
 \mathbf{Q_{pp}}\\
 \mathbf{Q_{fp}}
 \end{smallmatrix}\right)(I_L\otimes R)\right\}_{i_1,i_2-(j-k)-L}^{m_1,m_2}\mathbf{1}_{1\le i_2-(j-k)-L\le L},\\
\Big(\mathbf{B}_{i_1i_2}^{m_1m_2}(ff)\Big)_{j,k}=\dfrac{1}{N}\ex\left\{\left(\begin{smallmatrix}
  \mathbf{Q_{pp}}\\
 \mathbf{ Q_{fp}}
  \end{smallmatrix}\right)(I_L\otimes R)\right\}_{i_1,i_2-(j-k)}^{m_1,m_2}\mathbf{1}_{1\le i_2-(j-k)\le L},\\
\Big(\mathbf{B}_{i_1i_2}^{m_1m_2}(pp)\Big)_{j,k}=\dfrac{1}{N}\ex\left\{\left(\begin{smallmatrix}
   \mathbf{Q_{pf}}\\
   \mathbf{Q_{ff}}
   \end{smallmatrix}\right)(I_L\otimes R)\right\}_{i_1,i_2-(j-k)}^{m_1,m_2}\mathbf{1}_{1\le i_2-(j-k)\le L},\\
\Big(\mathbf{B}_{i_1i_2}^{m_1m_2}(pf)\Big)_{j,k}=\dfrac{1}{N}\ex\left\{\left(\begin{smallmatrix}
    \mathbf{Q_{pf}}\\
    \mathbf{Q_{ff}}
    \end{smallmatrix}\right)(I_L\otimes R)\right\}_{i_1,i_2-(j-k)+L}^{m_1,m_2}\mathbf{1}_{1\le i_2-(j-k)+L\le L}.
 \end{multline*}
 
Also for every $ML\times ML$ block matrix $\mathbf{D}$ we define the sequence $(\tau^{(M)}(\mathbf{D})(l))_{l=-L+1,\ldots,L-1}$ as
\begin{equation}
  \label{eq:def-tau}
\tau^{(M)}(\mathbf{D})(l)=\dfrac{1}{ML}\tr \mathbf{D}(J^l_L\otimes I_M) = \frac{1}{ML} \sum_{m=1}^{M} \sum_{i-i^{\prime}=l} {\bf D}_{i,i^{\prime}}^{m,m}
\end{equation}
and $N \times N$ Toeplitz matrix $\mathcal{T}^{(M)}_{N,L}(\mathbf{D})$ given by
\begin{align}
\mathcal{T}^{(M)}_{N,L}(\mathbf{D})=\sum_{l=-L+1}^{L-1}\tau^{(M)}(\mathbf{D})(l)J^{-l}_N.
\end{align}
In other words, the entries of $\mathcal{T}^{(M)}_{N,L}(\mathbf{D})$ are defined by the relation
\begin{equation}
  \label{eq:entries-tauMNL}
  \left[\mathcal{T}^{(M)}_{N,L}(\mathbf{D})\right]_{j_1,j_2} = \tau^{(M)}(\mathbf{D})(j_1-j_2) \, \mathbf{1}_{-(L-1) \leq j_1-j_2 \leq L-1}
\end{equation}
We observe that if $\mathbf{D}$ is block diagonal, i.e. if $\mathbf{D}^{m_1,m_2}_{i_1,i_2} = 0$ 
for each $m_1,m_2$ when $i_1 \neq i_2$, then, matrix $\mathcal{T}^{(M)}_{N,L}(\mathbf{D})$ coincides with
the diagonal matrix $\mathcal{T}^{(M)}_{N,L}(\mathbf{D}) = \left( \frac{1}{ML} \mathrm{Tr} \mathbf{D} \right) \, I_N$. 
It clear that
$$
\dfrac{1}{N}\sum_{i_3=1}^{L}\ex\left\{\left(\left(\begin{smallmatrix}
\mathbf{Q_{pp}}\\
\mathbf{Q_{fp}}
\end{smallmatrix}\right)(I_L\otimes R)\right)_{i_1i_3}^{m_1m_2}\delta_{i_3,i_2-(j-k)}\right\} = \Big(\mathbf{B}_{i_1i_2}^{m_1m_2}(ff)\Big)_{j,k}
$$
In order to rewrite the term
$$
\dfrac{1}{N}\sum_{m^{\prime},j^{\prime}}\sum_{i_3,i^{\prime}=1}^{L}\delta_{i_3+L+j,i^{\prime}+j^{\prime}}\\
\times\ex\left\{\left(\hat{A}_{i_1i_2}^{m_1m_2}(ff)\right)_{j^{\prime},k}(\mathbf{Q_{fp}}(I_L\otimes R))_{i^{\prime}i_3}^{m^{\prime}m^{\prime}}\right\}
$$
in a more convenient way, we put $l=i^{\prime}-i_3$, and remark that
\begin{multline*}
\dfrac{1}{N}\sum_{m^{\prime},j^{\prime}}\sum_{i_3,i^{\prime}=1}^{L}\delta_{i_3+L+j,i^{\prime}+j^{\prime}}
  \times\ex\left\{\left(\hat{A}_{i_1i_2}^{m_1m_2}(ff)\right)_{j^{\prime},k}(\mathbf{Q_{fp}}(I_L\otimes R))_{i^{\prime}i_3}^{m^{\prime}m^{\prime}}\right\} =  \\
  \frac{ML}{N} \sum_{m^{\prime}} \sum_{l=-(L-1)}^{L-1}  \ex\left\{\left(\hat{A}_{i_1i_2}^{m_1m_2}(ff)\right)_{L+j-l,k} \frac{1}{ML} \sum_{i^{\prime}-i_3=l}
  (\mathbf{Q_{fp}}(I_L\otimes R))_{i^{\prime}i_3}^{m^{\prime}m^{\prime}}\right\}
\end{multline*}
Using the definition (\ref{eq:def-tau}), this can be rewritten as
\begin{multline*}
\dfrac{1}{N}\sum_{m^{\prime},j^{\prime}}\sum_{i_3,i^{\prime}=1}^{L}\delta_{i_3+L+j,i^{\prime}+j^{\prime}}
  \times\ex\left\{\left(\hat{A}_{i_1i_2}^{m_1m_2}(ff)\right)_{j^{\prime},k}(\mathbf{Q_{fp}}(I_L\otimes R))_{i^{\prime}i_3}^{m^{\prime}m^{\prime}}\right\} =  \\
  c_N \, \sum_{l=-(L-1)}^{L-1}  \ex\left\{\left(\hat{A}_{i_1i_2}^{m_1m_2}(ff)\right)_{L+j-l,k} \tau^{M}\left( \mathbf{{Q}_{fp}}(I_L \otimes R)\right)(l) \right\}
\end{multline*}
We introduce $j^{\prime}=L+j-l$, and using (\ref{eq:entries-tauMNL}), we notice that
\begin{multline*}
\dfrac{1}{N}\sum_{m^{\prime},j^{\prime}}\sum_{i_3,i^{\prime}=1}^{L}\delta_{i_3+L+j,i^{\prime}+j^{\prime}}
  \times\ex\left\{\left(\hat{A}_{i_1i_2}^{m_1m_2}(ff)\right)_{j^{\prime},k}(\mathbf{Q_{fp}}(I_L\otimes R))_{i^{\prime}i_3}^{m^{\prime}m^{\prime}}\right\} =  \\
  c_N \, \ex\left\{ \sum_{j^{\prime}=1}^{N} \left[\mathcal{T}^{(M)}_{N,L}\left( \mathbf{{Q}_{fp}}(I_L \otimes R)\right)\right]_{L+j,j^{\prime}}  \left(\hat{A}_{i_1i_2}^{m_1m_2}(ff)\right)_{j^{\prime},k}  \right\}= \\
  c_N\ex\left\{\left(J_N^{L}\mathcal{T}_{N,L}^{(M)}(\mathbf{Q_{fp}}(I_L\otimes R))\mathbf{\hat{A}}_{i_1i_2}^{m_1m_2}(ff)\right)_{j,k}\right\}
\end{multline*}
We obtain similarly that
\begin{multline*}
  \dfrac{1}{N}\sum_{m^{\prime},j^{\prime}}\sum_{i_3,i^{\prime}=1}^{L}\delta_{i_3+j,i^{\prime}+j^{\prime}}\ex\left\{\left(\mathbf{\hat{A}}_{i_1i_2}^{m_1m_2}(pf)\right)_{j^{\prime},k}(\mathbf{Q_{pp}}(I_L\otimes R))_{i^{\prime}i_3}^{m^{\prime}m^{\prime}}\right\} = \\
  c_N\ex\left\{\left(\mathcal{T}_{N,L}^{(M)}(\mathbf{Q_{pp}}(I_L\otimes R))\mathbf{\hat{A}}_{i_1i_2}^{m_1m_2}(pf)\right)_{j,k}\right\}
  \end{multline*}
Therefore, matrix $\mathbf{A}_{i_1i_2}^{m_1m_2}(ff)$ is also defined by
 \begin{align*}
   &\Big(\mathbf{A}_{i_1i_2}^{m_1m_2}(ff)\Big)_{j,k}=\Big(\mathbf{B}_{i_1i_2}^{m_1m_2}(ff)\Big)_{j,k}-c_N\ex\left\{\left(J_N^{L}\mathcal{T}_{N,L}^{(M)}(\mathbf{Q_{fp}}(I_L\otimes R))\mathbf{\hat{A}}_{i_1i_2}^{m_1m_2}(ff)\right)_{j,k}\right\} \\
   &-c_N\ex\left\{\left(\mathcal{T}_{N,L}^{(M)}(\mathbf{Q_{pp}}(I_L\otimes R))\mathbf{\hat{A}}_{i_1i_2}^{m_1m_2}(pf)\right)_{j,k}\right\}
   \end{align*}
 Writing $\mathbf{Q_{fp}}$ and $\mathbf{Q_{pp}}$ as $\mathbf{Q_{fp}} = \mathbb{E}\left(\mathbf{Q_{fp}}\right) + \mathbf{Q^{\circ}_{fp}} =  \mathbf{Q^{\circ}_{fp}}$ (see
     (\ref{q-pf})) and
  $\mathbf{Q_{pp}} = \mathbb{E}\left(\mathbf{Q_{pp}}\right) + \mathbf{Q^{\circ} _{pp}}$, we obtain that 
       \begin{multline*}
   \Big(\mathbf{A}_{i_1i_2}^{m_1m_2}(ff)\Big)_{j,k}=\Big(\mathbf{B}_{i_1i_2}^{m_1m_2}(ff)\Big)_{j,k}  -c_N\ex\left\{\left(\mathcal{T}_{N,L}^{(M)}(\mathbf{Q_{pp}}(I_L\otimes R))\mathbf{A}_{i_1i_2}^{m_1m_2}(pf)\right)_{j,k}\right\}\\
 -c_N\ex\left\{\left(J_N^{L}\mathcal{T}_{N,L}^{(M)}(\mathbf{Q^{\circ}_{fp}}(I_L\otimes R))\mathbf{\hat{A}}_{i_1i_2}^{m_1m_2}(ff)\right)_{j,k}\right\}\\
 -c_N\ex\left\{\left(\mathcal{T}_{N,L}^{(M)}(\mathbf{Q^{\circ}_{pp}}(I_L\otimes R))\mathbf{\hat{A}}_{i_1i_2}^{m_1m_2}(pf)\right)_{j,k}\right\}
       \end{multline*}
       We define the $N \times N$ matrix  $\mathbf{\Delta}_{i_1i_2}^{m_1m_2}(ff)$ by
       \begin{multline*}
       \mathbf{\Delta}_{i_1i_2}^{m_1m_2}(ff) = -c_N\ex\left\{J_N^{L}\mathcal{T}_{N,L}^{(M)}(\mathbf{Q^{\circ}_{fp}}(I_L\otimes R))\mathbf{\hat{A}}_{i_1i_2}^{m_1m_2}(ff)\right\}\\
       -c_N\ex\left\{\mathcal{T}_{N,L}^{(M)}(\mathbf{Q^{\circ}_{pp}}(I_L\otimes R))\mathbf{\hat{A}}_{i_1i_2}^{m_1m_2}(pf)\right\}
       \end{multline*}
       Dropping the indices $i_1$, $i_2$, $m_1$, $m_2$, we eventually obtain that
 \begin{align*}
\mathbf{A_{ff}}=\mathbf{B_{ff}}-c_N\ex\left\{\mathcal{T}_{N,L}^{(M)}(\mathbf{Q_{pp}}(I_L\otimes R))\right\}\mathbf{A_{pf}}+\mathbf{\Delta_{ff}}.
 \end{align*}
Using similar calculations, it is possible to establish that:
 \begin{align*}
 &\mathbf{A_{pf}}=\mathbf{B_{pf}}-c_N\ex\left\{\mathcal{T}_{N,L}^{(M)}(\mathbf{Q_{ff}}(I_L\otimes R))\right\}\mathbf{A_{ff}}+\mathbf{\Delta_{pf}}\\
 &\mathbf{A_{fp}}=\mathbf{B_{fp}}-c_N\ex\left\{\mathcal{T}_{N,L}^{(M)}(\mathbf{Q_{pp}}(I_L\otimes R))\right\}\mathbf{A_{pp}}+\mathbf{\Delta_{fp}}\\
 &\mathbf{A_{pp}}=\mathbf{B_{pp}}-c_N\ex\left\{\mathcal{T}_{N,L}^{(M)}(\mathbf{Q_{ff}}(I_L\otimes R))\right\}\mathbf{A_{fp}}+\mathbf{\Delta_{pp}},
 \end{align*}
 where $\mathbf{\Delta_{pf}}$, $\mathbf{\Delta_{fp}}$, and $\mathbf{\Delta_{pp}}$ are defined as
 \begin{align*}
&\mathbf{\Delta_{pf}}=-c_N\ex\left\{\mathcal{T}_{N,L}^{(M)}(\mathbf{Q^{\circ}_{pf}}(I_L\otimes R))J_N^{*L}\mathbf{\hat{A}_{pf}}\right\}
-c_N\ex\left\{\mathcal{T}_{N,L}^{(M)}(\mathbf{Q^{\circ}_{ff}}(I_L\otimes R))\mathbf{\hat{A}_{ff}}\right\},\\
&\mathbf{\Delta_{fp}}=-c_N\ex\left\{J_N^L\mathcal{T}_{N,L}^{(M)}(\mathbf{Q^{\circ}_{fp}}(I_L\otimes R))\mathbf{\hat{A}_{fp}}\right\}
-c_N\ex\left\{\mathcal{T}_{N,L}^{(M)}(\mathbf{Q^{\circ}_{pp}}(I_L\otimes R))\mathbf{\hat{A}_{pp}}\right\},\\
&\mathbf{\Delta_{pp}}=-c_N\ex\left\{\mathcal{T}_{N,L}^{(M)}(\mathbf{Q^\circ_{pf}}(I_L\otimes R))J_N^{*L}\mathbf{\hat{A}_{pp}}\right\}
-c_N\ex\left\{\mathcal{T}_{N,L}^{(M)}(\mathbf{Q^{\circ}_{ff}}(I_L\otimes R))\mathbf{\hat{A}_{fp}}\right\}.
 \end{align*}	
 By Lemma \ref{le:symetries-E(Q)}, matrices $\mathbb{E}\left(\mathbf{Q_{ff}}\right)$ and  $\mathbb{E}\left(\mathbf{Q_{pp}}\right)$ are block diagonal.
 Therefore, matrices  $\ex\{\mathcal{T}_{N,L}^{(M)}(\mathbf{Q_{ff}}(I_L\otimes R))\}$ and  $\ex\{\mathcal{T}_{N,L}^{(M)}(\mathbf{Q_{pp}}(I_L\otimes R))\}$
 reduce to $\frac{1}{ML} \ex\{\tr\mathbf{Q_{ff}}(I_L\otimes R)\} \, I_N$ and  $\frac{1}{ML} \ex\{\tr\mathbf{Q_{pp}}(I_L\otimes R)\} \, I_N$ respectively.
 As $\ex\{\tr\mathbf{Q_{ff}}(I_L\otimes R)\} = \ex\{\tr\mathbf{Q_{pp}}(I_L\otimes R)\}$ (see (\ref{q-pp})), 
 we eventually obtain that
 \begin{align}\label{sys}
 \begin{pmatrix}
 I_N & \dfrac{c_N}{ML}\ex\left\{\tr \mathbf{Q_{pp}}(I_L\otimes R)\right\}I_N\\
 \dfrac{c_N}{ML}\ex\left\{\tr \mathbf{Q_{pp}}(I_L\otimes R)\right\}I_N&I_N
 \end{pmatrix}\mathbf{A}
 =
 \mathbf{B}
 +
 \mathbf{\Delta}.
 \end{align}
 In the following, we denote by ${\bs \alpha}(z)$ the function defined by
 \begin{equation}
   \label{eq:def-blodalpha}
         {\bs \alpha}(z) = \dfrac{1}{ML}\ex\left\{\tr \mathbf{Q_{pp}}(I_L\otimes R)\right\}.
         \end{equation}
  To find the expression of $\mathbf{A}$, we have to prove that the  matrix  governing the linear system (\ref{sys}), is invertible.
 For this, we recall that $\mathbf{Q_{pp}}(z) = z Q(z^{2})$, and introduce the function $\alpha(z)$ defined by
 $$
 \alpha(z) = \frac{1}{ML} \mathrm{Tr}\left( \mathbb{E}\{ Q(z)(I_L \times R)\} \right).
 $$
 $\alpha$ is clearly an element of $\mathcal{S}(\mathbb{R}^{+})$. In order to evaluate its associated positive measure $\overline{\mu}_N$, we denote by
 $\hat{\mu}_N$ the positive measure defined by
 \begin{equation}
   \label{eq:def-hatmu}
   d \hat{\mu}_N(\lambda) = \frac{1}{ML} \sum_{i=1}^{ML} \hat{f}_i^{*} (I_L \otimes R) \hat{f}_i \, \delta_{\hat{\lambda}_i},
 \end{equation}
 where we recall that $(\hat{\lambda}_i)_{i=1, \ldots, ML}$ and $(\hat{f}_i)_{i=1, \ldots, ML}$ represent the eigenvalues and
 eigenvectors of $W_f W_p^{*} W_p W_f^{*}$. 
 We remark that $\hat{\mu}$ is carried by $\mathbb{R}^{+}$ and that its mass $\hat{\mu}(\mathbb{R}^{+})$ coincides with $\frac{1}{M} \mathrm{Tr}R$.
 Then, measure $\overline{\mu}_N$ is defined by
 \begin{equation}
   \label{eq:def-overlinemu}
 \int_{\mathbb{R}^{+}} \phi(\lambda) \, d\overline{\mu}_N(\lambda) = \mathbb{E} \left( \int_{\mathbb{R}^{+}} \phi(\lambda) \, d \hat{\mu}_N(\lambda) \right).
 \end{equation}
 Moreover, we notice that
   $$
   {\bs \alpha}(z) = z \alpha(z^{2}).
   $$
   Therefore, Lemma \ref{beta} implies that ${\bs \alpha} \in \mathcal{S}(\mathbb{R})$ and that
   $$
   1 - c_N^{2}  {\bs \alpha}(z)^{2} \neq 0
   $$
   if $z \in \mathbb{C}^{+}$. This implies that the matrix governing the linear system (\ref{sys}) is invertible for $z \in \mathbb{C}^{+}$.
   Matrix $\mathbf{H}$ given by
\begin{align*}
\mathbf{H}=\begin{pmatrix}
I_N & c_N {\bs \alpha}(z) \, I_N\\
 c_N {\bs \alpha}(z) \, I_N&I_N
\end{pmatrix}^{-1}.
\end{align*}
is thus well defined for each $z \in \mathbb{C}^{+}$. 

The blocks of $\mathbf{H}$ are of course given by
\begin{align*}
&\mathbf{H_{pp}}=\mathbf{H_{ff}}=\frac{1}{1 - c_N^{2} {\bs \alpha(z)}^{2}} \, I_N\\
&\mathbf{H_{pf}}=\mathbf{H_{fp}}=-\frac{c_N {\bs \alpha(z)}}{1 - c_N^{2} {\bs \alpha(z)}^{2}} \, I_N.
\end{align*}
(\ref{sys}) implies that $\mathbf{A}=\mathbf{H}\mathbf{B}+\mathbf{H}\mathbf{\Delta}$.  (\ref{res_id}) implies that we only need to evaluate matrices $\mathbf{A_{pf}}$ and $\mathbf{A_{fp}}$. We obtain that these matrices are given by
\begin{align*}
&\mathbf{A_{pf}}=\mathbf{H_{pp}}\mathbf{B_{pf}}+\mathbf{H_{pf}}\mathbf{B_{ff}}+\mathbf{H_{pp}}\mathbf{\Delta_{pf}}+\mathbf{H_{pf}}\mathbf{\Delta_{ff}}\\
&\mathbf{A_{fp}}=\mathbf{H_{fp}}\mathbf{B_{pp}}+\mathbf{H_{ff}}\mathbf{B_{fp}}+\mathbf{H_{fp}}\mathbf{\Delta_{pp}}+\mathbf{H_{ff}}\mathbf{\Delta_{fp}}.
\end{align*}
This and the definition (\ref{matr_A}) of matrix $\mathbf{A}_{i_1i_2}^{m_1m_2}$ lead immediately to
\begin{multline*}
\left(\ex\left\{\mathbf{Q}\left(\begin{smallmatrix}
0& W_fW_p^*\\
W_pW_f^*&0
\end{smallmatrix}\right)\right\}\right)_{i_1i_2}^{m_1m_2}=
\tr \mathbf{A}_{i_1i_2}^{m_1m_2}(pf)\mathbf{1}_{i_2\le L}+
\tr \mathbf{A}_{i_1i_2-L}^{m_1m_2}(pf)\mathbf{1}_{i_2>L}=\\
\dfrac{1}{1-c_N^2\bm{\alpha}^2}\tr \Big( \mathbf{B_{pf}}-c_N\bm{\alpha} \mathbf{B_{ff}}
+ \mathbf{\Delta_{pf}}-c_N\bm{\alpha}\mathbf{\Delta_{ff}}\Big)_{i_1i_2}^{m_1m_2}\mathbf{1}_{i_2\le L}\\
+\dfrac{1}{1-c_N^2\bm{\alpha}^2}\tr \Big( \mathbf{B_{fp}}-c_N\bm{\alpha} \mathbf{B_{pp}}+ \mathbf{\Delta_{fp}}-c_N\bm{\alpha}\mathbf{\Delta_{pp}}\Big)_{i_1i_2-L}^{m_1m_2}\mathbf{1}_{i_2> L}
\end{multline*}
It is easy to notice that $\tr \left(\mathbf{B_{fp}}\right)_{i_1i_2}^{m_1m_2}=\tr \left(\mathbf{B_{pf}} \right)_{i_1i_2}^{m_1m_2}=0$, and $\tr \left(\mathbf{B_{pp}}\right)_{i_1i_2}^{m_1m_2}=\ex\{\left(\mathbf{Q}\Pi_{ff}(I_2L\otimes R)\right)_{i_1i_2+L}^{m_1m_2}\}$, $\tr \left(\mathbf{B_{ff}}\right)_{i_1i_2}^{m_1m_2}=\ex\{\left(\mathbf{Q}\Pi_{pp}(I_2L\otimes R)\right)_{i_1i_2}^{m_1m_2}\}$, where $\Pi_{ff}=\left(\begin{smallmatrix}
0&0\\
0&I_{ML}
\end{smallmatrix}\right)$ and $\Pi_{pp}=\left(\begin{smallmatrix}
I_{ML}&0\\0&0
\end{smallmatrix}\right)$. Hence,
\begin{multline*}
\left(\ex\left\{\mathbf{Q}\left(\begin{smallmatrix}
0& W_fW_p^*\\
W_pW_f^*&0
\end{smallmatrix}\right)\right\}\right)_{i_1i_2}^{m_1m_2}=-\dfrac{c_N\bm{\alpha}}{1-c_N^2\bm{\alpha}^2}\Big(\ex\{\mathbf{Q}\Pi_{pp}(I_{2L}\otimes R)\}\\
+\ex\{\mathbf{Q}\Pi_{ff}(I_{2L}\otimes R)\}\Big)_{i_1i_2}^{m_1m_2}+\mathcal{E}_{i_1i_2}^{m_1m_2}
=-\dfrac{c_N \bm{\alpha}}{1-c_N^2\bm{\alpha}^2}\Big(\ex\{\mathbf{Q}(I_{2L}\otimes R)\}\Big)_{i_1i_2}^{m_1m_2} +\mathcal{E}_{i_1i_2}^{m_1m_2},
\end{multline*}
where $\mathcal{E}_{i_1i_2}^{m_1m_2}$ represents the remaining terms depending on the entries of matrix $\mathbf{\Delta}_{i_1i_2}^{m_1m_2}$. Using the identity
(\ref{res id}), we obtain that 
\begin{align}
\label{eq:first-equation-E(Q)}
z\ex\{\mathbf{Q}\}+I_{2ML}=\ex\left\{\mathbf{Q}\left(\begin{smallmatrix}
0& W_fW_p^*\\
W_pW_f^*&0
\end{smallmatrix}\right)\right\}=-\dfrac{c_N\bm{\alpha}}{1-c_N^2\bm{\alpha}^2}\ex\{\mathbf{Q}\}(I_{2L}\otimes R)+\mathcal{E},
\end{align}
which immediately leads to 
\begin{align*}
- \ex\{\mathbf{Q}\} \left(\dfrac{c_N\bm{\alpha}}{1-c_N^2\bm{\alpha}^2}(I_{2L}\otimes R)+z\right)= I_{2ML} - \mathcal{E}
\end{align*}
As $\mathbb{E}({\bf Q})$ is block diagonal, (\ref{eq:first-equation-E(Q)}) implies 
that matrix $\mathcal{E}$ is also block diagonal, i.e. $\mathcal{E}_{fp} = \mathcal{E}_{pf}=0$. 
We apply Lemma \ref{beta} to $\beta(z) = \alpha(z)$, and conclude that matrix $-\left(\dfrac{c_N\bm{\alpha}}{1-c_N^2\bm{\alpha}^2}(I_{2L}\otimes R)+z\right)$ is invertible
for each $z \in \mathbb{C}^{+}$, and that matrix $\mathbf{S}_N(z)$, defined by
\begin{equation}
  \label{eq:def-boldS}
 \mathbf{S}_N(z)=-\left(\dfrac{c_N\bm{\alpha}(z)}{1-c_N^2\bm{\alpha}^2(z)}R+z\right)^{-1}
\end{equation}
belongs to $\mathcal{S}_{M}(\mathbb{R})$, and verifies $\| \mathbf{S}_N(z) \| \leq \frac{1}{\im z}$. We deduce from this that
$$
\ex\{\mathbf{Q}\} = -\left(\dfrac{c_N\bm{\alpha}}{1-c_N^2\bm{\alpha}^2}(I_{2L}\otimes R)+z\right)^{-1}  + \mathcal{E} \left(\dfrac{c_N\bm{\alpha}}{1-c_N^2\bm{\alpha}^2}(I_{2L}\otimes R)+z\right)^{-1}
$$
or equivalently that
\begin{equation}
  \label{eq:expre-E(boldQ)}
\ex\{\mathbf{Q}(z)\} = I_{2L} \otimes {\bf S}(z) -  \mathcal{E}(z) \,\left(I_{2L} \otimes  {\bf S}(z)\right)
\end{equation}
This allows to evaluate $\mathbb{E}(Q(z))$ by identification of the first diagonal blocks of the 
left and right hand sides of (\ref{eq:expre-E(boldQ)}). For this, we introduce 
the $M \times M$ matrix-valued function $S_N(z)$ defined by
\begin{equation}
\label{eq:def-S}
S_N(z)=-\left(z I_M + \dfrac{c_N z \alpha_N(z)}{1-c_N^2\alpha_N(z)^2}R_N \right)^{-1}
\end{equation}
Lemma \ref{beta} implies that $S$ belongs to $\mathcal{S}_{M}(\mathbb{R}^{+})$, verifies $\| S_N(z) \| \leq \frac{1}{\im z}$, and that $\mathbf{S}(z)$ and $S(z)$ are linked by the equation $\mathbf{S}(z) = z S(z^{2})$. 
As $\mathbb{E}({\bf Q}_{pp}(z)) = z \mathbb{E}(Q(z^{2}))$,  (\ref{eq:expre-E(boldQ)}) leads to 
\begin{equation}
  \label{eq:expre-EQz2}
\mathbb{E}(Q(z^{2})) = I_L \otimes S(z^{2}) - \mathcal{E}_{pp}(z) \, I_L \otimes S(z^{2})
\end{equation}
for each $z \in \mathbb{C}^{+}$. Therefore, $\mathcal{E}_{pp}(z)$ only depends on $z^{2}$. As 
the image of $\mathbb{C}^{+}$ by the transformation $z \rightarrow z^{2}$ is $\mathbb{C} - \mathbb{R}^{+}$, 
we obtain that $\mathcal{E}_{pp}(z) = E_{pp}(z^{2})$ for some function $E_{pp}$ analytic in 
$\mathbb{C} - \mathbb{R}^{+}$. This discussion leads to 
\begin{equation}
\label{eq:expre-E(Q)}
\mathbb{E}(Q(z)) = I_L \otimes S(z) - E_{pp}(z) \, \left( I_L \otimes S(z) \right)
\end{equation}
for each $z \in \mathbb{C} - \mathbb{R}^{+}$. 

In the following, we prove that 
\begin{equation}
  \label{eq:behaviour-error}
  \dfrac{1}{ML} \mathrm{Tr} \left( \mathbb{E}(Q_N(z)) - I_L \otimes S_N(z) \right) = -\dfrac{1}{ML}\tr( E_{pp}(z)  (I_L \otimes S_N(z)) ) = 
   \mathcal{O}_z(\frac{1}{N^{2}})
\end{equation}

\section{Evaluation the error term $\mathcal{E}$}\label{sec:evaluation_epsilon}

In order to establish (\ref{eq:behaviour-error}), we prove the following result.
\begin{proposition}
\label{prop:convergence-calE}
For each deterministic $ML \times ML$
sequence of matrices $(F_{1,N})_{N \geq 1}$ such that $\sup_{N \geq 1} \| F_{1,N}\| \leq \kappa$, then
\begin{equation}
  \label{eq:behaviour-error-general}
  \left| \dfrac{1}{ML}\tr(\mathcal{E}_{pp}(z) \, F_{1,N}) \right|  \leq \kappa \, \frac{1}{N^{2}} \, P_1(|z^{2}|) P_{2}(\frac{1}{\im  z^{2}})
\end{equation}
holds for each $z$ for which $\im z^{2} > 0$, where $P_1$ and $P_2$ are 2 nice polynomials. 
\end{proposition}
{\bf Proof.}
We define $F$ as the $2ML \times 2ML$ matrix  $F_N = \left( \begin{array}{cc} F_{1,N} & 0 \\ 0 & 0 \end{array} \right)$
and remark that $\dfrac{1}{ML} \tr \mathcal{E}F =  \dfrac{1}{ML}\tr(\mathcal{E}_{pp}(z) \, F_{1,N})$ can be written as
\begin{multline}
  \label{eq:expre-1-error}
\dfrac{1}{ML}\tr\mathcal{E}F=\dfrac{1}{1-c^{2} \bm{\alpha}^2}\sum_{\substack{i_1,i_2\\m_1,m_2}}\Big(\left(\tr\mathbf{\Delta}_{i_1i_2}^{m_1m_2}(pf)-c \bm{\alpha}\tr\mathbf{\Delta}_{i_1i_2}^{m_1m_2}(ff)\right)\mathbf{1}_{i_2\le L}\\
+\left(\tr\mathbf{\Delta}_{i_1i_2-L}^{m_1m_2}(fp)-c \bm{\alpha}\tr\mathbf{\Delta}_{i_1i_2-L}^{m_1m_2}(pp)\right)\mathbf{1}_{i_2> L}\Big)F_{i_2i_1}^{m_2m_1}.
\end{multline}
As matrix $F$ verifies $F_{i_2,i_1}^{m_2,m_1} = 0$ if $i_2 > L$, $\dfrac{1}{ML}\tr\mathcal{E}F$ is reduced to the
right hand side of (\ref{eq:expre-1-error}) that we now evaluate. 
\begin{multline*}
\sum_{\substack{i_1,i_2\\m_1,m_2}}\tr\mathbf{\Delta}_{i_1i_2}^{m_1m_2}(pf)F_{i_2i_1}^{m_2m_1}\mathbf{1}_{i_2\le L}
=c \sum_{\substack{i_1,i_2\\m_1,m_2}}\sum_{j,k}\ex\Big\{\mathcal{T}_{N,L}^M(\mathbf{Q^{\circ}_{ff}}(I_L\otimes R))_{jk}\Big(\mathbf{Q}\left(\begin{smallmatrix}
w_{f,k}\\0
\end{smallmatrix}\right)\Big)_{i_1}^{m_1}\\
\times\Big(w_{f,j}^*\Big)_{i_2}^{m_2}F_{i_2i_1}^{m_2m_1}
+(\mathcal{T}_{N,L}^M(\mathbf{Q^{\circ}_{pf}}(I_L\otimes R))J_N^{*L})_{jk}\Big(\mathbf{Q}\left(\begin{smallmatrix}
0\\w_{p,k}
\end{smallmatrix}\right)\Big)_{i_1}^{m_1}\Big(w_{f,j}^*\Big)_{i_2}^{m_2}F_{i_2i_1}^{m_2m_1}
\Big\}\mathbf{1}_{i_2\le L} \\
=c \, \tr\ex\Big\{\mathcal{T}_{N,L}^M(\mathbf{Q^{\circ}_{ff}}(I_L\otimes R))
\left(\begin{smallmatrix}
W_f\\0
\end{smallmatrix}\right)^*F\mathbf{Q}\left(\begin{smallmatrix}
W_f\\0
\end{smallmatrix}\right) 
+\mathcal{T}_{N,L}^M(\mathbf{Q^{\circ}_{pf}}(I_L\otimes R))J_N^{*L}
\left(\begin{smallmatrix}
W_f\\0
\end{smallmatrix}\right)^*F\mathbf{Q}\left(\begin{smallmatrix}
0\\W_p
\end{smallmatrix}\right)
\Big\} \\
=c \, \tr\ex\Big\{\mathcal{T}_{N,L}^M(\mathbf{Q^{\circ}_{ff}}(I_L\otimes R))
\left(\Pi_{pf}W\right)^*F\mathbf{Q}\left(\Pi_{pf}W\right) \\
+\mathcal{T}_{N,L}^M(\mathbf{Q^{\circ}_{pf}}(I_L\otimes R))J_N^{*L}
\left(\Pi_{pf}W\right)^*F\mathbf{Q}\left(\Pi_{fp}W\right)
\Big\}.
\end{multline*}
Similar calculations lead to the following expression of $\dfrac{1}{ML}\tr\mathcal{E}F$:
\begin{multline}
    \label{eq:expre-horrible-error}
 \dfrac{1}{ML}\tr\mathcal{E}F= \dfrac{c}{(1-c_N^2\bm{\alpha}^2)} 
 \frac{1}{ML} \tr\ex\Big\{\mathcal{T}_{N,L}^M(\mathbf{Q^{\circ}_{ff}}(I_L\otimes R))
\left(\Pi_{pf}W\right)^*F\mathbf{Q}\left(\Pi_{pf}W\right)\\
+\mathcal{T}_{N,L}^M(\mathbf{Q^{\circ}_{pf}}(I_L\otimes R))J_N^{*L}
\left(\Pi_{pf}W\right)^*F\mathbf{Q}\left(\Pi_{fp}W\right)
-c \bm{\alpha}\mathcal{T}_{N,L}^M(\mathbf{Q^{\circ}_{pp}}(I_L\otimes R))
\left(\Pi_{pf}W\right)^*F\mathbf{Q}\left(\Pi_{fp}W\right)\\
-c \bm{\alpha}J_N^{L}\mathcal{T}_{N,L}^M(\mathbf{Q^{\circ}_{fp}}(I_L\otimes R))
\left(\Pi_{pf}W\right)^*F\mathbf{Q}\left(\Pi_{pf}W\right)\Big\}
\end{multline}
We now evaluate the right hand side of (\ref{eq:expre-horrible-error}). The Schwartz inequality leads to 
\begin{multline*}
\left| \dfrac{1}{ML}\tr\ex\Big\{\mathcal{T}_{N,L}^M(\mathbf{Q^{\circ}_{ff}}(I_L\otimes R))
\left(\Pi_{pf}W\right)^*F\mathbf{Q}\left(\Pi_{pf}W\right)\Big\} \right|\\
=\left|\sum_{l=-L+1}^{L-1}\ex\Big\{\tau^{(M)}(\mathbf{Q_{ff}^{\circ}}(I_L\otimes R))(l)
\frac{1}{ML} \mathrm{Tr}\Big(J^{*l}_N \left(\Pi_{pf}W\right)^*F\mathbf{Q}\left(\Pi_{pf}W\right)\Big)\Big\} \right|\\
=\left|\sum_{l=-L+1}^{L-1}\ex\Big\{\frac{1}{ML}\tr(\mathbf{Q_{ff}^{\circ}}(I_L\otimes R)(J_L^l\otimes I_M))\dfrac{1}{ML}\tr\Big(J^{*l}_N\left(\Pi_{pf}W\right)^*F\mathbf{Q}\left(\Pi_{pf}W\right)\Big)^{\circ}\Big\}  \right|\\
\le\sum_{l=-L+1}^{L-1}\var\left \{\frac{1}{ML}\tr(\mathbf{Q_{ff}}(I_L\otimes R)(J_L^l\otimes I_M)) \right\}^{1/2} \\
\times\var\left\{\dfrac{1}{ML}\tr\Big(J^{*l}_{(N)}\left(\Pi_{pf}W\right)^*F\mathbf{Q}\left(\Pi_{pf}W\right)\Big)\right\}^{1/2}
\end{multline*}
Using Lemma~\ref{var}, we obtain that 
\begin{align*}
&\var\left \{\frac{1}{ML}\tr(\mathbf{Q_{ff}}(I_L\otimes R)(J_L^l\otimes I_M)) \right\} \le  \dfrac{1}{N^{2}}P_1(|z^{2}|)P_2\left(\dfrac{1}{\im  z^{2}}\right)
\end{align*}
and that
\begin{align*}
& \var\left\{\dfrac{1}{ML}\tr\Big(J^{*l}_{(N)}\left(\Pi_{pf}W\right)^*F\mathbf{Q}\left(\Pi_{pf}W\right)\Big)\right\} \le \kappa^{2} \dfrac{1}{N^{2}}P_1(|z^{2}|)P_2\left(\dfrac{1}{\im  z^{2}}\right)
\end{align*}
Since $L$ does not grow with $N$ this implies immediately 
$$
\left| \dfrac{1}{ML}\tr\ex\Big\{\mathcal{T}_{N,L}^M(\mathbf{Q^{\circ}_{ff}}(I_L\otimes R)) \left(\Pi_{pf}W\right)^*F\mathbf{Q}\left(\Pi_{pf}W\right)\Big\} \right| \leq  \kappa \dfrac{1}{N^{2}}P_1(|z^{2}|)P_2\left(\dfrac{1}{\im  z^{2}}\right)
$$
It can be shown similarly that the 3 other normalized traces can be upper bounded by the same kind of term. 
It remains to control the terms $\frac{1}{1 - (c_N \, \bm{\alpha}_N)^{2}}$ and $\frac{\bm{\alpha}_N}{1 - (c_N \, \bm{\alpha}_N)^{2}}$. For this, we use 
Lemma \ref{le:gestion-(1-z(cbeta)2} for the choice $\beta_N(z) = \alpha_N(z)$. 
It is sufficient to verify that the measures $(\overline{\mu}_N)_{N \geq 1}$ associated to 
functions $(\alpha_N(z))_{N \geq 1}$ verify (\ref{eq:mass-tau}) and (\ref{eq:first-moment-tau}).
For each $N$, it holds that
$$
\int_{0}^{+\infty}  d \, \overline{\mu}_N(\lambda)  =  \mathbb{E} \left( \int_{0}^{+\infty} d \, \hat{\mu}_N(\lambda) \right) = \frac{1}{M} \mathrm{Tr} R_N
$$
and
$$
\int_{0}^{+\infty}\lambda \,  d \, \overline{\mu}_N(\lambda)  = \mathbb{E} \left( \int_{0}^{+\infty}\lambda \, d \, \hat{\mu}_N(\lambda) \right) = \mathbb{E} \left( \frac{1}{ML} \mathrm{Tr}((I_L \otimes R) W_f W_p^{*} W_p W_f^{*}) \right)
$$
A straightforward calculation leads to $ \mathbb{E} \left( \frac{1}{ML} \mathrm{Tr}(W_f W_p^{*} W_p W_f^{*}) \right) = c_N \frac{1}{M} \mathrm{Tr} R_N \frac{1}{M} \mathrm{Tr} R_N^{2}$. Therefore, (\ref{eq:lowerbound-(1-z(cbeta)2-1}) implies that
$$
\frac{1}{|1 -  z  (c_N \alpha_N(z))^{2}|} \leq P_1(|z|) P_2(\frac{1}{\im z})
$$
for each $z \in \mathbb{C}^{+}$, and  if $z^{2} \in \mathbb{C}^{+}$, it holds that
$$
\frac{1}{|1 -  z^{2}  (c_N \alpha_N(z^{2}))^{2}|} \leq P_1(|z^{2}|) P_2(\frac{1}{\im z^{2}})
$$
As $\bs{\alpha}_N(z) = z \alpha_N(z^{2})$, this is equivalent to
$$
\frac{1}{1 - (c_N \, \bm{\alpha}_N)^{2}} \leq P_1(|z^{2}|) P_2(\frac{1}{\im z^{2}})
$$
Finally, we remark that $|\alpha_N(z)| \leq \frac{1}{M} \mathrm{Tr} R_N \, \frac{1}{\im z} \leq b \, \frac{1}{\im z}$ for each $z \in \mathbb{C}^{+}$. Therefore, if $z^{2} \in \mathbb{C}^{+}$, it holds that $|\alpha_N(z^{2})| \leq  b \, \frac{1}{\im z^{2}}$
and that $|\bm{\alpha}_N(z)| = |z| |\alpha_N(z^{2})|$ verifies
$$
|\bm{\alpha}_N(z)|\leq b |z| \,  \frac{1}{\im z^{2}} \leq b (1+|z|^{2})  \,  \frac{1}{\im z^{2}}
$$
This completes the proof of Proposition \ref{prop:convergence-calE}.   $\blacksquare$ \\

Proposition \ref{prop:convergence-calE} immediately leads to the following Corollary. 
\begin{corollary}
\label{coro:convergence-E(Q)-S}
For each sequence $(F_N)_{N \geq 1}$ of deterministic $ML \times ML$ matrices such that 
$\sup_{N \geq 1} \leq \kappa$, then, we have 
\begin{equation}
\label{eq:convergence-E(Q)-S}
\left| \frac{1}{ML} \mathrm{Tr} \left[ \left( \mathbb{E}(Q_N(z)) - I_L \otimes S_N(z) \right) F_N \right] \right|
\leq \kappa \,  \frac{1}{N^{2}} P_1(|z|) P_{2}(\frac{1}{\im  z}) 
\end{equation}
for each $z \in \mathbb{C}^{+}$. In particular, it holds that 
\begin{equation}
\label{eq:convergence-E(Q)-S-F=I}
\left| \frac{1}{ML} \mathrm{Tr} \left[ \left( \mathbb{E}(Q_N(z)) - I_L \otimes S_N(z) \right) \right] \right|
\leq \kappa \,  \frac{1}{N^{2}} P_1(|z|) P_{2}(\frac{1}{\im  z}) 
\end{equation}
\end{corollary}
    {\bf Proof.} (\ref{eq:expre-EQz2}) implies that
    $$
    \left| \frac{1}{ML} \mathrm{Tr} \left[ \left( \mathbb{E}(Q_N(z^{2})) - I_L \otimes S_N(z^{2}) \right) F_N \right] \right| =
    \left| \frac{1}{ML} \mathrm{Tr} \mathcal{E}_{pp}(z) S_N(z^{2}) F_N \right|
    $$
    As $\| S_N(z^{2}) \| \leq \frac{1}{\im z^{2}}$ if $z^{2} \in \mathbb{C}^{+}$, the application of Proposition
    \ref{prop:convergence-calE} to matrix $F_{1,N} = S_N(z^{2}) F_N$ implies that
    $$
    \left| \frac{1}{ML} \mathrm{Tr} \left[ \left( \mathbb{E}(Q_N(z^{2})) - I_L \otimes S_N(z^{2}) \right) F_N \right] \right| \leq
    \kappa \frac{1}{N^{2}} P_1(|z^{2}|) P_2(\frac{1}{\im z^{2}})
    $$
    for each $z$ such that $z^{2} \in \mathbb{C}^{+}$. Exchanging $z^{2}$ by $z$ eventually establishes (\ref{eq:convergence-E(Q)-S}). 

\section{Deterministic equivalent of $\ex\{Q\}$}
\label{sec:deterministic-equivalent}

\subsection{The canonical equation}

\begin{proposition}
\label{prop:existence-uniqueness}
	If $z \in \mathbb{C}^{+}$, there exists a unique solution of the equation 
	\begin{align}\label{t=}
	t_N(z)=\dfrac{1}{M}\tr R_N\left(-zI_M-\dfrac{zc_Nt_N(z)}{1- zc_N^2t_N^2(z)}R_N\right)^{-1}
	\end{align}
        satisfying $t_N(z) \in \mathbb{C}^{+}$ and $z t_N(z) \in \mathbb{C}^{+}$. Function $z \rightarrow t_N(z)$ is an element of $\mathcal{S}(\mathbb{R}^+)$, and the associated positive measure, denoted by $\mu_N$, verifies 
\begin{equation}
\label{eq:masse-first-moment-mu}
\mu_N(\mathbb{R}^{+}) = \frac{1}{M} \mathrm{Tr} R_N, \; \int_{\mathbb{R}^{+}} \lambda \, d \, \mu_N(\lambda) = c_N \, \frac{1}{M} \mathrm{Tr} R_N \, \frac{1}{M} \mathrm{Tr} R_N^{2}
\end{equation}
Moreover, it exists nice constants $\beta$ and $\kappa$ such that
\begin{equation}
  \label{eq:control-inverse-1-z(ct)2}
  \frac{1}{\left|1 - z \left(c_N \, t_N(z)\right)^{2}\right|} \leq \frac{\kappa \, (\beta^{2} + |z|^{2})^{2}}{ \left( \im  z \right)^{3}}
    \end{equation}
for each $N$. 
Finally, the $M \times M$ valued function $T_N(z)$ defined by 
\begin{equation}
\label{eq:def-T}
T_N(z) = - \left(zI_M+ \dfrac{zc_Nt_N(z)}{1- zc_N^2t_N^2(z)}R_N\right)^{-1}
\end{equation}
belongs to $\mathbb{S}_M(\mathbb{R}^{+})$. The associated $M \times M$ positive matrix-valued measure, denoted
$\nu_N^{T}$, verifies 
\begin{equation}
\label{eq:mass-nuT}
\nu_N^{T}(\mathbb{R}^{+}) = I_M
\end{equation}
as well as 
\begin{equation}
\label{eq:lien-muT-nu}
\mu_N = \frac{1}{M} \mathrm{Tr} R_N  \nu_N^{T}
\end{equation}
\end{proposition}

\textbf{Proof.} As $N$ is assumed to be fixed in the statement of the Proposition, we omit to mention that $t_N, T_N, \mu_N, \ldots$
depend on $N$ in the course of the proof. We first prove the existence of a solution such that  $z \rightarrow t(z)$ is an element of
$\mathcal{S}(\mathbb{R}^+)$. For this, we use the classical fixed point equation scheme. We define $t_0(z)=-\frac{1}{z}$, which is of course
an element of $\mathcal{S}(\mathbb{R}^+)$, 
and generate sequence $(t_n(z))_{n \geq 1}$ by the formula
\begin{align*}
t_{n+1}(z)=\dfrac{1}{M}\tr R\left(-zI_M-\dfrac{zct_{n}(z)}{1- zc^2t_{n}^2(z)}R\right)^{-1}.
\end{align*}
We establish by induction that for each $n$, $t_n \in \mathcal{S}(\mathbb{R^+})$, 
and that its associated measure $\mu_n$ verifies $\mu_n(\mathbb{R^+})=\frac{1}{M}\tr R$
and 
\begin{align}\label{tight}
\int_{0}^{+\infty}\lambda\mu_n(d\lambda) = c \frac{1}{M} \tr(R) \frac{1}{M} \tr(R^{2}) 
\end{align} 
Thank's to (\ref{eq:hypothesis-R-bis}), this last property will imply that sequence $(\mu_n)_{n \geq 1}$ is tight. We assume that $t_n$ indeed satisfies the above conditions, and prove that $t_{n+1}(z)$ 
also meets these requirements. Lemma \ref{beta} implies that function  $T_n(z)=\left(-zI_M-\dfrac{zct_n(z)}{1-zc^2t_n^2}R\right)^{-1}$ is an element
of $\mathcal{S}_M(\mathbb{R^+})$. According to Proposition~\ref{stil_tran}, to prove that $t_{n+1}(z)\in \mathcal{S}(\mathbb{R}^+)$, we need to check that $\im  t_{n+1}(z), \im  zt_{n+1}(z)>0$ if $z \in \mathbb{C}^{+}$, as well as that $\lim_{y\rightarrow+\infty}iyt_{n+1}(iy)$ exists. As $T_n \in \mathcal{S}_M(\mathbb{R^+})$ and $t_{n+1}(z)=\frac{1}{M}\tr R T_n(z)$,
it is clear  that  $\im  t_{n+1}(z), \im  zt_{n+1}(z)>0$. 
Finally, it holds that
\begin{align*}
 -iyt_{n+1}(iy)=\dfrac{1}{M}\tr R\left(I_M+\dfrac{ciyt_n(iy)}{iy-(ciyt_n(iy))^2}R\right)^{-1}.
\end{align*}
Since $t_n(z)$ is a Stieltjes transform we have $-iyt_n(iy)\rightarrow\mu_n(\mathbb{R}^+)$, which
implies that  $-iyt_{n+1}(iy)\rightarrow\frac{1}{M}\tr R$, i.e. that $\mu_{n+1}(\mathbb{R}^+)=\frac{1}{M}\tr R$.

We finally check that $\mu_{n+1}$ satisfies (\ref{tight}). For this, we follow \cite{hachem-loubaton-najim-2007}.
\begin{align*}
\int_{0}^{+\infty}\lambda\mu_{n+1}(d\lambda)=\lim_{y\rightarrow+\infty}\Re\left(-iy(iy\dfrac{1}{M}\tr RT_n(iy)+\dfrac{1}{M}\tr R)\right).
\end{align*}
Using twice the resolvent identity we can express $T_{n}$ as
\begin{align*}
T_{n}=-\dfrac{1}{z}\left(I_M+\dfrac{ct_n}{1-zc^2t_n^2}R\right)^{-1}=-\dfrac{1}{z}+\dfrac{R}{z}\dfrac{ct_n}{1-zc^2t_n^2}-\left(\dfrac{ct_n}{1-zc^2t_n^2}\right)^2R^2T_{n},
\end{align*}
from which it follows that 
\begin{align*}
-z \left(\dfrac{1}{M}\tr(zRT_n(z))+\dfrac{1}{M}\tr R)\right)=-\dfrac{c zt_n}{1-zc^2t_n^2}\dfrac{1}{M}\tr R^2+\left(\dfrac{czt_n}{1-zc^2t_n^2}\right)^2\dfrac{1}{M}\tr R^3T_{n}.
\end{align*}
Since $-iyt_n(iy)\rightarrow\frac{1}{M}\tr R$ and $t_n(iy)\rightarrow0$ we can conclude that $-iy(iy\frac{1}{M}\tr RT_n(iy)+\frac{1}{M}\tr R)\rightarrow\frac{c}{M^2}\tr R\tr R^2$
as expected. \\

We now prove that sequence $t_n$ converges towards a function $t \in  \mathcal{S}(\mathbb{R^+})$
verifying equation (\ref{t=}).  For this we evaluate $\theta_n=t_{n+1}-t_n$
\begin{align*}
&\theta_n=\dfrac{1}{M}\tr R(T_n-T_{n-1})=\dfrac{1}{M}\tr RT_n\dfrac{zc(t_n-t_{n-1})(1+zc^2t_nt_{n-1})}{(1-zc^2t_n^2)(1-zc^2t_{n-1}^2)}RT_{n-1}\\
&=\theta_{n-1}\dfrac{zc(1+zc^2t_nt_{n-1})}{(1-zc^2t_n^2)(1-zc^2t_{n-1}^2)}
\dfrac{1}{M}\tr RT_nRT_{n-1}.
\end{align*} 
We denote by $f_n(z)$ the term defined by
\begin{equation}
\label{eq:def-fn}
f_n(z) = \dfrac{zc(1+zc^2t_nt_{n-1})}{(1-zc^2t_n^2)(1-zc^2t_{n-1}^2)}
\dfrac{1}{M}\tr RT_nRT_{n-1}
\end{equation}
Lemma~\ref{beta} implies that $\|T_k\|\le\frac{1}{\im  z}$ and that $|t_k| \leq \frac{b}{\im  z}$ for each 
$k \geq 1$ and each $z \in \mathbb{C}^{+}$. Therefore, it holds that
$$
\left| zc(1+zc^2t_nt_{n-1}) \dfrac{1}{M}\tr RT_nRT_{n-1} \right| \leq \kappa \left(\frac{|z|}{(\im  z)^{2}} \left(1 + \frac{|z|}{(\im  z)^{2}}\right) \right)
$$
Moreover, it is clear that for each $k$, $|1-zc^2t_k^2| \geq (1 - c^{2} \frac{|z|}{(\im  z)^{2}})$. For each $\epsilon > 0$ small enough, we consider the domain $\mathcal{D}_{\epsilon}$ defined
by 
\begin{equation}
\label{eq:def-Depsilon}
\mathcal{D}_{\epsilon} = \{z \in \mathbb{C}^{+},  \frac{|z|}{(\im  z)^{2}} < \epsilon \}
\end{equation}
Then, for $z \in \mathcal{D}_{\epsilon}$, it holds that
$$
\frac{1}{|1-zc^2t_{n}^2|} \frac{1}{|1-zc^2t_{n-1}^2|} \leq \frac{1}{(1 - c^{2} \epsilon)^{2}}
$$
and that
$$
|f_n(z)| \leq \frac{\kappa}{(1 - c^{2} \epsilon)^{2}} \left(\epsilon + \epsilon^{2} \right)
$$
We choose $\epsilon$ in such a way that $ \frac{\kappa}{(1 - c^{2} \epsilon)^{2}} \left(\epsilon + \epsilon^{2} \right) < 1/2$. Then, for
each $z \in \mathcal{D}_{\epsilon}$, it holds that
$$
|\theta_n|\le \frac{1}{2} |\theta_{n-1}|
$$
Therefore, for each $z$ in  $\mathcal{D}_{\epsilon}$, $(t_n(z))_{n \geq 1}$ is a Cauchy sequence. We denote by $t(z)$ its limit. $(t_n(z))_{n \geq 1}$ is uniformly bounded  on every compact set of $\mathbb{C} - \mathbb{R}^+$.  This implies that $(t_n(z))_{n \geq 1}$ is a normal family on
 $\mathbb{C} - \mathbb{R}^+$. We consider 
a converging subsequence extracted from $(t_n(z))_{n \geq 1}$. The corresponding limit $t_*(z)$ is analytic over $\mathbb{C} - \mathbb{R}^{+}$. If $z \in \mathcal{D}_{\epsilon}$, $t_*(z)$ must be equal to $t(z)$. Therefore,  the limits of all converging
subsequences extracted from $(t_n(z))_{n \geq 1}$ must coincide on $\mathcal{D}_{\epsilon}$, and therefore on $\mathbb{C} - \mathbb{R}^{+}$.
This implies that $t_n(z)$ converges uniformly on
each compact subset towards a function  which is analytic $\mathbb{C} - \mathbb{R}^+$, and that  we also denote by $t(z)$. It is clear that $t(z)$ verifies (\ref{t=}) and that $t \in \mathcal{S}(\mathbb{R}^{+})$
and verifies (\ref{eq:masse-first-moment-mu}). Moroever, Lemma \ref{beta} implies that $T \in \mathcal{S}_M(\mathbb{R}^{+})$, while (\ref{eq:lien-muT-nu}) and (\ref{eq:mass-nuT}) are obtained immediately. \\

As (\ref{eq:masse-first-moment-mu}) holds, (\ref{eq:control-inverse-1-z(ct)2}) is a consequence of the application of Lemma \ref{le:gestion-(1-z(cbeta)2}
to the function $\beta_N(z) = t_N(z)$.  \\

We now prove that if $z \in \mathbb{C}^{+}$ and $t_1(z)$ and $t_2(z)$ are 2 solutions of (\ref{t=}) such that
$t_i(z)$ and $z t_i(z)$ belong to $\mathbb{C}^{+}$, $i=1,2$, then $t_1(z) = t_2(z)$. In order to prove this,
we first establish the following useful Lemma.
\begin{lemma}
\label{le:properties-u-v}
If $z \in \mathbb{C}^{+}$ and if $t(z)$ verifies the conditions of Proposition \ref{prop:existence-uniqueness}, then, it holds
that
\begin{equation}
  \label{eq:signe-1-u}
  1 - u(z) > 0
\end{equation}
and
\begin{equation}
  \label{eq:det}
  \det(\mathbf{I}-\mathbf{D}) > 0
\end{equation}
where
\begin{eqnarray}
\label{eq:def-D}
\mathbf{D} & = & \begin{pmatrix}
u(z)&&v(z)\\
|z|^2v(z)&&u(z)
\end{pmatrix}\\
  \label{eq:def-u}
  u(z) & = &  c \; \frac{| c z t(z)|^{2} \, \frac{1}{M} \mathrm{Tr}(R T(z) (T(z))^{*} R)}{|1 - z (c t(z))^{2}|^{2}} \\
    \label{eq:def-v}
    v(z) & = &  c \; \frac{ \frac{1}{M} \mathrm{Tr}(R T(z) (T(z))^{*} R)}{|1 - z (c t(z))^{2}|^{2}}
\end{eqnarray}
\end{lemma}
    {\bf Proof.} Using the equation $t(z) = \frac{1}{M} \mathrm{Tr} R T(z)$, we obtain immediately after some algebra that
\begin{equation}
\label{eq:linear-system-imt-imzt}
\left( \begin{array}{c} \frac{\im (t(z))}{\im (z)}  \\  \frac{\im (z t(z))}{\im (z)} \end{array} \right) = 
\mathbf{D} \; \left( \begin{array}{c} \frac{\im (t(z))}{\im (z)}  \\  \frac{\im (z t(z))}{\im (z)} \end{array} \right) + \left( \begin{array}{c} \frac{1}{M} \mathrm{Tr}(R T(z) (T(z))^*) \\ 0  \end{array} \right)
\end{equation}
The first component of (\ref{eq:linear-system-imt-imzt}) implies that
$$
(1 - u(z)) \, \frac{\im (t(z))}{\im (z)} = v(z) \, \frac{\im (z t(z))}{\im (z)} +  \frac{1}{M} \mathrm{Tr}(R T(z) (T(z))^*)
$$
Therefore, it holds that $(1 - u(z)) >  0$. Plugging the equality
$$
\frac{\im (t(z))}{\im (z)} = \frac{v(z)}{1 - u(z)} \,  \frac{\im (z t(z))}{\im (z)} + \frac{1}{1 - u(z)} \,  \frac{1}{M} \mathrm{Tr}(R T(z) (T(z))^*)
$$
into the second component of (\ref{eq:linear-system-imt-imzt}) leads to
$$
\left( 1 - u(z) - \frac{|z|^{2} v^{2}(z)}{1 - u(z)} \right) \, \frac{\im (z t(z))}{\im (z)} =  \frac{|z|^{2} v(z)}{1 - u(z)} \,  \frac{1}{M} \mathrm{Tr}(R T(z) (T(z))^*) > 0
$$
and to (\ref{eq:det}). \\

To complete the proof of the uniqueness, we assume that equation (\ref{t=}) has 2 solutions $t_1(z)$ and
$t_2(z)$ such that $t_i(z)$ and $z t_i(z)$ belong to $\mathbb{C}^{+}$ for $i=1,2$. The proof of Lemma 
\ref{beta} (see in particular (\ref{eq:factorizarion-1-z(cbeta)2})) implies that for $i=1,2$, then $1 - z (c t_i(z))^{2} \neq 0$ and 
matrix $-zI-\dfrac{zct_i(z)}{1- zc^2t_i^2(z)}R$ is invertible. We denote by
$T_1(z)$ and $T_2(z)$ the matrices defined by (\ref{eq:def-T}) when $t(z) = t_1(z)$ and  $t(z) = t_2(z)$
respectively. $u_i(z)$ and $v_i(z)$, $i=1,2$, are defined similarly from (\ref{eq:def-u}) and (\ref{eq:def-v}) when $t(z) = t_1(z)$ and  $t(z) = t_2(z)$. Using that $t_i(z) = \frac{1}{M} \mathrm{Tr}(R T_i(z))$ for $i=1,2$, we obtain immediately that
$$
t_1(z) - t_2(z) = \left( u_{1,2}(z) + z v_{1,2}(z) \right) \, (t_1(z) - t_2(z))
$$
where 
\begin{equation}
\label{eq:def-u12}
u_{1,2}(z) =  c \, \frac{c z t_1(z) c z t_2(z) \, \frac{1}{M} \mathrm{Tr}(R T_1(z) R T_2(z))}{\left(1 - z (c t_1(z))^{2}\right) \left(1 - z (c t_2(z))^{2}\right)}
\end{equation}
 and
\begin{equation}
\label{eq:def-v12}
v_{1,2}(z) =  c \, \frac{ \frac{1}{M} \mathrm{Tr}(R T_1(z) R T_2(z))}{\left(1 - z (c t_1(z))^{2}\right) \left(1 - z (c t_2(z))^{2}\right)}
\end{equation}
In order to prove
that $t_1(z) = t_2(z)$, it is sufficient establish that $1 -  u_{1,2}(z) - z v_{1,2}(z) \neq 0$. For this, we prove the following
inequality:
\begin{equation}
  \label{eq:inequality-uniqueness}
  | 1 - u_{1,2}(z) - z v_{1,2}(z) | > \sqrt{(1 -u_1(z)) - |z| v_1(z)} \, \sqrt{(1 -u_2(z)) - |z| v_2(z)}
\end{equation}
which, by Lemma \ref{le:properties-u-v}, implies  $1 -  u_{1,2}(z) - z v_{1,2}(z)) \neq 0$. For this, we remark that
the Schwartz inequality leads to $|u_{1,2}(z)| \leq \sqrt{u_1(z)} \sqrt{u_2(z)}$ and  $|v_{1,2}(z)| \leq \sqrt{v_1(z)} \sqrt{v_2(z)}$. Therefore,
$$
| 1 - u_{1,2}(z) - z v_{1,2}(z) | \geq 1 -  \sqrt{u_1(z)} \sqrt{u_2(z)} -  \sqrt{|z| v_1(z)} \sqrt{|z| v_2(z)}
$$
We now use the inequality
\begin{equation}
  \label{eq:useful-inequality}
  \sqrt{a b} - \sqrt{c d} \geq \sqrt{a - c} \, \sqrt{b - d}
\end{equation}
where $a,b,c,d$ are positive real numbers such that $a \geq c$ and $b \geq d$. (\ref{eq:useful-inequality}) for $a=b=1$ and
$c = u_1(z)$, $d=u_2(z)$ implies that $1 -  \sqrt{u_1(z)} \sqrt{u_2(z)} \geq \sqrt{1 - u_1(z)}\sqrt{1 - u_2(z)}$. Therefore, it holds that
$$
| 1 - u_{1,2}(z) - z v_{1,2}(z) | \geq  \sqrt{1 - u_1(z)}\sqrt{1 - u_2(z)} -  \sqrt{|z| v_1(z)} \sqrt{|z| v_2(z)}
$$
(\ref{eq:useful-inequality}) for $a=1 - u_1(z)$,  $b=1 - u_2(z)$, $c = |z| v_1(z)$ and $d = |z| v_2(z)$ eventually
    leads to (\ref{eq:inequality-uniqueness}). This completes the proof of the uniqueness of the solution of
    (\ref{t=}) and Proposition~\ref{prop:existence-uniqueness}.  $\blacksquare$

\begin{remark}
  \label{re:properties-u-v-z-real-negative}
  (\ref{eq:signe-1-u}) and (\ref{eq:det}) are still valid if $z$ belongs to $\mathbb{R}^{-*}$. To check this,
  it is sufficient to remark if $z=x \in \mathbb{R}^{-*}$, the fundamental equation (\ref{eq:linear-system-imt-imzt}) is still valid, but $\frac{\im (t(z))}{\im (z)}$ and  $\frac{\im (z t(z))}{\im (z)}$ have to
  be replaced by $t^{'}(x)$ and $(x t(x))^{'}$ where $^{'}$ denotes the differentiation operator w.r.t. $x$.
  The same conclusions are obtained because  $t^{'}(x) > 0$ and $(x t(x))^{'} > 0$ if $x \in \mathbb{R}^{-*}$.
\end{remark}


    \subsection{Convergence}
    \label{subsec:convergence}

 In this paragraph, we establish that the empirical eigenvalue distribution $\hat{\nu}_N$ of 
matrix $W_{f,N}W_{p,N}^{*}W_{p,N}W_{f,N}^{*}$ has almost surely the same deterministic behaviour 
than the probability measure $\nu_N$ defined by 
\begin{equation}
\label{eq:def-nu}
\nu_N = \frac{1}{M} \mathrm{Tr} \nu_N^{T}
\end{equation}
where we recall that $\nu_N^{T}$ represents the positive matrix valued measure associated to $T_N(z)$. For this, we first establish the following Proposition. 
\begin{proposition}
\label{prop:convergence-E(Q)-T}
For each sequence $(F_N)_{N \geq 1}$ of deterministic $ML \times ML$ matrices such that 
$\sup_{N \geq 1}\|F_N\| \leq \kappa$, then, 
\begin{equation}
\label{eq:convergence-E(Q)-T}
\frac{1}{ML} \mathrm{Tr} \left[ \left( \mathbb{E}(Q_N(z)) - I_L \otimes T_N(z) \right) F_N \right] \rightarrow 0
\end{equation}
holds for each $z \in \mathbb{C} - \mathbb{R}^{+}$.
\end{proposition}
{\bf Proof.} 
Corollary \ref{coro:convergence-E(Q)-S} implies that
\begin{align*}
\dfrac{1}{ML}\tr (\ex\{Q_N\}-(I_L\otimes S_N))F_N=\mathcal{O}\left(\dfrac{1}{N^{2}}\right)
\end{align*}
We have therefore to show
 that $\frac{1}{ML} \tr\left(I_L \otimes (S_N-T_N)\right) F_N \rightarrow0$. 
It is easy to check that 
\begin{multline}\label{eq:tr_s-t}
\dfrac{1}{ML}\tr\left( I_L \otimes (S-T) \right) F
  =  \dfrac{1}{ML}\tr (I_L \otimes S)\left(\dfrac{zc_N\alpha}{1-zc_N^2\alpha^2}-\dfrac{zc_Nt}{1-zc_N^2t^2}\right)
(I_L \otimes RT) F \\
 =  \dfrac{zc_N(\alpha-t)(1+zc_N^2\alpha t)}{(1-zc_N^2\alpha^2)(1-zc_N^2t^2)}\dfrac{1}{ML}\tr (I_L \otimes SRT) F.
\end{multline}
We express $\alpha - t$ as $\alpha - \frac{1}{M} \mathrm{Tr}RS + \frac{1}{M} \mathrm{Tr}R(S-T)$, and deduce 
from (\ref{eq:tr_s-t}) that 
\begin{multline}
\label{eq:(S-T)F}
\dfrac{1}{ML}\tr\left( I_L \otimes (S-T) \right) F  =  \left( \alpha - \frac{1}{M} \mathrm{Tr}RS  \right) \dfrac{zc_N (1+zc_N^2\alpha t)}{(1-zc_N^2\alpha^2)(1-zc_N^2t^2)}\\
\times\dfrac{1}{ML}\tr (I_L \otimes SRT) F  
+ \frac{1}{M} \mathrm{Tr}R(S-T) \, \dfrac{zc_N (1+zc_N^2\alpha t)}{(1-zc_N^2\alpha^2)(1-zc_N^2t^2)}\dfrac{1}{ML}\tr (I_L \otimes SRT) F
\end{multline}
(\ref{eq:convergence-E(Q)-S}) implies that $\alpha - \frac{1}{M} \mathrm{Tr}RS = \mathcal{O}_z(\frac{1}{N^{2}})$. 
Therefore, in order to establish (\ref{eq:convergence-E(Q)-T}), it is sufficient to prove that 
$\frac{1}{M} \mathrm{Tr}R(S-T) \rightarrow 0$. For this, we take $F = I_L \otimes R$ in (\ref{eq:(S-T)F}) 
and get that 
\begin{equation}
\label{eq:(S-T)R}
\frac{1}{M} \mathrm{Tr}R(S(z)-T(z)) = f_N(z) \, \frac{1}{M} \mathrm{Tr}R(S(z)-T(z)) +   \mathcal{O}_z(\frac{1}{N^{2}})
\end{equation}
where $f_N(z)$ is defined by 
$$
f_N(z) = \dfrac{zc_N (1+zc_N^2\alpha t)}{(1-zc_N^2\alpha^2)(1-zc_N^2t^2)} \frac{1}{M} \mathrm{Tr}(RS(z)RT(z))
$$
$f_N(z)$ is similar to the term defined in (\ref{eq:def-fn}). Using the arguments of the proof of 
Proposition \ref{prop:existence-uniqueness}, we obtain that it is possible to find $\epsilon > 0$ 
for which,  $\sup_{N \geq N_0} |f_N(z)| < \frac{1}{2}$ for each $z \in \mathcal{D}_{\epsilon}$ for some 
large enough integer $N_0$. We recall that $\mathcal{D}_{\epsilon}$ is defined by (\ref{eq:def-Depsilon}). 
We therefore deduce from (\ref{eq:(S-T)R}) that $\frac{1}{M} \mathrm{Tr}R(S(z)-T(z)) \rightarrow 0$ and 
$\dfrac{1}{ML}\tr\left( I_L \otimes (S(z)-T(z)) \right) F$ converge towards $0$ for each $z \in \mathcal{D}_{\epsilon}$. As functions $z \rightarrow \dfrac{1}{ML}\tr\left( I_L \otimes (S_N(z)-T_N(z)) \right) F_N$
are holomorphic on $\mathbb{C} - \mathbb{R}^{+}$ and are uniformly bounded on each compact subset of
 $\mathbb{C} - \mathbb{R}^{+}$, we deduce from Montel's theorem that $\dfrac{1}{ML}\tr\left( I_L \otimes (S_N(z)-T_N(z)) \right) F_N$ converges towards $0$ for each $z \in \mathbb{C} - \mathbb{R}^{+}$.  $\blacksquare$ \\

We deduce the following Corollary. 
\begin{corollary}
\label{coro:weak-convergence-hatnu}
The empirical eigenvalue distribution $\hat{\nu}_N$ of $W_{f,N}W_{p,N}^{*}W_{p,N}W_{f,N}^{*}$
verifies 
\begin{equation}
\label{eq:weak-convergence-hatnu}
\hat{\nu}_N - \nu_N \rightarrow 0
\end{equation}
weakly almost surely.
\end{corollary}
{\bf Proof.} Proposition \ref{prop:convergence-E(Q)-T} implies that 
$\mathbb{E}\left(\frac{1}{ML} \tr Q_N(z) \right) - \frac{1}{M} \mathrm{Tr}(T_N(z)) \rightarrow 0$ for each 
$z \in \mathbb{C} - \mathbb{R}^{+}$. The Poincaré-Nash inequality and the Borel Cantelli Lemma 
imply that $\frac{1}{ML} \mathrm{Tr}(Q_N(z)) - \mathbb{E}\left(\frac{1}{ML} \tr Q_N(z) \right) \rightarrow 0$ a.s. for each $z \in \mathbb{C} - \mathbb{R}^{+}$. Therefore, it holds that 
\begin{equation}
\label{eq:convergence-empirical-distribution}
\frac{1}{ML} \mathrm{Tr}(Q_N(z)) - \frac{1}{M} \mathrm{Tr}(T_N(z)) \rightarrow a.s.
\end{equation}
for each $z \in \mathbb{C} - \mathbb{R}^{+}$. Corollary 2.7 of \cite{hachem-loubaton-najim-2007}
implies that $\hat{\nu}_N - \nu_N \rightarrow 0$ weakly almost surely provided we verify that 
$(\hat{\nu}_N)_{N \geq 1}$ is almost surely tight and that $(\nu_N)_{N \geq 1}$ is tight. 
It is clear that 
$$
\int_{\mathbb{R}^{+}} \lambda \, d \, \hat{\nu}_N(\lambda) = \frac{1}{ML} \mathrm{Tr} W_{f,N}W_{p,N}^{*}W_{p,N}W_{f,N}^{*} \leq \| W_N \|^{4}
$$
where we recall that 
$$
W_N = \left( \begin{array}{c} W_{p,N} \\ W_{f,N} \end{array} \right)
$$
It holds that $\| W_N \| \leq \sqrt{b} \, \| W_{iid,N} \|$ where $W_{iid,N}$ is defined 
by (\ref{eq:def-Wiid}). As $\| W_{iid,N} \| \rightarrow (1+\sqrt{c_*})$ almost surely 
(see \cite{L:15}), we obtain that $\frac{1}{ML} \mathrm{Tr} W_{f,N}W_{p,N}^{*}W_{p,N}W_{f,N}^{*}$ 
is almost surely bounded for $N$ large enough. This implies that $(\hat{\nu}_N)_{N \geq 1}$ is almost 
surely tight. As for sequence  $(\nu_N)_{N \geq 1}$, we have shown that 
$\sup_{N} \int_{\mathbb{R}^{+}} \lambda \, d \, \mu_N(\lambda) < +\infty$. 
As $\mu_N = \frac{1}{M} \tr R_N \nu_N^{T}$, the condition $R_N > a I$ for each $N$ leads to 
$$
\int_{\mathbb{R}^{+}} \lambda \, d \, \mu_N(\lambda)  \geq a \, \int_{\mathbb{R}^{+}} \lambda \, d \, \nu_N(\lambda)
$$
Therefore, it holds that $\sup_{N} \int_{\mathbb{R}^{+}} \lambda \, d \, \nu_N(\lambda) < +\infty$, a condition 
which implies that $(\nu_N)_{N \geq 1}$ is tight.  $\blacksquare$

\section{Detailed study of $\nu_N$.}
\label{sec:study-nuN}
In this section, we study the properties of $\nu_N$. (\ref{eq:hypothesis-R-bis}) implies that $\mu_N$ and $\nu_N$ are absolutely continuous one with respect each other. Hence, they share the same 
properties, and the same support denoted $\mathcal{S}_N$ in the following. We thus study $\mu_N$ 
and deduce the corresponding results related to $\nu_N$. As in the context of other models, $\mu_N$ can be characterized by studying theStieltjes transform $t_N(z)$ near the real axis. In the following, we denote by $\overline{M}$ the
number of distinct eigenvalues $(\overline{\lambda}_{l,N})_{l=1, \ldots, \overline{M}}$ 
arranged in the decreasing order,
and by $(m_{l,N})_{l=1, \ldots, \overline{M}}$ their multiplicities. It of course holds
that $\sum_{l=1}^{\bar{M}} m_{l,N} = M$. 

\subsection{Properties of $t(z)$ near the real axis.}
In this paragraph, we establish that if $x_0 \in \mathbb{R}^{+*}$, then, $\lim_{z \rightarrow x_0, z \in \mathbb{C}^{+}} t(z)$
exists and is finite. It will be denoted by $t(x_0)$ in order to simplify the notations. Moreover, when $c \leq 1$, $\lim_{z \rightarrow 0, z \in \mathbb{C}^{+} \cup \mathbb{R}^{*}} |t(z)| = +\infty$,
and  $\lim_{z \rightarrow 0, z \in \mathbb{C}^{+} \cup \mathbb{R}^{*}} z t(z) = 0$. 
The results of \cite{silverstein-choi-1995} will imply that measure $\mu_N$ is absolutely continuous w.r.t. the Lebesgue measure, and that the corresponding density is equal to $\frac{1}{\pi} \im (t(x))$ for each $x \in \mathbb{R}^{+*}$.
When $c > 1$, a Dirac mass appears at $0$. \\

We first address the case where $x_0 \neq 0$, and, in order to establish the existence of $\lim_{z \rightarrow x_0, z \in \mathbb{C}^{+}} t(z)$, we prove the following properties:
\begin{itemize}
\item If $(z_n)_{n \geq 1}$ is a sequence of $\mathbb{C}^{+}$ converging towards $x_0$, then ${|t(z_n)|}_{n \geq 1}$ is
  bounded
\item If $(z_{1,n})_{n \geq 1}$ and $(z_{2,n})_{n \geq 1}$ are two sequences of $\mathbb{C}^{+}$ converging towards $x_0$
  and verifying $\lim_{z_{i,n} \rightarrow x_0} = t_i$ for $i=1,2$, then $t_1 = t_2$.
\end{itemize}
\begin{lemma}
  \label{le:t(zn)-bounded}
  If $x_0 \in \mathbb{R}^{+*}$, and if $(z_n)_{n \geq 1}$ is a sequence of $\mathbb{C}^{+}$ such that
  $\lim_{n \rightarrow +\infty} z_n = x_0$, then the set ${|t(z_n)|}_{n \geq 1}$ is
  bounded.
\end{lemma}
    {\bf Proof.} We assume that $|t(z_n)| \rightarrow +\infty$. Equation (\ref{t=}) can be written as
    \begin{equation}
\label{eq:t=eigenvalues}
    t(z_n) = \frac{1}{M} \sum_{l=1}^{\overline{M}} \frac{m_l \, \overline{\lambda}_l}{-z_n(1 + \frac{c t(z_n) \overline{\lambda_l}}{1 - z \, (ct(z_n))^{2}})}  
    \end{equation}
    As $x_0 \neq 0$, the condition $|t(z_n)| \rightarrow +\infty$ implies that it exists $l_0$ for which
    $$
    (1 + \frac{c t(z_n) \overline{\lambda}_{l_0}}{1 - z \, (ct(z_n))^{2}}) \rightarrow 0
    $$
    or equivalently
    $$
    z_n c t(z_n) - \frac{1}{c t(z_n)} \rightarrow \overline{\lambda}_{l_0}
    $$
    As  $|t(z_n)| \rightarrow +\infty$, it holds that $ z_n c t(z_n)  \rightarrow \overline{\lambda}_{l_0}$, a contradiction.  $\blacksquare$ \\

    \begin{lemma}
      \label{le:t1=t2}
      Consider $(z_{1,n})_{n \geq 1}$ and $(z_{2,n})_{n \geq 1}$ two sequences of $\mathbb{C}^{+}$ converging towards $x_0 \in \mathbb{R}^{+*}$  and verifying $\lim_{z_{i,n} \rightarrow x_0} t(z_{i,n}) = t_i$ for $i=1,2$. Then, it holds that $t_1 = t_2$.
    \end{lemma}
        {\bf Proof.} The statement of the Lemma is obvious if $x_0$ does not belong to $\mathcal{S}$. Therefore, we assume
        that $x_0 \in \mathcal{S}-\{0 \}$. We first observe that if $\lim_{n \rightarrow +\infty} z_n = x_0$ ($z_n \in \mathbb{C}^{+}$) and
        $t(z_n) \rightarrow t_0$, then
        \begin{align}
          \label{eq:property1-limit-t(zn)}
          &1 - x_0 \, (c t_0)^{2}  \neq 0 \\
          \label{eq:property2-limit-t(zn)}
          &1 + \frac{ c t_0  \,  \overline{\lambda}_l}{1 - x_0 \, (c t_0)^{2}}  \neq 0, \;  l=1, \ldots, \overline{M}
        \end{align}
        Indeed, if (\ref{eq:property1-limit-t(zn)}) does not hold, Eq. (\ref{eq:t=eigenvalues}) leads to 
$t_0 = 0$, a contradiction because $ 1 - x_0 \, (c t_0)^{2}$ was assumed equal to $0$. Similarly, 
if (\ref{eq:property2-limit-t(zn)}) does not hold, the limit of $t(z_n)$ cannot be finite. Therefore, matrix $T_0$ defined by
        \begin{equation}
          \label{eq:def-T0}
          T_0 = - \left( x_0 \left[ I + \frac{ c t_0}{1 - x_0 \, (c t_0)^{2}} \, R \right] \right)^{-1}
        \end{equation}
        is well defined, and it holds that $T(z_n) \rightarrow T_0$ and that $t_0 = \frac{1}{M} \mathrm{Tr} R T_0$. In particular, for $i=1,2$, $T(z_{i,n}) \rightarrow T_i$
        where $T_i$ is defined by (\ref{eq:def-T0}) when $t_0 = t_i$, $i=1,2$, and $t_i =  \frac{1}{M} \mathrm{Tr} R T_i$.
        Using the equation (\ref{t=}) for $z=z_{i,n}$,
        we obtain immediately that
        \begin{eqnarray}
          \nonumber
         &  & \left( \begin{array}{c} t(z_{1,n}) - t(z_{2,n}) \\ z_{1,n} t(z_{1,n}) - z_{2,n} t(z_{2,n}) \end{array} \right)  = 
          \left( \begin{array}{cc} u_0(z_{1,n},z_{2,n}) & v_0(z_{1,n},z_{2,n}) \\ z_{1,n} z_{2,n} v_0(z_{1,n},z_{2,n}) & u_0(z_{1,n},z_{2,n})
          \end{array} \right) \\
          \label{eq:linear-system-z1-z2}
          &   &\times \left( \begin{array}{c} t(z_{1,n}) - t(z_{2,n}) \\ z_{1,n} t(z_{1,n}) - z_{2,n} t(z_{2,n}) \end{array} \right)
            + \left( \begin{array}{c}  (z_{1,n} - z_{2,n}) \, \frac{1}{M} \mathrm{Tr} T(z_{1,n}) R 
T(z_{2,n}) \\ 0 \end{array} \right)
        \end{eqnarray}
        where $u_0(z_1,z_2)$ and $v_0(z_1,z_2)$ are defined by
        \begin{equation}
\label{eq:def-u_0(z1,z2)}
u_0(z_1,z_2) =  c \, \frac{c z_1 t(z_1) c z_2 t(z_2) \, \frac{1}{M} \mathrm{Tr}(R T(z_1) R T(z_2))}{\left(1 - z_1 (c t(z_1))^{2}\right) \left(1 - z_2 (c t(z_2))^{2}\right)}
\end{equation}
 and
\begin{equation}
\label{eq:def-v_0(z1,z2)}
v_0(z_1,z_2) =  c \,  \frac{\frac{1}{M} \mathrm{Tr}(R T(z_1) R T(z_2))}{\left(1 - z_1 (c t(z_1))^{2}\right) \left(1 - z_2 (c t(z_2))^{2}\right)}
\end{equation}
for $z_i \in \mathbb{C}^{+}$, $i=1,2$. Taking the limit, we obtain that
\begin{equation}
          \left( \begin{array}{c} t_1 - t_2 \\ x_0( t_1 - t_2) \end{array} \right)  = 
          \left( \begin{array}{cc} u_0(x_0,x_0) & v_0(x_0,x_0) \\ x_0^{2} v_0(x_0,x_0) & u_0(x_0,x_0)
          \end{array} \right) \;  \left( \begin{array}{c} t_1 - t_2 \\ x_0( t_1 - t_2) \end{array} \right)
\end{equation}
where $u_0(x_0,x_0)$ and  $v_0(x_0,x_0)$ are defined by replacing $z_i, t(z_i), T(z_i)$ by $x_0, t_i, T_i$  in
(\ref{eq:def-u_0(z1,z2)}, \ref{eq:def-v_0(z1,z2)}) for $i=1,2$. If the determinant 
$(1 - u_0(x_0,x_0))^{2} - x_0^{2} v_0(x_0,x_0)^{2} \neq 0$ of the above linear system is non zero, it of course holds
that $t_1 = t_2$. \\

We now consider the case where $(1 - u_0(x_0,x_0))^{2} - x_0^{2} v_0(x_0,x_0)^{2} = 0$.
We denote by $u_i(x_0)$ and $v_i(x_0)$, $i=1,2$ the limits of $u(z_{i,n})$ and $v(z_{i,n})$, $i=1,2$ 
when $n \rightarrow +\infty$. We recall that $u(z)$ and $v(z)$ are defined by (\ref{eq:def-u}) and 
(\ref{eq:def-v}) respectively. It is clear that $u_i(x_0)$ and $v_i(x_0)$ coincide with 
(\ref{eq:def-u}) and (\ref{eq:def-v}) when $(z,t(z),T(z))$ are replaced by $(x_0, t_i, T_i)$ respectively. 
(\ref{eq:det}) thus implies that 
\begin{equation}
\label{eq:det-limit}
(1 - u_i(x_0))^{2} - x_0^{2} v_i(x_0)^{2} \geq 0
\end{equation}
for $i=1,2$.  
Using the Schwartz inequality and (\ref{eq:useful-inequality}) as in the uniqueness
proof of the solutions of Eq. (\ref{t=}) (see Proposition \ref{prop:existence-uniqueness}), it is easily seen that
\begin{multline}
\label{eq:inequalities-limit}
|(1 -u_0(x_0,x_0))^{2} - x_0^{2} (v_0(x_0,x_0))^{2}|  \geq  (1 -\sqrt{u_1(x_0)} \sqrt{u_2(x_0)})^{2} - x_0^{2} v_1(x_0) v_2(x_0) \\
   \geq  (1 - u_1(x_0)) (1 - u_2(x_0)) - x_0^{2} v_1(x_0) v_2(x_0) \\
    \geq  \sqrt{(1 - u_1(x_0))^{2} - x_0^{2} v_1(x_0)^{2}}  \sqrt{(1 - u_2(x_0))^{2} - x_0^{2} v_2(x_0)^{2}} \geq 0 
\end{multline}
Therefore,  $(1 - u_0(x_0,x_0))^{2} - x_0^{2} v_0(x_0,x_0)^{2} = 0$ implies that the Schwartz inequalities 
and the inequalities  (\ref{eq:useful-inequality}) used to establish (\ref{eq:inequalities-limit}) are equalities. 
Hence, it holds that $|u_0(x_0,x_0)|^{2} = u_1(x_0) u_2(x_0)$, or equivalently 
$| \frac{1}{M} \mathrm{Tr}(RT_1 R T_2)| = (\frac{1}{M} \mathrm{Tr}(RT_1   T_1^{*}R))^{1/2}  (\frac{1}{M} \mathrm{Tr}(RT_2  T_2^{*}R ))^{1/2}$. This implies that $T_1 = a T_2^{*}$ for some constant $a \in \mathbb{C}$. Moreover, as $t_i = \frac{1}{M} \mathrm{Tr}(R T_i)$  for $i=1,2$, it must hold that $t_1 = a t_2^{*}$. In order to prove (\ref{eq:inequalities-limit}) we use (\ref{eq:useful-inequality}) twice, for set $\{a=b=1,\; c=u_1(x_0),\;d=u_2(x_0)\}$ and set $\{a=(1-u_1(x_0))^2,\;b=(1-u_2(x_0))^2,\;c=x_0^2v_1^2,\;d=x_0^2v_2^2\}$. Since all
  these terms are positive real numbers, $\sqrt{ab} -\sqrt{cd} = \sqrt{a - c} \sqrt{b -d}$ if and only if $ad = bc$. It gives us
\begin{align}
u_1(x_0)=u_2(x_0)\\
(1-u_1(x_0))^2x_0^2v_2(x_0)^2=(1-u_2(x_0))^2x_0^2v_1(x_0)^2
\end{align}
Since $x_0\ne0$ and $(1 - u_1(x_0))^{2} - x_0^{2} v_1(x_0)^{2}\ge0$, if $u_1(x_0)=1$ it follows that $v_1(x_0)=0$ which is impossible. Hence, $u_1(x_0)\ne1$ and we have $v_1(x_0)=v_2(x_0)$. From the definition of $u_i$ and $v_i$ one can notice that $u_i(x_0)=c^2x_0^2|t_i|^2v_i(x_0)$. Which gives us immediately $|t_1|^2=|t_2|^2$ and, as a consequence, $|a|=1$. Using once again the fact that $v_1(x_0)=v_2(x_0)$ and $T_1=aT_2^*$, we obtain that 
\begin{align*}
\dfrac{|a|^2\frac{1}{M}\tr(T_2^*RRT_2)}{|1-x_0c^2a^2(t_2^*)^2|^2}=\dfrac{\frac{1}{M}\tr(RT_2T_2^*R)}{|1-x_0c^2t_2^2|^2}
\end{align*}
The numerators of both sides are equal and non zero, from what follows that the denominators are also equal, i.e.
\begin{align*}
|1-x_0c^2a^2(t_2^*)^2|=|1-x_0c^2t_2^2|
\end{align*}
We remark that if $w$ and $z$ satisfy  $|1-w|=|1-z|$ and $|w|=|z|$, then, either $w=z$, either $w=\bar{z}$. We use this remark
for $w = x_0c^2t_2^2$ and $z=x_0c^2a^2(t_2^*)^2$. If $w=z$, it holds that $a^2(t_2^*)^2=t_2^2\Rightarrow t_1^2=t_2^2$ and since $\im t_i\ge0$ we conclude $t_1=t_2$. If $w=\bar{z}$,
we have $a^2(t_2^*)^2=(t_2^*)^2$. If $t_2=0$ then it also holds that $t_1=0$. Otherwise, we have $a=\pm1$. If $a=1$, the condition $\im t_i\ge 0$, leads to
the conclusion that $t_1$ and $t_2$ are real and coincide. We finally consider the case $a=-1$. We recall $T_1=aT_2^* = - T_2^*$. Therefore, it holds that 
\begin{align*}
x_0I_M-\dfrac{x_0t_2^*}{1-x_0c^2(t_2^*)^2}R=-x_0I_M-\dfrac{x_0t_2^*}{1-x_0c^2(t_2^*)^2}R,
\end{align*}
which is impossible, since $x_0\ne0$. This completes the proof of Lemma~(\ref{le:t1=t2}).  $\blacksquare$
\\

Lemmas \ref{le:t1=t2} and \ref{le:t(zn)-bounded}, and their corresponding proofs imply the following result. 
\begin{proposition}
\label{eq:prop-t(x)}
For each $x > 0$, $\lim_{z \rightarrow x, z \in \mathbb{C}^{+}} t(z) = t(x)$ exists. Moreover, 
$1 - x (c t(x))^{2} \neq 0$, and matrix $(I + \frac{c t(x)}{1 - x (c t(x))^{2}} \, R)$ 
is invertible. Therefore, $\lim_{z  \rightarrow x, z \in \mathbb{C}^{+}} T(z) = T(x)$ 
where $T(x)$ represents matrix $T(x) = \left(-x(I + \frac{c t(x)}{1 - x (c t(x))^{2}} \, R)\right)^{-1}$. 
Moreover, $t(x)$ is solution of the equation 
\begin{equation}
\label{eq:t=x-real}
t(x) = \frac{1}{M} \mathrm{Tr}( R T(x))
\end{equation}
If $u(x)$ and $v(x)$ represent the terms defined by (\ref{eq:def-u}) and 
(\ref{eq:def-v}) for $z=x$, then it holds that
\begin{equation}
  \label{eq:1-u-positive-on-R-S}
  1 - u(x) > 0
\end{equation}
and 
\begin{equation}
  \label{eq:det-on-R-S}
  (1 - u(x))^{2} - x^{2} (v(x))^{2} \geq 0
\end{equation}
for each $x \neq 0$. Moreover, the inequality (\ref{eq:det-on-R-S}) is strict 
if $x \in \mathbb{R}^{+} - \mathcal{S}$. If moreover $\im (t(x)) > 0$, then, we have
\begin{equation}
\label{eq:det=0}
1 - u(x) - x v(x) = 0
\end{equation}
\end{proposition}
It just remains to justify  (\ref{eq:1-u-positive-on-R-S}), (\ref{eq:det-on-R-S}), and (\ref{eq:det=0}).  
As function $z \rightarrow t(z)$ is analytic
on $\mathbb{C} - \mathcal{S}$, $x \rightarrow t(x)$ is differentiable on $\mathbb{R}^{+} - \mathcal{S}$.
As $(t(x))^{'} > 0$ and $(xt(x))^{'} > 0$ hold on $\mathbb{R}^{+} - \mathcal{S}$, 
the arguments used in the context of Remark \ref{re:properties-u-v-z-real-negative} are also valid on  $\mathbb{R}^{+} - \mathcal{S}$, thus justifying there (\ref{eq:1-u-positive-on-R-S}) and the strict inequality in (\ref{eq:det-on-R-S}). 
$1 - u(x) \geq 0$ and inequality (\ref{eq:det-on-R-S}) also hold on $\mathcal{S}-\{0 \}$ by letting $z \rightarrow x$, $z \in \mathbb{C}^{+}$ in Proposition  \ref{re:properties-u-v-z-real-negative}. As $v(x) > 0$ for each $x \neq 0$, the strict inequality (\ref{eq:1-u-positive-on-R-S}) is a consequence of  (\ref{eq:det-on-R-S}). \\

In order to prove (\ref{eq:det=0}), we use the second component of (\ref{eq:linear-system-imt-imzt}), 
and remark that it implies that 
$$
\im (t(x)) = \left(u(x) + x v(x)\right) \, \im (t(x))
$$
Therefore, $\im (t(x)) > 0$ leads to (\ref{eq:det=0}).  $\blacksquare$ \\

We also add the following useful result which shows that the real part of $t(x)$ is negative for each 
$x > 0$. 
\begin{proposition}
\label{prop:extra-property}
For each $x \in \mathbb{R}^{+*}$, it holds that $\mathrm{Re}(t(x)) < 0$.
\end{proposition}
    {\bf Proof.}
    It is easily checked that
    \begin{equation}
\label{eq:linear-system-ret-rezt}
\left( \begin{array}{c} \mathrm{Re}(t(z)   \\  \mathrm{Re}(z t(z))  \end{array} \right) = 
\left( \begin{array}{cc} u(z) & - v(z) \\ - |z|^{2} v(z) & u(z) \end{array} \right) \; \left( \begin{array}{c}  \mathrm{Re}(t(z)   \\  \mathrm{Re}(z t(z)) \end{array} \right) + \left( \begin{array}{c} - \mathrm{Re}(z) \frac{1}{M} \mathrm{Tr}(R T(z) (T(z))^*) \\   - |z|^{2} \frac{1}{M} \mathrm{Tr}(R T(z) (T(z))^*\end{array} \right)
    \end{equation}
    for each $z \in \mathbb{C} - \mathcal{S}$. Moreover, as all the terms coming into play in (\ref{eq:linear-system-ret-rezt}) have a finite limit when $z \rightarrow x$ when $x \neq 0$, (\ref{eq:linear-system-ret-rezt}) remains 
valid on $\mathbb{R}^{*}$. For $z = x$, the first component of (\ref{eq:linear-system-ret-rezt}) leads to
    \begin{equation}
      \label{eq:expre-ret(x)-1}
      \mathrm{Re}(t(x)) ( 1 - u(x) + x v(x)) = - x \frac{1}{M} \mathrm{Tr}(RT(x)T(x)^{*})
    \end{equation}
    Proposition \ref{eq:prop-t(x)} implies that $1 - u(x) > 0$, when $x \in \mathbb{R}^{*}$. Therefore,
    $1 - u(x) + x v(x)$ is strictly positive as well, and it holds that
    \begin{equation}
      \label{eq:expre-ret(x)-2}
      \mathrm{Re}(t(x)) = - x \,  \frac{1}{1 - u(x) + x v(x)} \, \frac{1}{M} \mathrm{Tr}(R T(x) T(x)^{*})
    \end{equation}
    Therefore, $x > 0$ implies that $\mathrm{Re}(t(x)) < 0$ as expected.  $\blacksquare$ \\

We now study the behaviour of $t(z)$ when $z \rightarrow 0$. We first establish that
$\lim_{z \rightarrow 0, z \in \mathbb{C}^{+} \cup \mathbb{R}^{*}} |t(z)| = +\infty$, and then that
$\lim_{z \rightarrow 0, z \in \mathbb{C}^{+} \cup \mathbb{R}^{*}} z t(z) = 0$ if $c \leq 1$ and is strictly negative if $c > 1$. We recall that $t(x)$ for $x > 0$ is defined by $t(x) =  \lim_{z \rightarrow x, z \in \mathbb{C}^{+}} t(z)$.
For this, we establish various lemmas.
\begin{lemma}
  \label{le:t(z)-unbounded-0}
  It holds that $\lim_{z \rightarrow 0, z \in \mathbb{C}^{+} \cup \mathbb{R}^{*}} |t(z)| = +\infty$. 
\end{lemma}
We assume that the statement of the Lemma does not hold, i.e. that it exists a sequence of
elements of $\mathbb{C}^{+}  \cup \mathbb{R}^{*}$ $(z_n)_{n \geq 1}$  such that
  $\lim_{n \rightarrow +\infty} z_n = 0$ and  $t(z_n) \rightarrow t_0$. (\ref{t=}) and (\ref{eq:t=x-real}) 
imply that
\begin{equation}
  \label{eq:zt(z)=}
z_n t(z_n) = - \frac{1}{M} \sum_{l=1}^{\overline{M}} \frac{m_l \overline{\lambda}_l}{1 + \frac{c t(z_n) \overline{\lambda}_l}{1 - z_n (ct(z_n))^{2}}}
\end{equation}
$1 + \frac{c t(z_n) \overline{\lambda}_l}{1 - z_n (ct(z_n))^{2}}$ clearly converges towards $1 + c t_0 \overline{\lambda}_l$. As the left hand side of (\ref{eq:zt(z)=}) converges towards $0$, for each $l$, $1 +  c t_0 \overline{\lambda}_l$ cannot vanish. Therefore, matrix $I + c t_0 R$ is invertible,
and taking the limit of (\ref{eq:zt(z)=}) gives
$$
\frac{1}{M} \mathrm{Tr} R(I + c t_0 R)^{-1} = 0
$$
As $\im \frac{1}{M} \mathrm{Tr} R(I + c t_0 R)^{-1}$ cannot be zero if $t_0$ is not real,
$t_0$ must be real. We now use the observation that $|z_n| v(z_n) \leq 1$ for each $n$ (see
Lemma \ref{le:properties-u-v} if $z_n \in \mathbb{C}^{+} \cup \mathbb{R}^{+*}$, and Remark \ref{re:properties-u-v-z-real-negative} if $z_n \in \mathbb{R}^{-*}$). As $|1 - z_n (ct(z_n))^{2}|^{2} \rightarrow 1$,   
$|z_n| v(z_n)$ bounded implies that $|z_n| \frac{1}{M} \mathrm{Tr}(R T(z_n) R T(z_n)^{*})$
is bounded. It is easy to check that
$$
|z_n| \frac{1}{M} \mathrm{Tr}(R T(z_n) R T(z_n)^{*}) = \frac{1}{|z_n|}
\frac{1}{M} \mathrm{Tr}(R(I + c t_0 R)^{-1}R(I + c t_0 R)^{-1}) + \mathcal{O}(1)
$$
Therefore, the boundedness of $|z_n| \frac{1}{M} \mathrm{Tr}(R T(z_n) R T(z_n)^{*})$
implies that $\frac{1}{M} \mathrm{Tr}(R(I + c t_0 R)^{-1}R(I + c t_0 R)^{-1}) = 0$
which is of course impossible.  $\blacksquare$

\begin{lemma}
\label{le:zt(z)-bounded}
Consider a sequence $(z_n)_{n \geq 1}$ of elements of $\mathbb{C}^{+}  \cup \mathbb{R}^{*}$ such that 
$\lim_{n \rightarrow +\infty} z_n = 0$. Then, the set $(z_n t(z_n))_{n \geq 1}$ is bounded. 
\end{lemma}
{\bf Proof.} We assume that $(z_n t(z_n))_{n \geq 1}$ is not bounded. Therefore, one can extract
from $(z_n)_{n \geq 1}$ a subsequence, still denoted  $(z_n)_{n \geq 1}$, such that 
$\lim_{n \rightarrow +\infty} |z_n t(z_n)| = +\infty$. Then, 
$$
\frac{c t(z_n)}{1 - z_n (c t(z_n))^{2}} = \frac{1}{\frac{1}{c t(z_n)} - z_n t(z_n)} \rightarrow 0
$$
Therefore, 
$$
-\frac{1}{M} \mathrm{Tr} R \left(I + \frac{c t(z_n)}{1 - z_n (c t(z_n))^{2}} R\right)^{-1} \rightarrow 
-\frac{1}{M} \mathrm{Tr}R
$$
This is a contradiction because the above term coincides with $z_n t(z_n)$ which cannot converge towards 
a finite limit.  $\blacksquare$ \\

\begin{lemma}
\label{le:delta1=delta2}
Assume that $(z_{1,n})_{n \geq 1}$ and $(z_{2,n})_{n \geq 1}$ are sequences of elements of 
$\mathbb{C}^{+}  \cup \mathbb{R}^{*}$ such that $\lim_{n \rightarrow +\infty} z_{i,n} = 0$ and
$\lim_{n \rightarrow +\infty} z_{i,n} t(z_{i,n}) = \delta_i$ for $i=1,2$. 
Then, $\delta_1 = \delta_2$. 
\end{lemma}
We first remark that $|t(z_{i,n})| \rightarrow +\infty$ for $i=1,2$. Equation (\ref{t=}) implies immediately that  
\begin{equation}
\label{eq:zt=}
z t(z) = \left( z c t(z) - \frac{1}{c t(z)} \right) \frac{1}{M} \mathrm{Tr}R \left(R + \frac{1}{c t(z)} -  z c t(z)\right)^{-1} 
\end{equation}
As $\frac{1}{c t(z_{i,n})} \rightarrow 0$, $z_{i,n} c t(z_{i,n}) - \frac{1}{c t(z_{i,n})} \rightarrow c \delta_i$ 
for $i=1,2$. If $\delta_i \neq 0$, Eq. (\ref{eq:zt=}) thus implies that $c \frac{1}{M} \mathrm{Tr}R \left(R + \frac{1}{c t(z_{i,n})} -  z_{i,n} c t(z_{i,n})\right)^{-1}$ converges towards $1$, which implies that 
matrix  $R - c \delta_i I$ is invertible. Therefore, either $\delta_i = 0$, either $\delta_i$ is a solution of the equation
\begin{equation}
  \label{eq:deltai-non-zero}
1 = c \frac{1}{M} \mathrm{Tr}R(R - c \delta_i I)^{-1}
\end{equation}
or equivalently, $\delta_i$ verifies
\begin{equation}
  \label{eq:deltai}
\delta_i = c \delta_i \frac{1}{M} \mathrm{Tr}R(R - c \delta_i I)^{-1}
\end{equation}
We note that the solutions of this equation are real, so that $\delta_i \in \mathbb{R}$ for $i=1,2$.  
Eq. (\ref{eq:linear-system-z1-z2}) leads to 
\begin{multline*}
z_{1,n} t(z_{1,n}) - z_{2,n} t(z_{2,n}) = z_{1,n} z_{2,n} v_0(z_{1,n}, z_{2,n}) ( t(z_{1,n}) -  t(z_{2,n}))\\
 + u_0(z_{1,n}, z_{2,n}) (z_{1,n} t(z_{1,n}) - z_{2,n} t(z_{2,n}))
\end{multline*}
It is straightforward to check that $z_{1,n} z_{2,n} v_0(z_{1,n}, z_{2,n}) ( t(z_{1,n}) -  t(z_{2,n})) \rightarrow 0$ and that $u_0(z_{1,n}, z_{2,n}) \rightarrow u_0(0,0) = c \frac{1}{M} \mathrm{Tr}R(R-c \delta_1 I)^{-1} R (R - c \delta_2 I)^{-1}$. Therefore, we obtain that 
\begin{equation}
  \label{eq:delta1-delta2}
\delta_1 - \delta_2 = u_0(0,0) (\delta_1 - \delta_2)
\end{equation}
We recall that $|u_0(z_{1,n}, z_{2,n})| \leq \sqrt{u(z_{1,n})} \sqrt{u(z_{2,n})} \leq 1$. Moreover, 
we observe that $u(z_{i,n}) \rightarrow u_i(0) = c \frac{1}{M}  \mathrm{Tr}R(R-c \delta_i I)^{-1} R (R - c \delta_i I)^{-1}$ and that $0 < u_i(0) \leq 1$. The Schwartz inequality leads to 
\begin{equation}
\label{eq:schwartz-zt-0}
|u_0(0,0)| \leq \sqrt{u_1(0)} \sqrt{u_2(0)} \leq 1
\end{equation}
If the Schwartz inequality (\ref{eq:schwartz-zt-0}) is strict, $|u_0(0,0)| < 1$, and $\delta_1 = \delta_2$. 
We now assume that $u_0(0,0) = \sqrt{u_1(0)} \sqrt{u_2(0)} = 1$. This implies that 
$$
R - c \delta_1 I = \kappa (R - c \delta_2 I)
$$
for some real constant $\kappa$, or equivalently, $\overline{\lambda}_l - c \delta_1 = \kappa(\overline{\lambda}_l - c \delta_2)$ for each $l=1, \ldots, \overline{M}$. If $R$ is not a multiple of $I$, 
$\kappa$ must be equal to 1, since otherwise, we would have $\overline{\lambda}_l = \overline{\lambda}_{l'}$ 
for each $l,l'$. $\kappa =1$ implies immediately that $\delta_1 = \delta_2$. We finally consider the case 
where $R = \sigma^{2} I$. Then, (\ref{eq:deltai}) implies that $\delta_i$ is solution of $\delta_i \frac{\sigma^2 c}{\sigma^{2} - c \delta_i} = \delta_i$, i.e. $\delta_i = 0$ or
\begin{equation}
  \label{eq:deltai-R=I}
\delta_i =   \sigma^{2} \left( \frac{1}{c} - 1 \right)
\end{equation}
We now check that $\delta_1 = 0, \delta_2 = \sigma^{2} \left( \frac{1}{c} - 1 \right)$ or
$\delta_2 = 0, \delta_1 = \sigma^{2} \left( \frac{1}{c} - 1 \right)$ is impossible.
If this holds, $u_1(0)$ and $u_2(0)$ cannot be both equal to 1, and $|u_0(0,0)| < 1$. Therefore, (\ref{eq:delta1-delta2})
leads to a contradiction, and $\delta_1 = \delta_2$ is equal either to $0$, either to  $\sigma^{2} \left( \frac{1}{c} - 1 \right)$.  $\blacksquare$ \\

Lemmas \ref{le:zt(z)-bounded} and \ref{le:delta1=delta2} imply the following corollary. 
\begin{corollary}
\label{coro:behaviour-zt(z)-c-leq-1}
If $c \leq 1$, it holds that 
\begin{equation}
\label{eq:delta=0-if-c<1}
\lim_{z \rightarrow 0, z \in \mathbb{C}^{+} \cup \mathbb{R}^{*}} z t(z) = 0
\end{equation}
and that 
\begin{equation}
\label{eq:no-mass-0}
\mu(\{0 \}) = 0
\end{equation}
\end{corollary}
{\bf Proof.} Lemmas \ref{le:zt(z)-bounded} and \ref{le:delta1=delta2} lead to the conclusion 
that $\lim_{z \rightarrow 0, z \in \mathbb{C}^{+} \cup \mathbb{R}^{*}} z t(z) = \delta$ 
where $\delta$ is either equal
to $0$, either coincides with a solution of the equation (\ref{eq:deltai}). In order to precise this,
we remark that $t(x) > 0$ if $x < 0$ implies that $\delta \leq 0$. 
Therefore, $\delta$ coincides with a non positive solution of equation (\ref{eq:deltai}).
If $c \leq 1$, it is clear that (\ref{eq:deltai}) has no strictly negative solutions. Therefore,
(\ref{eq:delta=0-if-c<1}) is established. (\ref{eq:no-mass-0}) is a direct consequence of the 
identity
$$
\mu(\{0\}) = \lim_{z \rightarrow 0, z \in \mathbb{C}^{+} \cup \mathbb{R}^{*}} - z t(z) 
$$
 $\blacksquare$

In order to address the case where $c >1$ and to precise the behaviour of $\im (t(z))$ when
$z \rightarrow 0, z \in \mathbb{C}^{+} \cup \mathbb{R}^{*}$ if $c \leq 1$, we have to evaluate 
$z (t(z))^{2}$ when $z \rightarrow 0$. The following Lemma holds.
\begin{lemma}
  \label{le:z(t(z))^{2}}
  \begin{itemize}
  \item  If $c=1$, it holds that
    $\lim_{z \rightarrow \mathbb{C}^{+} \cup \mathbb{R}^{*}} |z (t(z))^{2}| = +\infty$.
  \item If $c < 1$, 
\begin{equation}
\label{eq:z(t(z))2-0}
\lim_{z \rightarrow \mathbb{C}^{+} \cup \mathbb{R}^{*}} z (t(z))^{2} = - \frac{1}{c(1-c)}
\end{equation}
  \item If $c > 1$, the assumption $\lim_{z \rightarrow 0, z \in \mathbb{C}^{+} \cup \mathbb{R}^{*}} z t(z)  = \delta = 0$ implies that
    $\lim_{z \rightarrow \mathbb{C}^{+} \cup \mathbb{R}^{*}} z (t(z))^{2} = - \frac{1}{c(1-c)}$, a contradiction
    because the above limit is necessarily negative. Hence, $\delta$ is non zero and coincides with the 
strictly negative solution of Eq. (\ref{eq:deltai}), and $\mu(\{0\}) = - \delta$.
    \end{itemize}
\end{lemma}
    {\bf Proof.} (\ref{t=}) implies that
    \begin{equation}
      \label{eq:z(t(z)^{2}=}
    z (t(z))^{2} = - \frac{1}{M} \mathrm{Tr}R \left( \frac{I}{t(z)} + \frac{c}{1 - z (ct(z))^{2}} R \right)^{-1}
    \end{equation}
    We assume in the course of this proof that $\delta = 0$ (if $c \leq 1$, this property holds). We first establish the
    first item of Lemma \ref{le:z(t(z))^{2}}. We assume that $c=1$ and that
    there exists a sequence $(z_n)_{n \in \mathbb{C}^{+} \cup \mathbb{R}^{*}}$ such that $z_n \rightarrow 0$ and $z_n t(z_n)^{2} \rightarrow \alpha$. 
    As $|t(z_n)| \rightarrow +\infty$, (\ref{eq:z(t(z)^{2}=}) leads to $\alpha = \alpha - 1$, a contradiction. Therefore, if $c=1$,
    $\lim_{z \rightarrow 0,  \rightarrow \mathbb{C}^{+} \cup \mathbb{R}^{*}} |z t(z)^{2}| = +\infty$ as expected.

    We now establish the 2 last items. For this, we establish that if $c \neq 1$, then,  $|z t(z)^{2}|$ is 
bounded when
    $z \in \mathbb{C}^{+} \cup \mathbb{R}^{*}$ and $z$ is close from $0$. For this, 
    we assume the existence of a sequence $(z_n)_{n \geq 1}$ of elements of  $\mathbb{C}^{+} \cup \mathbb{R}^{*}$
    such that $z_n \rightarrow 0$ and $|z_n t(z_n)^{2}| \rightarrow +\infty$. Then, it holds that
    $$
    1 =  - \frac{1}{M} \mathrm{Tr}R \left(z_nt(z_n) I + \frac{c z_n  t(z_n)^{2}}{1 - z_n (ct(z_n))^{2}} R \right)^{-1}
    $$
    As $|z_n t(z_n)^{2}| \rightarrow +\infty$, $\frac{c z_n  t(z_n)^{2}}{1 - z_n (ct(z_n))^{2}} \rightarrow -\frac{1}{c}$. Condition
    $z_n t(z_n) \rightarrow 0$ thus implies that $c = 1$, a contradiction. Using again  (\ref{eq:z(t(z)^{2}=}),
    we obtain immediately that if $z_n (t(z_n))^{2} \rightarrow \alpha$, then $\alpha = - \frac{1}{c(c-1)}$. As  $|z t(z)^{2}|$ remains bounded when $z \in \mathbb{C}^{+} \cup \mathbb{R}^{*}$ is close from $0$,
    this implies that $\lim_{z \rightarrow 0, z \in \mathbb{C}^{+} \cup \mathbb{R}^{*}} z (t(z))^{2} = - \frac{1}{c(1-c)}$ as expected. Taking $z \in \mathbb{R}^{-*}$
    leads to the conclusion that the above limit is negative. When $c > 1$, this is a contradiction 
because $- \frac{1}{c(1-c)}$ is positive. Therefore, if $c > 1$, $\delta$, the limit of $z t(z)$, cannot be equal to $0$. Hence, $\delta$ coincides with the strictly negative solution of (\ref{eq:deltai}) and 
$\mu(\{0\}) = - \delta > 0$. This completes the proof of the Lemma.  $\blacksquare$ \\

Putting all the pieces together, we obtain the following characterization of $\mu_N$ when $c \leq 1$. 
\begin{theorem}
  \label{th:characterization-mu-c<1}
  The density $f_N(x)$ of $\mu_N$ w.r.t. the Lebesgue measure is a continuous function on $\mathbb{R}^{+*}$,
  and is given by $f_N(x) = \frac{1}{\pi} \im (t_N(x))$ for each $x > 0$. If $c_N \leq 1$,
  $\mu_N$ is absolutely continuous, and if $c_N > 1$, then $d \mu_N(x) = f_N(x) dx + \mu_N(\{0\}) \delta_0$.
  $0 \in \mathcal{S}_N$, and the interior  $\mathcal{S}_N^{\circ}$ of $\mathcal{S}_N$ is given by 
\begin{equation}
\label{eq:initial-characterization-support-mu}
\mathcal{S}_N^{\circ} = \{ x \in \mathbb{R}^{+}, \im (t(x)) > 0 \}
\end{equation}
If moreover $c_N < 1$, it holds that
\begin{equation}
\label{eq:behaviour-f-0}
f_N(x) \simeq  \frac{1}{\pi} \, \frac{1}{\sqrt{x \, c_N (1-c_N)}}
\end{equation}
when $x \rightarrow 0^{+}$, while if $c_N =1$,
\begin{equation}
\label{eq:behaviour-f-0-c=1}
f_N(x) \simeq  \frac{1}{\pi} \, \frac{\sqrt{3}}{2}  \left( \frac{1}{M}\mathrm{Tr}\, R^{-1}\right)^{-1/3} \frac{1}{x^{2/3}}
\end{equation}
\end{theorem}
    {\bf Proof.} $t(z)$ is not analytic in a neighbourhood of $0$; hence, 
$0 \in \mathcal{S}$. As $\lim_{z \rightarrow x, z \in \mathbb{C}^{+}} t(z) = t(x)$ exists for $x \neq 0$, Theorem 2.1 of
    \cite{silverstein-choi-1995} implies that if $\mathcal{A} \in \mathbb{R}^{+*}$ is a Borel set of zero Lebesgue measure,
    then $\mu(\mathcal{A}) = \int_{\mathcal{A}} f(x) dx = 0$. The continuity of $f$ on $\mathbb{R}^{+*}$ is a also a consequence
    of  \cite{silverstein-choi-1995}. \\

We now prove (\ref{eq:behaviour-f-0}). For this, we remark that (\ref{eq:z(t(z))2-0}) implies that
\begin{equation}
\label{eq:x(t(x))2-0}
\lim_{x \rightarrow 0, x > 0} x (t(x))^{2} = - \frac{1}{c(1-c)}
\end{equation}
As $\im (t(x)) \geq 0$ for each $x \neq 0$, (\ref{eq:x(t(x))2-0}) 
implies that $t(x) \simeq \frac{i}{\sqrt{x} \, \sqrt{c(1-c}}$ when $x \rightarrow 0^{+}$, or equivalently
that $\frac{1}{\pi} \im (t(x)) \simeq  \frac{1}{\pi} \, \frac{1}{\sqrt{x \, c (1-c)}}$. \\

It remains to establish (\ref{eq:behaviour-f-0-c=1}). For this, we first prove that
\begin{equation}
  \label{eq:behaviour-x2t3}
\lim_{x \rightarrow 0, x > 0} x^{2} (t(x))^{3} = \left(\frac{1}{M} \mathrm{Tr}\, R_N^{-1}\right)^{-1}
\end{equation}
For this, we write (\ref{eq:t=x-real}) as
\begin{equation}
\label{eq:new-t=x-real}
\frac{1}{M} \mathrm{Tr} R \left( -x t(x) I + \frac{1}{1-\frac{1}{x (t(x))^{2}}} R \right)^{-1} = 1
\end{equation}
As $c = 1$, $x t(x) \rightarrow 0$ and $|x (t(x))^{2}| \rightarrow +\infty$ when $x \rightarrow 0, x > 0$.
The left hand side of (\ref{eq:new-t=x-real}) can be expanded as
\begin{multline*}
\frac{1}{M} \mathrm{Tr} R \left( -x t(x) I + \frac{1}{1-\frac{1}{x (t(x))^{2}}} R \right)^{-1}  = 1 - \frac{1}{x (t(x))^{2}} \\
+ 
\frac{1}{M} \mathrm{Tr}\, R^{-1}  \, x t(x) + x t(x) \epsilon_1(x) + \frac{1}{x (t(x))^{2}} \epsilon_2(x)
\end{multline*}
where $\epsilon_1(x)$ and $\epsilon_2(x)$ converge towards $0$ when $x \rightarrow 0, x > 0$. Therefore,  (\ref{eq:new-t=x-real}) implies that
$$
 \frac{1}{M} \mathrm{Tr}\, R^{-1} \, x t(x) - \frac{1}{x (t(x))^{2}}   = x t(x) \tilde{\epsilon}_1(x) + \frac{1}{x (t(x))^{2}} \tilde{\epsilon}_2(x)
$$
where  $\tilde{\epsilon}_1(x)$ and $\tilde{\epsilon}_2(x)$ converge towards $0$ when $x \rightarrow 0, x > 0$. This leads immediately to
(\ref{eq:behaviour-x2t3}). As function $x \rightarrow  x^{2} (t(x))^{3}$ is continuous on $\mathbb{R}^{+*}$, it holds that
$$
\lim_{x \rightarrow 0, x > 0} x^{2/3} t(x) = e^{2i k \pi/3}  \left(\frac{1}{M} \mathrm{Tr}\, R^{-1} \right)^{-1/3} 
$$
where $k$ is equal to $0,1$ or $2$. If $k=0$, the real part of $t(x)$ must be positive if $x$ is close 
enough from $0$. Lemma \ref{prop:extra-property} thus leads to a contradiction.  If $k=2$, $\im (t(x)) < 0$ for $x$ small
enough, a contradiction as well. Hence, $k$ is equal to $1$. Therefore,
\begin{equation}
  \label{eq:behav-t-near-0-c=1}
\lim_{x \rightarrow 0, x > 0} x^{2/3} \im (t(x)) = \sin 2 \pi/3 \,   \left(\frac{1}{M} \mathrm{Tr}\, R^{-1} \right)^{-1/3} 
\end{equation}
This completes the proof of (\ref{eq:behaviour-f-0-c=1}).  $\blacksquare$ \\

We now show that function $x \rightarrow t(x)$ and 
$x \rightarrow f(x)$ possess a power series expansion in a neighbourhood of each point of 
$\mathcal{S}_N^{\circ}$. More precisely: 
\begin{proposition}
\label{prop:t(x)-analytic}
If $x_0 > 0$ and $\im (t(x_0)) > 0$, then, $t$ and $f$ can be expanded as 
$$
t(x) = \sum_{k=0}^{+\infty} a_k (x-x_0)^{k}, f(x) = \sum_{k=0}^{+\infty} b_k (x-x_0)^{k}
$$
when $|x-x_0|$ is small enough. 
\end{proposition}
As in \cite{silverstein-choi-1995} and \cite{dozier-silverstein-2004}, the proof is based on the 
holomorphic implicit function theorem (see \cite{cartan}). We denote $t(x_0)$ by $t_0$. Then, Eq. (\ref{eq:t=x-real}) at point $x_0$ can be written as $h(x_0,t_0) = 0$ where function $h(z,t)$ is defined 
by 
$$
h(z,t) = t - \frac{1}{M} \mathrm{Tr}\left( R \left(-z(I + \frac{c t}{1 - z (c t)^{2}} \, R) \right)^{-1} \right)
$$
As $x_0 > 0$ and $\im (t_0) > 0$, function $(z,t) \rightarrow h(z,t)$ is holomorphic 
in a neighbourhood of $(x_0,t_0)$. It is easy to check that 
\begin{equation}
\label{eq:expre-partial-derivative-h}
\left(\frac{\partial h}{\partial t}\right)_{x_0,t_0} = 1 - u_{0}(x_0,x_0) - x_0^{2} v_{0}(x_0,x_0)
\end{equation}
where we recall that functions $u_0$ and $v_0$ are given by (\ref{eq:def-u_0(z1,z2)}) and
(\ref{eq:def-v_0(z1,z2)}). Following the proof of Lemma \ref{le:t1=t2}, we obtain immediately 
that $ 1 - u_{0}(x_0,x_0) - x_0^{2} v_{0}(x_0,x_0) = 0$ implies that $T(x_0) = a T(x_0)^{*}$, 
and that $t_0 = a t_0^{*}$ for some $a \in \mathbb{C}$. The arguments of the above proof then 
lead to the conclusion that $t_0 = t_0^{*}$, a contradiction because  $\im (t(x_0)) > 0$. 
Hence, $\left(\frac{\partial h}{\partial t}\right)_{x_0,t_0} \neq 0$. The holomorphic implicit 
function theorem thus implies that it exists a function $z \rightarrow \tilde{t}(z)$, holomorphic 
in a neigbourhood $N$ of $x_0$, verifying $\tilde{t}(x_0) = t_0$ and 
$h(z,\tilde{t}(z)) = 0$ for each $z \in N$. Moreover, condition $\im (t_0) = 
\im (\tilde{t}(x_0)) > 0$ implies that $\im (\tilde{t}(z)) > 0$ and 
$\im (z \tilde{t}(z)) > 0$ if $|z - x_0| < \epsilon$ for $\epsilon$ small enough. 
Therefore, if $z \in \mathbb{C}^{+}$ and $|z - x_0| < \epsilon$, it must hold that 
$\tilde{t}(z) = t(z)$ (see Proposition \ref{prop:existence-uniqueness}). Hence, 
$t(x) = \lim_{z \rightarrow x, z \in \mathbb{C}^{+}} t(z)$ must coincide with $\tilde{t}(x)$ when 
if $|x - x_0| < \epsilon$. As 
$\tilde{t}(z)$ is holomorphic in a neighbourhood of $x_0$, function $x \rightarrow t(x)$ 
can be expanded as 
$$
t(x) = \sum_{k=0}^{+\infty} a_k (x-x_0)^{k}
$$
when $|x - x_0| < \epsilon$. This immediately implies that $f$ possesses a power series expansion 
in the interval $(x_0 - \epsilon, x_0 + \epsilon)$.  $\blacksquare$ \\

We finally use the above results in order the study measure $\nu_N$ associated to the Stieltjes transform
\begin{equation}
  \label{eq:def-s}
  t_{N,\nu}(z) = \frac{1}{M} \mathrm{Tr}T_N(z)
\end{equation}
As $\nu_N$ and $\mu_N$ are absolutely continuous one with respect each other, $d\nu_N(x)$
can also be written as $d\nu_N(x) = g_N(x) dx + \nu_N(\{0\}) \delta_0$. 
Using the identity
$$
\frac{1}{M} \mathrm{Tr} \left[ -z\left( I + \frac{c t(z)}{1 - z (c t(z))^{2}} R \right) \right] \, T(z) = 1
$$
we obtain immediately that
\begin{equation}
  \label{eq:expre-s}
  t_{\nu}(z) = -\frac{1}{z} -  \frac{c (t(z))^{2}}{1 - z (c t(z))^{2}}
\end{equation}
If $x > 0$, $t_{\nu}(x) = \lim_{z \rightarrow x, z \in \mathbb{C}^{+}}$ exists, and is given by the righthandside of
(\ref{eq:expre-s}) when $z=x$. Hence, for $x > 0$, $g(x) = \frac{1}{\pi} \im (t_{\nu}(x))$, i.e.
\begin{equation}
  \label{eq:expre-g}
  g(x) = - \frac{1}{\pi} \, \frac{ c \, \im ((t(x))^{2})}{|1 - x(ct(x))^{2}|^{2}}
\end{equation}
If $c > 1$, $|z t(z)^{2}| \rightarrow +\infty$ if $z \rightarrow 0$. (\ref{eq:expre-s}) thus implies that
$\nu_N(\{0 \}) = \lim_{z \rightarrow 0} - z t_{\nu}(z)$ coincides with $1 - \frac{1}{c}$, which, of course, is not surprising. 
We now evaluate the behaviour of $g$ when $x \rightarrow 0, x > 0$ and $c \leq 1$. 
\begin{proposition}
  \label{prop:behaviour-g-0}
  If $c < 1$, it holds that
  \begin{equation}
    \label{eq:g-near-0-c<1}
    g(x) \simeq_{x \rightarrow 0} \frac{1}{\pi} \, \frac{1}{ \sqrt{c \, (1 - c)}} \, \frac{1}{M} \mathrm{Tr}(R^{-1}) \, \frac{1}{\sqrt{x}}
  \end{equation}
  while if $c=1$, it holds that
  \begin{equation}
    \label{eq:g-near-0-c=1}
    g(x) \simeq_{x \rightarrow 0} \frac{1}{\pi} \,  \frac{\sqrt{3}}{2} \, \left( \frac{1}{M} \mathrm{Tr}(R^{-1}) \right)^{2/3} \frac{1}{x^{2/3}}
  \end{equation}
  \end{proposition}
    {\bf Proof.} Using Eq. (\ref{eq:z(t(z)^{2}=}), we obtain after some algebra that
    $$
    z (t(z))^{2} + \frac{1}{c(1-c)} \simeq_{z \rightarrow 0}  \, \frac{1}{M} \mathrm{Tr}R^{-1} \, \frac{1}{c^{2} (1 -c)^{3}} \frac{1}{t(z)}
    $$
    As $t(x) \simeq_{x \rightarrow 0, x > 0} \frac{i}{\sqrt{x} \sqrt{c(1-c)}}$, we get that
    $$
    \im ((t(x))^{2}) \simeq - i \,  \frac{1}{M} \mathrm{Tr}R^{-1} \, \frac{1}{1-c} \, \frac{1}{(c(1-c))^{3/2}} \, \frac{1}{\sqrt{x}}
    $$
    Therefore, (\ref{eq:expre-g}) immediately leads to (\ref{eq:g-near-0-c<1}). (\ref{eq:g-near-0-c=1}) is an immediate
    consequence of (\ref{eq:behav-t-near-0-c=1}).  $\blacksquare$ \\

Proposition \ref{prop:behaviour-g-0} means in practice that if $c_N \leq 1$, a number of eigenvalues of matrix $W_{f,N} W_{p,N}^{*}  W_{p,N} W_{f,N}^{*}$ are close from
$0$. Moreover, the rate of convergence of $g_N$ towards $+\infty$ is higher if $c_N = 1$, showing that 
in this case, the proportion of eigenvalues close to $0$ is even larger than if $c_N < 1$.  \\

    We finally mention that $t_{\nu}(x)$ and $g(x)$ possess a power expansion around eachpoint $x_0 \in \mathcal{S}^{\circ}$. This is an obvious
    consequence of Proposition \ref{prop:t(x)-analytic} and of the above expressions of $s_{\nu}(x)$ and of $g(x)$ in terms of $t(x)$.

\subsection{Characterization of $\mathcal{S}_N$.}

We denote by $w_N(z)$ the function defined by 
\begin{equation}
\label{eq:def-w}
w_N(z) = -\frac{(1-z (c_Nt_N(z))^{2})}{c_Nt_N(z)} = z c_N t_N(z) - \frac{1}{c_N t_N(z)}
\end{equation}
It is clear that $w$ is analytic on $\mathbb{C} - \mathcal{S}$, that $\im (w(z)) > 0$ if $z \in \mathbb{C}^{+}$, 
that $w(x) = \lim_{z \rightarrow x, z \in \mathbb{C}^{+}} w(z)$ exists for each $x \in \mathbb{R}^*$, and that the limit still exists
if $x=0$. If we denote this limit by $w(0)$, then, it holds that $w(0)=0$ if $c \leq 1$ and that $w(0) = c \delta$
if $c >1$, where we recall that $\delta$ is defined as the solution of (\ref{eq:deltai-non-zero})· Moreover, $w(x)$ is real if and only if $t(x)$ is real. Therefore, the interior $\mathcal{S}^{o}$ of $\mathcal{S}$ is also given by 
\begin{equation}
\label{eq:interior-support-w} 
\mathcal{S}^{o} = \{ x \in \mathbb{R}^{+}, \im (w(x)) > 0 \}
\end{equation}
Moreover, as $t(x)^{\prime}$ and $(x t(x))^{\prime}$ are strictly positive  if $x \in \mathbb{R} - \mathcal{S}$,  the derivative $w^{\prime}(x)$ of $w(x)$ w.r.t. $x$ is also strictly positive on $ \mathbb{R} - \mathcal{S}$. $t(z)$ can be expressed in terms of $w(z)$ as 
\begin{equation}
\label{eq:relation-t-w-1}
t(z) = \frac{1}{z} \, w(z) \, \frac{1}{M} \mathrm{Tr} R \left( R - w(z) I \right)^{-1}
\end{equation}
(\ref{eq:def-w}) implies that 
\begin{equation}
\label{eq:relation-t-w-2}
1 + c t(z) w(z) - z (c t(z))^{2} = 0
\end{equation}
Plugging (\ref{eq:relation-t-w-1}) into (\ref{eq:relation-t-w-2}), we obtain immediately that $w_N(z)$ verifies the equation 
\begin{equation}
\label{eq:inversion}
\phi_N(w_N(z)) = z
\end{equation}
where $\phi_N(w)$ is defined by 
\begin{equation}
\label{eq:def-phi}
\phi_N(w) = c_N w^{2} \, \frac{1}{M} \mathrm{Tr} R_N \left( R_N - w I \right)^{-1} \, \left( c_N \, \frac{1}{M} \mathrm{Tr} R_N \left( R_N - w I \right)^{-1} - 1 \right)
\end{equation}
Observe that (\ref{eq:inversion}) holds not only on $\mathbb{C} - \mathcal{S}$, but also for each $x \in \mathcal{S}$. Therefore, 
it holds that $\phi(w(x)) = x$ for each $x \in \mathbb{R}$. For each $x \in \mathbb{R} - \mathcal{S}$, it thus holds that 
$\phi^{'}(w(x)) \, w^{'}(x) = 1$. Therefore, as $w^{'}(x) > 0$ if $x \in \mathbb{R} - \mathcal{S}$, $w(x)$ satisfies $\phi^{'}(w(x)) > 0$ for each  $x \in \mathbb{R} - \mathcal{S}$. This implies that if 
$x \in \mathbb{R} - \mathcal{S}$, then $w(x)$ is a real solution of the polynomial equation $\phi(w) = x$ for which 
$\phi^{'}(w) > 0$. Moroever, Proposition \ref{prop:extra-property} implies that if $x \in \mathbb{R}^{+} - \mathcal{S}$, then, $t(x) = \mathrm{Re}(t(x))$ is strictly negative. Eq. (\ref{eq:relation-t-w-1}) for $z=x$ thus leads to
the conclusion that if $x > 0$ does not belong to $\mathcal{S}$, then $w(x)$ also verifies
$w(x) \frac{1}{M} \mathrm{Tr} R \left( R - w(x) I \right)^{-1} < 0$. If $x < 0$, then, $t(x)$ is this time 
strictly positive and $w(x)$ still verifies $w(x) \frac{1}{M} \mathrm{Tr} R \left( R - w(x) I \right)^{-1} < 0$. This discussion leads to the following Proposition.
\begin{proposition}
  \label{prop:characterization-w(x)}
  If $x \in \mathbb{R} - \mathcal{S}$, then $w(x)$ verifies the following properties:
  \begin{equation}
\label{eq:proprietes-w(x)} 
  \phi(w(x))  = x, \; \phi^{'}(w(x)) >  0, \; w(x) \, \frac{1}{M} \mathrm{Tr} R \left( R - w(x) I \right)^{-1} < 0
\end{equation}
\end{proposition}
As shown below, if $x \in \mathbb{R} - \mathcal{S}$, the properties (\ref{eq:proprietes-w(x)}) characterize $w(x)$ among the set of all solutions of the equation $\phi(w) = x$ and allow to identify the support 
as the subset of $\mathbb{R}^{+}$ for which the equation $\phi(w) = x$ has no real solution satisfying 
the conditions (\ref{eq:proprietes-w(x)}). These results follow directly from an elementary study 
of function $w \rightarrow \phi(w)$. \\

We first consider the case $c \leq 1$, and identify the values of $x > 0$ for which the equation 
$\phi(w(x))=x$ has a real solution verifying (\ref{eq:proprietes-w(x)}), and those for which 
such a solution does not exist.  It is easily seen that if $x > 0$, all the real solutions of 
the equation $\phi(w)=x$ are strictly positive. Therefore, the third condition in (\ref{eq:proprietes-w(x)}) 
is equivalent to $ \frac{1}{M} \mathrm{Tr} R \left( R - w(x) I \right)^{-1} < 0$. 
We denote $\omega_{1,N} < \omega_{2,N} < \ldots < \omega_{\overline{M},N}$ the 
(necessarily real) $\overline{M}$
roots of $\frac{1}{M} \mathrm{Tr}R_N(R_N - w I)^{-1} = \frac{1}{c_N}$ and by $\mu_{1,N} < \mu_{2,N} < \ldots < \mu_{\overline{M}-1,N}$ the roots of
$\frac{1}{M} \mathrm{Tr}R_N(R_N - w I)^{-1} = 0$.
As $c \leq 1$, it is easily seen that $\omega_1 \geq 0$, and that
$\omega_{1} < \overline{\lambda}_{\overline{M}} < \mu_{1} < \omega_{2} < \overline{\lambda}_{\overline{M}-1} < \ldots < \mu_{\overline{M}-1} < \omega_{\overline{M}} < \overline{\lambda}_1$. It is clear
$\frac{1}{M} \mathrm{Tr}R(R-wI)^{-1} < 0$ if and only if
$w \in (\overline{\lambda}_{\overline{M}},\mu_{1}) \cup \ldots \cup  (\overline{\lambda}_2, \mu_{\overline{M}-1})\cup (\overline{\lambda}_1, +\infty)$. \\

For $x > 0$, the equation $\phi(w) = x$ is easily seen to be a polynomial equation of degree $2\overline{M}+1$. Therefore, $\phi(w) = x$ has $2\overline{M}+1$ solutions. For each $x > 0$, this equation has at least $2\overline{M} - 1$ real solutions that cannot coincide with $w(x)$ if  $x \in (\mathcal{S}^{\circ})^{c}$:
\begin{itemize}
\item $\overline{M}$ solutions belong to $]\omega_1, \overline{\lambda}_{\overline{M}}[, \ldots, ]\omega_{\overline{M}}, \overline{\lambda}_{1}[$.
    None of these solutions may correspond to $w(x)$ if $x \in (\mathcal{S}^{\circ})^{c}$ because $\frac{1}{M} \mathrm{Tr}R(R-wI)^{-1} > 0$ at these points.
  \item On each interval  $]\overline{\lambda}_{\overline{M}}, \mu_{1}[, \ldots, ]\overline{\lambda}_2, \mu_{\overline{M}-1}[$, the equation
    $\phi(w) = x$ has a real solution at which $\phi^{\prime}$ is negative. Therefore, $\phi(w) = x$ has $\overline{M}-1$ extra real
    solutions that are not equal to $w(x)$ if  $x \in (\mathcal{S}^{\circ})^{c}$.
\end{itemize}
As $\phi_N(w) \rightarrow +\infty$ if $w \rightarrow \overline{\lambda}_{1,N}, w > \overline{\lambda}_{1,N}$ and that
$\phi_N(w) \rightarrow +\infty$ if $w \rightarrow +\infty$, it exists at least a point in $]\overline{\lambda}_{1,N}, +\infty[$ at which $\phi_N^{'}$ vanishes. This point is moreover unique because otherwise, $\phi_N(w) = x$ would have more than $2 \overline{M} + 1$ solutions for certain values of $x$. We denote by  $w_{+,N}$ this point, and remark that if $x > x_{+,N} = \phi_N(w_{+,N})$,  $\phi_N(w) = x$ has $2 \overline{M} + 1$
    real solutions: the $2 \overline{M} - 1$ solutions that were introduced below, and 2 extra solutions that belong to
    $]\overline{\lambda}_1, w_{+}[$ and $]w_{+}, +\infty[$ respectively. Therefore, $w(x)$ is real, and it is easily seen that $w(x)$ coincides with the solution that belongs to $]w_{+}, +\infty[$. This implies that $]x_{+}, +\infty[  \subset \mathbb{R} - \mathcal{S}$. \\

If $\phi^{'}(w)$ does not vanish on  $]\overline{\lambda}_{\overline{M}}, \mu_{1}[ \cup \ldots \cup]\overline{\lambda}_2, \mu_{\overline{M}-1}[$, 
for each $x \in ]0, x_{+}[$, $\phi$ is decreasing on these intervals. Therefore, none of the real solutions of $\phi(w) = x$ match with the properties of 
$w(x)$ when $x \in \mathbb{R}^{+} - \mathcal{S}$. Therefore, $w(x)$ must be a complex number: $\phi(w) = x$ 
has thus $2 \overline{M} - 1$ real solutions, and a pair of complex conjugate roots: $w(x)$ is the positive 
imaginary part solution. In this case, $x \in \mathcal{S}^{\circ}$, and the support $\mathcal{S}$ coincides with 
$[0, x_{+}]$. \\

\begin{figure}[ht!]
	\centering
	\par
        \includegraphics[scale=0.5]{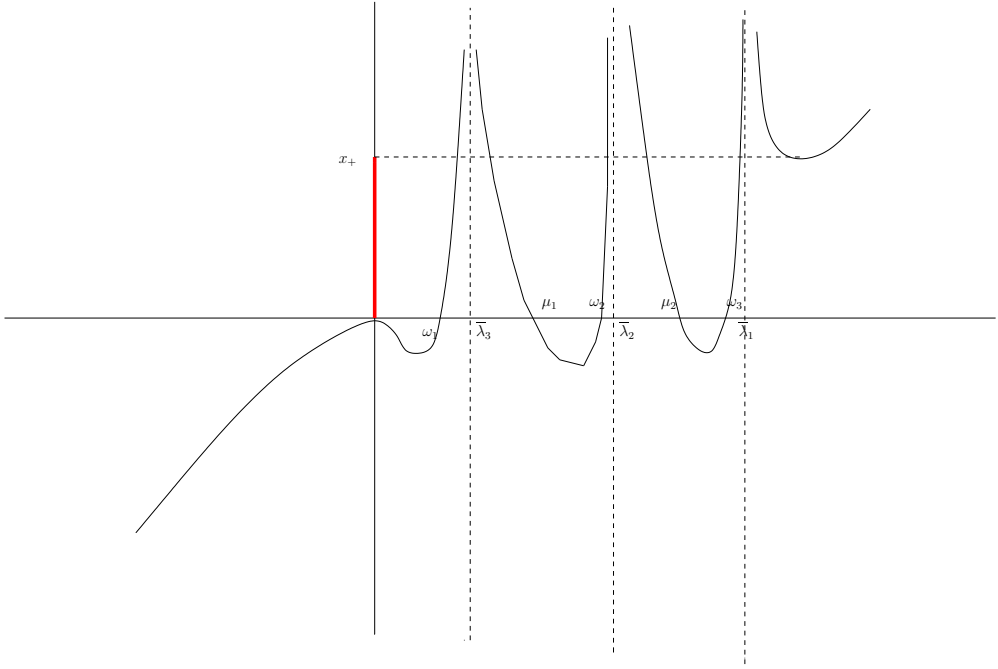}
	\caption
	    {Typical representation of $\phi\left(w\right)$ as a function of $w$ for $\overline{M}=3$. There is
              no local maximum on $[\overline{\lambda}_3, \mu_1]$ and on  $[\overline{\lambda}_2, \mu_2]$,
              so that 
              $\mathcal{S} = [0,x_{+}]$.}
	\label{fig:graph-phi-nolocal}
\end{figure}

We illustrate such a behaviour when $\overline{M} = 3$. In the context of
Fig. \ref{fig:graph-phi-nolocal}, the support is reduced to the single interval $[0,x_{+}]$
because $\phi^{'}(w) \neq 0$ for $w \in  [\overline{\lambda}_3, \mu_1] \cup  [\overline{\lambda}_2, \mu_2]$.

In order to precise the support when $\phi^{'}$ vanishes in  $]\overline{\lambda}_{\overline{M}}, \mu_{1}[ \cup \ldots \cup]\overline{\lambda}_2, \mu_{\overline{M}-1}[$, we need to characterize the corresponding zeros. 
For this, we first justify that $\phi^{'}$ cannot have a multiplicity 2 zero. Assume for example 
that $\phi^{'}$ has a multiplicity 2 zero in $]\overline{\lambda}_{\overline{M}+1-l}, \mu_{l}[$, and 
denote by $w_{l}$ this zero. Then, if $x_l = \phi(w_l)$, the equation $\phi(w) = x_l$ 
has $2 \overline{M} - 1$ simple real roots, and the multiplicity 3 root $w_l$. Therefore, 
the equation  $\phi(w) = x_l$ has $2 \overline{M} + 2$ roots (counting multiplicities), a
contradiction. We now establish the following useful result.  
\begin{proposition}
\label{prop:local-extrema-phi}
The number of local extrema of $\phi_N$ in  $]\overline{\lambda}_{\overline{M}}, \mu_{1}[ \cup \ldots \cup]\overline{\lambda}_2, \mu_{\overline{M}-1}[$ is an even number, say $2q$, with $0 \leq q \leq \overline{M}-1$. 
If $q \geq 1$, we denote the arguments of these extrema by $w_{1,N}^{+} < w_{2,N}^{-} < w_{2,N}^{+} < \ldots < w_{q-1,N}^{+} < w_{q,N}^{-}$, then $x_{1,N}^{+} = \phi_N(w_{1,N}^{+}), x_{2,N}^{-}=\phi_N(w_{2,N}^{-}), \ldots, 
x_{q-1,N}^{+}= \phi_N(w_{q-1,N}^{+}), x_{q,N}^{-} = \phi_N(w_{q,N}^{-})$ verify 
\begin{equation}
\label{eq:entrelacement-x}
x_{1,N}^{+} < x_{2,N}^{-} < x_{2,N}^{+} <\ldots < x_{q-1,N}^{+} < x_{q,N}^{-}
\end{equation}
Moreover, for each $l$, the interval $]\overline{\lambda}_{\overline{M}-(l-1)}, \mu_{l}[$
contains at most one interval $[w_{p,N}^{+}, w_{p+1,N}^{-}]$, and $x_{p,N}^{+}$ (resp. $x_{p+1,N}^{-}$)
is a local minimum (resp. local maximum) of $\phi_N$. 
\end{proposition}
{\bf Proof.} We establish that if $w_1,\,w_2\in\{w^+_1,\,w^-_2,\,\ldots,\,w^{+}_{q-1},\,w^-_{q}\}$ such that $w_1>w_2$, the images $x_1=\phi(w_1)$ and $x_2=\phi(w_2)$ are also satisfy $x_1>x_2$. The goal is to show that ratio $(x_1-x_2)/(w_1-w_2)$ is always positive. For more convenience we put $f_n=\frac{c_N}{M}\tr R_N(R_N-w_nI_M)^{-1}=\frac{c_N}{M}\sum_{1}^{\bar{M}}\frac{\overline{\lambda}_im_i}{\overline{\lambda}_i-w_2}$ for $n=1,2$.  With this and (\ref{eq:def-phi}) we can rewrite
		\begin{align}\label{equation:phi-f}
		x_n=\phi(w_n)=w_n^2f_n(f_n-1)=w_n^2p_n(p_n-1),
		\end{align}
		where $p_n=1-f_n$. Let us notice that extremes $w_1$ and $w_2$ are by definition such that $f_1$ and $f_2$ are negative. Using directly (\ref{equation:phi-f}) for $x_1$ and $x_2$ we can write 
		\begin{multline}\label{expen_x1-x2/w1-w2}
		\dfrac{x_1-x_2}{w_1-w_2}=\dfrac{(w_1^2p_1^2-w_2^2p_2^2)-(w_1^2p_1-w_2^2p_2)}{w_1-w_2}\\
		=(w_1p_1+w_2p_2)\dfrac{w_1p_1-w_2p_2}{w_1-w_2}-\dfrac{w_1^2p_1-w_2^2p_2}{w_1-w_2}
		\end{multline}
		With the definition of $f_{1,2}$ the first term of (\ref{expen_x1-x2/w1-w2}) can be expended as
	\begin{multline*}
\dfrac{w_1p_1-w_2p_2}{w_1-w_2}
=1+\dfrac{c}{M}\sum_{1}^{\bar{M}}\dfrac{\overline{\lambda}_im_1}{w_1-w_2}\left(\dfrac{w_2}{\overline{\lambda}_i-w_2}-\dfrac{w_1}{\overline{\lambda}_i-w_1}\right)\\
=1-\dfrac{c}{M}\sum_{1}^{\bar{M}}\dfrac{\overline{\lambda}_i^2m_i}{(\overline{\lambda}_i-w_1)(\overline{\lambda}_i-w_2)}
\end{multline*}
	And similarly the second one as
	\begin{multline*}
\dfrac{w_1^2p_1-w_2^2p_2}{w_1-w_2}
=(w_1+w_2)+\dfrac{c}{M}\sum_{1}^{\bar{M}}\dfrac{\overline{\lambda}_im_1}{w_1-w_2}\left(\dfrac{w_2^2}{\overline{\lambda}_i-w_2}-
\dfrac{w_1^2}{\overline{\lambda}_i-w_1}\right)\\
=(w_1+w_2)\left(1-\dfrac{c}{M}\sum_{1}^{\bar{M}}\dfrac{\overline{\lambda}_i^2m_i}{(\overline{\lambda}_i-w_1)(\overline{\lambda}_i-w_2)}\right)+w_1w_2\dfrac{c}{M}\sum_{1}^{\bar{M}}\dfrac{\overline{\lambda}_im_i}{(\overline{\lambda}_i-w_1)(\overline{\lambda}_i-w_2)}
\end{multline*}
	Putting the last two equation in (\ref{expen_x1-x2/w1-w2}) we obtain
	\begin{multline*}
	  \dfrac{x_1-x_2}{w_1-w_2}=(w_1p_1+w_2p_2-w_1-w_2)\left(1-\dfrac{c}{M}\sum_{1}^{\bar{M}}\dfrac{\overline{\lambda}_i^2m_i}{(\overline{\lambda}_i-w_1)
            (\overline{\lambda}_i-w_2)}\right)\\
	-w_1w_2\dfrac{c}{M}\sum_{1}^{\bar{M}}\dfrac{\overline{\lambda}_im_i}{(\overline{\lambda}_i-w_1)(\overline{\lambda}_i-w_2)}
	=-(w_1f_1+w_2f_2)\\
	\times\left(1-\dfrac{c}{M}\sum_{1}^{\bar{M}}\dfrac{\overline{\lambda}_i^2m_i}{(\overline{\lambda}_i-w_1)(\overline{\lambda}_i-w_2)}\right)
	-w_1w_2\dfrac{c}{M}\sum_{1}^{\bar{M}}\dfrac{\overline{\lambda}_im_i}{(\overline{\lambda}_i-w_1)(\overline{\lambda}_i-w_2)}
	\end{multline*}
	Now we recall that $-f_n$ is positive as well as $w_1,\,w_2>0$ from what we have $-(w_1f_1+w_2f_2)>0$. That allows us to use the inequality
	\begin{align*}
	\dfrac{1}{(\overline{\lambda}_i-w_1)(\overline{\lambda}_i-w_2)}\le\dfrac{1}{2}\left(\dfrac{1}{(\overline{\lambda}_i-w_1)^2}+\dfrac{1}{(\overline{\lambda}_i-w_2)^2}\right)
	\end{align*}
	and to write
	\begin{multline*}
	\dfrac{x_1-x_2}{w_1-w_2}
	\ge -(w_1f_1+w_2f_2)\left(1-\dfrac{c}{2M}\sum_{1}^{\bar{M}}\dfrac{\overline{\lambda}_i^2m_i}{(\overline{\lambda}_i-w_1)^2}-\dfrac{c}{2M}\sum_{1}^{\bar{M}}\dfrac{\overline{\lambda}_i^2m_i}{(\overline{\lambda}_i-w_2)^2}\right)\\
	-w_1w_2\dfrac{c}{M}\sum_{1}^{\bar{M}}\dfrac{\overline{\lambda}_im_i}{(\overline{\lambda}_i-w_1)(\overline{\lambda}_i-w_2)}
		\end{multline*}
	It is easy to check that $\frac{c}{M}\sum\frac{\overline{\lambda}_i^2m_i}{(\overline{\lambda}_i-w)^2}=f(w)+wf^\prime(w)$. Using this we can rewrite last inequality as
	\begin{multline}
	\label{ineq:x1-x2}
	\dfrac{x_1-x_2}{w_1-w_2}
	\ge -\dfrac{1}{2}(w_1f_1+w_2f_2)\left(2-f_1-w_1f^\prime_1-f_2-w_2f^\prime_2\right)\\
	-w_1w_2\dfrac{c}{M}\sum_{1}^{\bar{M}}\dfrac{\overline{\lambda}_im_i}{(\overline{\lambda}_i-w_1)(\overline{\lambda}_i-w_2)}
	\end{multline}	
	Taking the derivatives of the expression (\ref{equation:phi-f}), we obtain that $\phi^\prime(w_n)=2w_nf_n^2-2w_nf_n+2w_n^2f_nf_n^\prime-w_n^2f_n^\prime$. 	By definition, $w_{1,2}$ are extremes of function $\phi(w)$, i.e. $\phi^\prime(w_{1,2})=0$. This gives immediately $f_n+w_nf^\prime_n-1=\frac{w_nf^\prime_n}{2f_n}$. After putting this into (\ref{ineq:x1-x2}) and regrouping terms we obtain
		\begin{align*}
		&\dfrac{x_1-x_2}{w_1-w_2}
		\ge \dfrac{1}{4}(w_1f_1+w_2f_2)\left(\dfrac{w_1f^\prime_1}{f_1}+\dfrac{w_2f^\prime_2}{f_2}\right)
		-w_1w_2\dfrac{c}{M}\sum_{1}^{\bar{M}}\dfrac{\overline{\lambda}_im_i}{(\overline{\lambda}_i-w_1)(\overline{\lambda}_i-w_2)}\\
		&=\dfrac{1}{4}(w_1^2f^\prime_1+w_2^2f^\prime_2)+\dfrac{1}{4}w_1w_2\left(f^\prime_1\dfrac{f_2}{f_1}+f^\prime_2\dfrac{f_1}{f_2}\right)-w_1w_2\dfrac{c}{M}\sum_{1}^{\bar{M}}\dfrac{\overline{\lambda}_im_i}{(\overline{\lambda}_i-w_1)(\overline{\lambda}_i-w_2)}
		\end{align*}
		Finally, we denote by $I_1,I_2,I_3$ the three parts of r.h.s and show that $I_1+\frac{1}{2}I_3$ and $I_2+\frac{1}{2}I_3$ can be presented as the sum of positive terms.
		Using again the definition of $f_{1,2}$ we expend $I_1+\frac{1}{2}I_3$ as
		 \begin{multline*}
		 \dfrac{1}{4}\left(w_1^2f^\prime_1+w_2^2f^\prime_2-2w_1w_2\dfrac{c}{M}\sum_{1}^{\bar{M}}\dfrac{\overline{\lambda}_im_i}{(\overline{\lambda}_i-w_1)(\overline{\lambda}_i-w_2)}\right)\\
		 =\dfrac{c}{4M}\sum\overline{\lambda}_im_i\Big(\dfrac{w_1^2}{(\overline{\lambda}_i-w_1)^2}+\dfrac{w_2^2}{(\overline{\lambda}_i-w_2)^2}
		 -\dfrac{2w_1w_2}{(\overline{\lambda}_i-w_1)(\overline{\lambda}_i-w_2)}\Big)\\
		 =\dfrac{c}{4M}\sum\overline{\lambda}_im_i\left(\dfrac{w_1}{\overline{\lambda}_i-w_1}-\dfrac{w_2}{\overline{\lambda}_i-w_2}\right)^2
		 \end{multline*}
		 Similarly, $I_2+\frac{1}{2}I_3$ can be written as
		\begin{multline*}
		  \dfrac{1}{4}w_1w_2\left(f^\prime_1\dfrac{f_2}{f_1}+f^\prime_2\dfrac{f_1}{f_2}-2\dfrac{c}{M}\sum_{1}^{\bar{M}}
                  \dfrac{\overline{\lambda}_im_i}{(\overline{\lambda}_i-w_1)(\overline{\lambda}_i-w_2)}\right)\\
		=w_1w_2\dfrac{c}{4M}\sum\overline{\lambda}_im_i
		\left(\dfrac{f_2/f_1}{(\overline{\lambda}_i-w_1)^2}+\dfrac{f_1/f_2}{(\overline{\lambda}_i-w_2)^2}-\dfrac{2}{(\overline{\lambda}_i-w_1)(\overline{\lambda}_i-w_2)}\right)\\
		=w_1w_2\dfrac{c}{4M}\sum\overline{\lambda}_im_i\left(\dfrac{\sqrt{f_2/f_1}}{\overline{\lambda}_i-w_1}-\dfrac{\sqrt{f_1/f_2}}
                {\overline{\lambda}_i-w_2}\right)^2
		\end{multline*}
This shows that $x_1 - x_2 > 0$, and that (\ref{eq:entrelacement-x}) holds. It remains to justify that each interval $(]\overline{\lambda}_{\overline{M}-(l-1)}, \mu_{l}[)_{l=1, \ldots, \overline{M}-1}$ contains at most one interval $[w_{p,N}^{+}, w_{p+1,N}^{-}]$. 
Assume that the interval $]\overline{\lambda}_{\overline{M}-(l-1)}, \mu_{l}[$ contains 2 intervals 
$[w_{p_1,N}^{+}, w_{p_1+1,N}^{-}]$ and $[w_{p_2,N}^{+}, w_{p_2+1,N}^{-}]$ with $p_1 < p_2$. Then, 
it also holds that $[w_{p_1+1,N}^{+}, w_{p_1+2,N}^{-}] \subset ]\overline{\lambda}_{\overline{M}-(l-1)}, \mu_{l}[$.
$x_{p_1,N}^{+}$ is necessarily a local minimum because $x_{p_1,N}^{+} < x_{p_1+1,N}^{-}$ while 
$x_{p_1+1,N}^{-}$ must be a local maximum. The same property holds for $x_{p_1+1,N}^{+}$ and $x_{p_1+2,N}^{-}$. 
However, this contradicts the property $x_{p_1+1,N}^{-} < x_{p_1+1,N}^{+}$. This completes the proof of 
Proposition \ref{prop:local-extrema-phi}.  $\blacksquare$ \\

Proposition \ref{prop:local-extrema-phi} allows to identify the support $\mathcal{S}_N$. 
\begin{corollary}
\label{coro:identification-support}
When $c_N \leq 1$, the support $\mathcal{S}_N$ is given by
\begin{equation}
\label{eq:identification-support}
\mathcal{S}_N = [0,x_{1,N}^{+}] \cup [x_{2,N}^{-}, x_{2,N}^{+}] \cup \ldots [x_{q,N}^{-}, x_{+,N}]
\end{equation}
\end{corollary}
{\bf Proof.} If $x$ belongs to the interior of the righthandside of (\ref{eq:identification-support}), 
$\phi(w) = x$ has only $2 \overline{M} -1$ real solutions. This implies that the 2 remaining 
roots are complex valued, i.e. that $x \in \mathcal{S}^{\circ}$. This leads to the conclusion 
that 
$$
]0,x_{1,N}^{+}[ \cup ]x_{2,N}^{-}, x_{2,N}^{+}[ \cup \ldots ]x_{q,N}^{-}, x_{+,N}[ \subset \mathcal{S}^{\circ}
$$ 
and that 
$$
[0,x_{1,N}^{+}] \cup [x_{2,N}^{-}, x_{2,N}^{+}] \cup \ldots [x_{q,N}^{-}, x_{+,N}] \subset \mathcal{S}
$$
Conversely, if $x \in \mathbb{R}^{+} - \left([0,x_{1,N}^{+}] \cup [x_{2,N}^{-}, x_{2,N}^{+}] \cup \ldots [x_{q,N}^{-}, x_{+,N}]\right)$, the equation $\phi(w)=x$ has $2 \overline{M}+1$ real solutions, which implies that 
$w(x)$ is real. Therefore, 
$$
\mathbb{R}^{+} - \left([0,x_{1,N}^{+}] \cup [x_{2,N}^{-}, x_{2,N}^{+}] \cup \ldots [x_{q,N}^{-}, x_{+,N}]\right)
\subset \mathbb{R}^{+} - \mathcal{S}
$$
or equivalently, 
$$
\mathcal{S} \subset [0,x_{1,N}^{+}] \cup [x_{2,N}^{-}, x_{2,N}^{+}] \cup \ldots [x_{q,N}^{-}, x_{+,N}]
$$
This completes the proof of Corollary (\ref{coro:identification-support}).  $\blacksquare$ \\

\begin{figure}[ht!]
	\centering
	\par
        \includegraphics[scale=0.5]{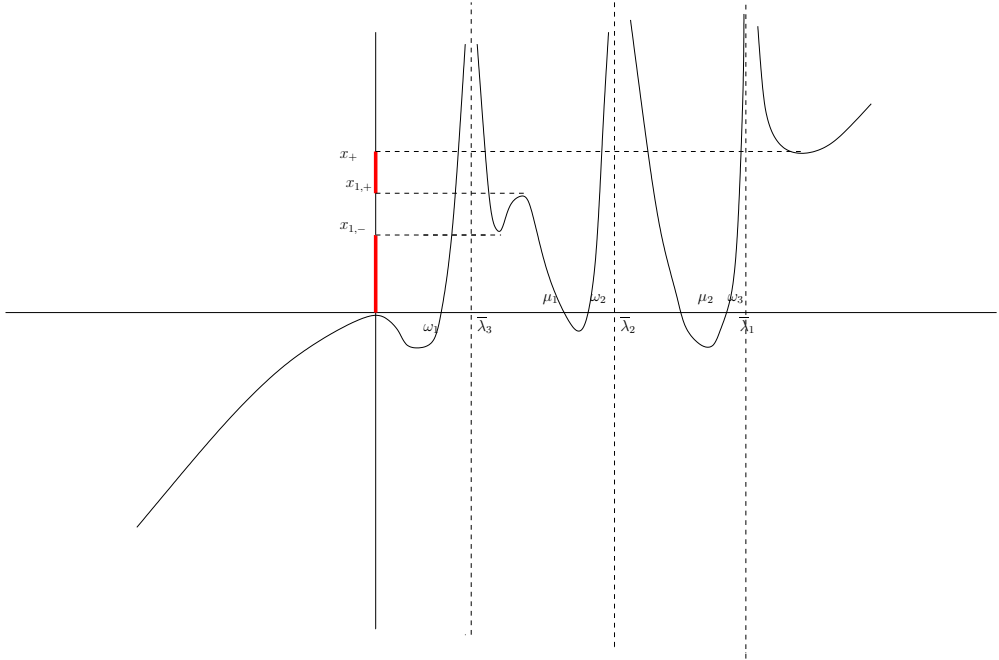}
	\caption
	    {Typical representation of $\phi\left(w\right)$ as a function of $w$ for $\overline{M}=3$. There are
              2 local extrema on  $[\overline{\lambda}_3, \mu_1]$ and no local maximum on $[\overline{\lambda}_2, \mu_2]$, so that  $\mathcal{S} = [0,x_1^{-}]\cup[x_1^{+},x_{+}]$.}
	\label{fig:graph-phi}
\end{figure}

We illustrate the above behaviour when $\overline{M} = 3$. 
In the context of Fig. \ref{fig:graph-phi}, $\phi^{'}$ vanishes on $[\overline{\lambda}_3, \mu_1]$ and not on  $[\overline{\lambda}_2, \mu_2]$. The support thus coincides with
$\mathcal{S} =  [0,x_1^{+}]\cup[x_2^{-},x_{+}]$.

\vspace{0.3cm}
\noindent
When matrix $R$ is reduced to $R = \sigma^{2} I$, i.e. $\overline{M} = 1$ and $\overline{\lambda}_1 = \sigma^{2}$, the support of course coincides with
$\mathcal{S} = [0, x_{+}]$, and $x_{+}$ is given by
\begin{equation}
  \label{eq:expre-x+-sigma2-I}
  x_{+} = \sigma^{4} c \, \left( 1 + \frac{1}{\frac{1+\sqrt{1+8c}}{2}} \right)^{2} \; \left(c + \frac{1+\sqrt{1+8c}}{2} \right)
\end{equation}
Moroever, $w_{+}$ is equal to
\begin{equation}
  \label{eq:expre-w+-sigma2-I}
  w_{+} = \sigma^{2} \left(1 + \frac{1+\sqrt{1+8c}}{2} \right)
\end{equation}
(\ref{eq:expre-x+-sigma2-I}) and (\ref{eq:expre-w+-sigma2-I}) are in accordance with the results of \cite{li-pan-yao-jmva-2015}. \\

We now briefly address the case $c_N > 1$. The behaviour of $\phi_N$ is essentially the same 
as if $c_N \leq 1$, except that the first root $\omega_{1,N}$ of the equation 
$\frac{1}{M} \mathrm{Tr}R_N(R_N - w I)^{-1} = \frac{1}{c_N}$ is now strictly negative. 
As $\phi_N(0) = 0$, this implies that it exists $\omega_{1,N} < w_{N,-} < 0$ for 
which $\phi_N^{'}(w_{N,-})=0$. Moreover, this point is unique, otherwise, the 
equation $\phi_N(w)=x$ would have more than $2 \overline{M}+1$ roots for certain values of $x>0$. 
$x_{-,N} = \phi_N(w_{-,N}) > 0$ is thus a local maximum of $\phi_N$ whose argument is strictly negative. 
We also notice that $\phi_N(w) > 0$ if $0<w<\overline{\lambda}_{\overline{M}}$. Apart these differences, 
the behaviour of $\phi_N$ for $w > \overline{\lambda}_{\overline{M}}$ remains the same as if $c_N \leq 1$. 
In particular, Proposition \ref{prop:local-extrema-phi} still holds true. However, we remark that if 
$0 < x < x_{-,N}$, the equation $\phi_N(w)=x$ has still $2 \overline{M} -1$ real solutions that are
strictly positive, and 2 extra real roots, the smallest one being less than $w_{-,N}$ and the other one 
being negative and largest that $w_{-,N}$. This implies that $w_N(x)$ is real. We also notice that $w_N(x)$ 
coincides with the smallest extra negative root because it satisfies conditions (\ref{eq:proprietes-w(x)}). 
Hence, the interval $]0, x_{-,N}[$ is included into $\mathbb{R}^{+} - \mathcal{S}_N$. If $\phi^{'}_N$ does 
not vanish on $]\overline{\lambda}_{\overline{M}}, \mu_{1}[ \cup \ldots \cup]\overline{\lambda}_2, \mu_{\overline{M}-1}[$, for $x \in ]x_{-,N}, x_{+,N}[$, the equation $\phi_N(w) = x$ has only $2 \overline{M} -1$ real solutions that do not satisfy conditions (\ref{eq:proprietes-w(x)}) and 2 extra complex conjugates solutions. Therefore,   
$]x_{-,N}, x_{+,N}[ \subset \mathcal{S}_N^{\circ}$ and $[x_{-,N}, x_{+,N}] \subset \mathcal{S}_N$. Conversely, 
$]0, x_{-,N}[ \cup ]x_{+,N}, +\infty[ \subset \mathbb{R}^{+} - \mathcal{S}_N$, which implies that 
$\mathcal{S}_N \subset \{ 0 \} \cup [x_{-,N}, x_{+,N}]$. As it was established above that 
$\{ 0 \} \subset \mathcal{S}_N$, we deduce that $\mathcal{S}_N = \{ 0 \} \cup [x_{-,N}, x_{+,N}]$
if $\phi^{'}_N$ does not vanish on $]\overline{\lambda}_{\overline{M}}, \mu_{1}[ \cup \ldots \cup]\overline{\lambda}_2, \mu_{\overline{M}-1}[$. If $\phi^{'}_N$ vanishes on $]\overline{\lambda}_{\overline{M}}, \mu_{1}[ \cup \ldots \cup]\overline{\lambda}_2, \mu_{\overline{M}-1}[$, i.e. if $q \geq 1$ (we recall that $q$ is defined in Proposition 
\ref{prop:local-extrema-phi}), the support is given by 
\begin{equation}
\label{eq:support-c>1}
\mathcal{S}_N = \{ 0 \} \cup [x_{-,N},x_{1,N}^{+}] \cup [x_{2,N}^{-}, x_{2,N}^{+}] \cup \ldots [x_{q,N}^{-}, x_{+,N}]
\end{equation}
To justify this, we just need to establish that $x_{-,N} < x_{1,N}^{+}$, and to use the same arguments 
as in the proof of Corollary \ref{coro:identification-support}. To justify $x_{-,N} < x_{1,N}^{+}$, 
we put $w_1 = w_{-,N}, w_2 = w_{1,N}^{+}$, and follow step by step the arguments used 
to evaluate $\phi(w_2) - \phi(w_1) > 0$. We notice that in contrast with the context of 
the proof of  Corollary \ref{coro:identification-support}, $w_1 < 0$ and $f_1 > 0$. 
However, $f_1 w_1$ is still negative, so that $-(w_1 f_1 + w_2 f_2)$ is still positive. 
This allows to conclude that all the inequalities used in the course of the proof 
of Corollary \ref{coro:identification-support} remain valid, except the evaluation 
of the term $I_2 + I_3/2$ that needs the following simple modification: we express 
$I_2 + I_3/2$ as 
$$
-w_1w_2\dfrac{c}{4M}\sum\lambda_im_i\times
	\left(\dfrac{-f_2/f_1}{(\lambda_i-w_1)^2}+\dfrac{-f_1/f_2}{(\lambda_i-w_2)^2}+\dfrac{2}{(\lambda_i-w_1)(\lambda_i-w_2)}\right)
$$
As $-f_2/f_1$ and $-f_1/f_2$ are positive, it holds that 
$$
I_2 + I_3/2 = -w_1w_2\dfrac{c}{4M}\sum\lambda_im_i\left(\dfrac{\sqrt{-f_2/f_1}}{\lambda_i-w_1}+\dfrac{\sqrt{-f_1/f_2}}{\lambda_i-w_2}\right)^2
$$
Therefore, $I_2 + I_3/2 > 0$, and $\phi(w_2) - \phi(w_1) > 0$ holds. \\

In order to unify the cases $c_N \leq 1$ and $c_N > 1$, we define $x_{-,N}$ for $c_N \leq 1$ by 
$x_{-,N} = 0$, and summarize the above discussion by the following result. 
\begin{theorem}
\label{th:final-characterization-support}
The support $\mathcal{S}_N$ is given by 
\begin{equation}
\label{eq:final-characterization-support}
\mathcal{S}_N =  \{ 0 \} \mathbb{I}_{c_N > 1} \cup [x_{-,N},x_{1,N}^{+}] \cup [x_{2,N}^{-}, x_{2,N}^{+}] \cup \ldots [x_{q,N}^{-}, x_{+,N}]
\end{equation}
\end{theorem}

We now establish that sequences $(w_{+,N})_{N \geq 1}$ and $(x_{+,N})_{N \geq 1}$ are bounded. In other words, 
for each $N$, the support $\mathcal{S}_N$ is included into a compact interval that does not depend on $N$. 
\begin{lemma}
\label{le:w+-x+-bounded}
\begin{equation}
\label{eq:w+-x+-bounded}
\sup_{N \geq 1} w_{+,N} < +\infty, \; \sup_{N \geq 1} x_{+,N} < +\infty
\end{equation}
\end{lemma}
In order to prove this lemma, we use that $w_{+,N} > \lambda_{1,N}$ and that $\phi_N^{'}(w_{+,N}) = 0$. It is easy to check that
\begin{multline*}
\phi_N^{'}(w) = 2 c_N^{2} w \frac{1}{M} \mathrm{Tr}R(w I - R)^{-1} - (c_N w)^{2} \frac{1}{M} \mathrm{Tr}R(w I - R)^{-2} \\ - 2 c_N^{2}  w \left( \frac{1}{M} \mathrm{Tr}R(w I - R)^{-1} \right)^{2} - 2 (c_N w)^{2} \frac{1}{M} \mathrm{Tr}R(w I - R)^{-2} \frac{1}{M} \mathrm{Tr}R(w I - R)^{-1}
\end{multline*}
For $w > b > \lambda_{1,N}$, it is clear that $\| (w I - R)^{-1} \| \leq \frac{1}{w-b}$. Writing that 
$w \frac{1}{M} \mathrm{Tr}R(w I - R)^{-1} = \frac{1}{M} \mathrm{Tr}R + \frac{1}{M} \mathrm{Tr}R^{2}(w I - R)^{-1}$
and $w^{2} \frac{1}{M} \mathrm{Tr}R(w I - R)^{-2} =  \frac{1}{M} \mathrm{Tr}R  +  w \left( \frac{1}{M} \mathrm{Tr}R(w I - R)^{-2} \right) -  \frac{1}{M} \mathrm{Tr}R^{2}(w I - R)^{-1}$, we obtain immediately that $\phi_N^{'}(w)$ can be written as
$$
\phi_N^{'}(w) = c_N^{2} \, \frac{1}{M} \mathrm{Tr} R + \delta_N(w)
$$
where $\delta_N(w)$ verifies $|\delta_N(w)| \leq \delta(w)$ and $w \rightarrow \delta(w)$ is a rational 
function of $w$ that does not depend on $N$ and which converges towards $0$ when $w \rightarrow +\infty$. 
Therefore, for each $\eta  > 0$, it exists $w_1 > b$ such that $\phi_N^{'}(w) >  c_N^{2} \, \frac{1}{M} \mathrm{Tr} R
- \eta$ for each $w \geq w_1$. As $c_N \rightarrow c_*$ and that $\frac{1}{M} \mathrm{Tr} R \geq a$, 
we obtain that  $\phi_N^{'}(w) > \frac{c_*^{2}}{2} \, a$ for $w \geq w_1$. As $\phi_N^{'}(w_{+,N}) = 0$, 
we deduce from this that $w_{+,N} < w_1$. As $w_1$ does not depend on $N$, this establishes that 
$\sup_{N \geq 1} w_{+,N} < +\infty$. To prove that $x_{+,N}$ is bounded, we observe that 
$x_{+,N} = \phi_N(w_{+,N}) < \phi_N(w_1)$. As $w_1 > b$, it is easily seen that
$$
\phi_N(w_1) < 2 c_N^{2} w_1^{2} \left( \frac{b}{(w_1 - b)^{2}} +  \frac{b}{(w_1 - b)} \right)
$$
Therefore, sequences $(\phi_N(w_1))_{N \geq 1}$ and $(x_{+,N})_{N \geq 1}$ are bounded. This completes the proof of Lemma \ref{le:w+-x+-bounded}.
 $\blacksquare$ \\

We finally provide a sufficient condition under which the support is reduced to $\mathcal{S}_N = [0,x_{+,N}]$ if 
$c_N < 1$ and to  $\mathcal{S}_N = \{0 \} \cup [x_{-,N},x_{+,N}]$ if $c_N > 1$. More precisely, the following result holds.
\begin{proposition}
  \label{prop:sufficient-condition-support}
  Assume that it exist $\kappa > 0$ such that for each $M$ large enough, the following condition holds:  
  \begin{equation}
    \label{eq:sufficient-condition-support}
    |\lambda_{k,N} - \lambda_{l,N}| \leq \kappa \left( \frac{|k-l|}{M} \right)^{1/2}
  \end{equation}
  for each pair $(k,l)$, $1 \leq k \leq l \leq M$. Then, for each $M$ large enough, 
 $\mathcal{S}_N = [0,x_{+,N}]$ if 
$c_N \leq 1$ and to  $\mathcal{S}_N = \{0 \} \cup [x_{-,N},x_{+,N}]$ if $c_N > 1$.
\end{proposition}
    {\bf Proof.} We assume that (\ref{eq:sufficient-condition-support}) holds, and that $\mathcal{S}$ does not coincide with $[0,x_{+}]$ or $\mathcal{S} = \{0 \} \cup [x_{-},x_{+}]$ ,
    i.e. $\phi^{'}(w)$ vanishes at a point $w_0$ such that $\lambda_{1} < w_0 < \lambda_M$ and $\frac{1}{M} \mathrm{Tr}R(R-w_0 I)^{-1} < 0$. After some algebra,
    we obtain that $w_0$ satisfies:
    $$
    \frac{1}{M} \mathrm{Tr}\left(R(R-w_0I)^{-1}\right)^{2} = \frac{-\frac{1}{M} \mathrm{Tr}R(R-w_0 I)^{-1}}{1 - 2 c \frac{1}{M} \mathrm{Tr}R(R-w_0 I)^{-1}}
    $$
    As  $\frac{1}{M} \mathrm{Tr}R(R-w_0 I)^{-1} < 0$, this implies that
    \begin{multline*}
    \frac{1}{M} \mathrm{Tr}\left(R(R-w_0I)^{-1}\right)^{2} = \frac{1}{M} \sum_{k=1}^{M} \left(\frac{\lambda_k}{\lambda_k-w_0}\right)^{2}
    < -\frac{1}{M} \mathrm{Tr}R(R-w_0 I)^{-1} \\
    \leq \frac{1}{M} \sum_{k=1}^{M} \frac{\lambda_k}{|\lambda_k-w_0|}
    \end{multline*}
    Jensen's inequality leads to $ \left(\frac{1}{M} \sum_{k=1}^{M} \frac{\lambda_k}{|\lambda_k-w_0|}\right)^{2} \leq \frac{1}{M} \sum_{k=1}^{M} \left(\frac{\lambda_k}{\lambda_k-w_0}\right)^{2}$. Therefore, we obtain that $ \frac{1}{M} \sum_{k=1}^{M} \frac{\lambda_k}{|\lambda_k-w_0|} < 1$, and that
    \begin{equation}
      \label{eq:second-derivative-bounded}
      \frac{1}{M} \sum_{k=1}^{M} \left(\frac{\lambda_k}{\lambda_k-w_0}\right)^{2}  < 1
    \end{equation}
    We assume that $\lambda_{j_0} < w_0 < \lambda_{j_0+1}$. Then, hypothesis (\ref{eq:hypothesis-R-bis}) and condition 
    (\ref{eq:sufficient-condition-support}) imply that
    $$
    \left(\frac{\lambda_k}{\lambda_k-w_0}\right)^{2} > \frac{a^{2}}{\kappa^{2}} \frac{M}{(|k-j_0| + 1)}
    $$
    Hence, it must hold that
    $$
    \frac{a^{2}}{\kappa^{2}}  \, \sum_{k=1}^{M} \frac{1}{(|k-j_0| + 1)} < 1
    $$
    for each $M$ large enough, a contradiction because $\sum_{k=1}^{M} \frac{1}{(|k-j_0| + 1)}$ is easily seen to 
be an unbounded term.  $\blacksquare$

    \section{No eigenvalues outside the support.}
\label{sec:localization-eigenvalues}
    In this paragraph, we establish the following result:
\begin{theorem}\label{th:no_eigen}
	Assume that there exists $\epsilon>0$, $\kappa_1 \in \mathbb{R}$, $\kappa_2 \in \mathbb{R} \cup \{+\infty\}$ and an integer $N_0$ such that
	\begin{align}\label{condition_a,b}
	(\kappa_1-\epsilon,\,\kappa_2+\epsilon)\cap \mathcal{S}_N=\varnothing\qquad\forall \,N\ge N_0.
	\end{align}
	Then with probability one, no eigenvalues of $W_{f,N}W_{p,N}^{*}W_{p,N}W_{f,N}^{*}$ appears in $[\kappa_1,\, \kappa_2]$ for all $N$ large enough.
\end{theorem}
We first remark that it is sufficient to consider the case where $\kappa_2 < +\infty$. To justify this claim, 
we recall that $\cup_{N \geq 1} \mathcal{S}_N$ is a compact subset (see Lemma \ref{le:w+-x+-bounded}), and notice that $\| W_{f,N}W_{p,N}^{*}W_{p,N}W_{f,N}^{*} \| \leq \| W_N \|^{4}$ where matrix $W_N$ is defined 
by (\ref{eq:def-WN}). Moreover, (\ref{eq:def-Wiid}) 
implies that almost surely, for $N$ large enough, $\| W_N \|^{2} \leq b \, (1 + \delta + \sqrt{c_*})^{2}$ 
where $\delta > 0$. Therefore, almost surely, the largest eigenvalue of $W_{f,N}W_{p,N}^{*}W_{p,N}W_{f,N}^{*}$
is for each $N$ large enough upperbounded by the nice constant $ b^{2} \, (1 + \delta + \sqrt{c_*})^{4}$. 
This justifies that it is sufficient to assume that $\kappa_2 < +\infty$ in the following. \\

In order to establish Theorem \ref{th:no_eigen}, we use the Haagerup-Thornbjornsen approach (\cite{HT:05}, see also \cite{capitaine-donati-martin-feral-2007}). 
The crucial step of the proof is the following Proposition.

\begin{proposition}\label{prop:expQ-T}
	$\forall z\in\mathbb{C^+},$ we have for $N$ large enough,
	\begin{align}\label{eq:EQ=T+r}
	\mathbb{E}\left\{\dfrac{1}{ML}\tr Q_N(z)\right\}=\dfrac{1}{M}\tr T_N(z)+\dfrac{1}{N^2}r_N(z)
	\end{align}
	where $r_N$ is holomorphic in $\mathbb{C^+}$ and satisfies
	\begin{align}\label{ineq:condition_r}
	|r_N(z)|\le P_1(|z|)P_2\left(\dfrac{1}{\im  z}\right)
	\end{align}
	for each $z\in\mathbb{C^+}$, where $P_1$ and $P_2$ are nice polynomials.
\end{proposition} 

{\bf Proof.} To prove (\ref{eq:EQ=T+r}) we write 
\begin{multline*}
\mathbb{E}\left\{\dfrac{1}{ML}\tr Q_N(z)\right\}-\dfrac{1}{M}\tr T_N(z)=\dfrac{1}{ML}\tr\left[\mathbb{E}\left\{ Q_N(z)\right\}-I_L\otimes S_N(z)\right]\\
+\dfrac{1}{M}\tr \left[S_N(z)-T_N(z)\right]
\end{multline*}
As (\ref{eq:convergence-E(Q)-S-F=I}) holds, it is sufficient to establish that
\begin{equation}
  \label{eq:evaluation-rate-S-T}
  \left| \frac{1}{M}\tr[S_N(z)-T_N(z)] \right| \leq  \frac{1}{N^2}P_1(|z|)P_2(\im ^{-1} z)
  \end{equation}
for some nice polynomial $P_1$ and $P_2$. In the following, we denote by $s_N(z)$ the function defined by
\begin{equation}
  \label{eq:def-s}
  s_N(z) = \frac{1}{M} \mathrm{Tr}R_N S_N(z)
\end{equation}
It is clear that $s_N \in \mathcal{S}(\mathbb{R}^{+})$. Moreover,  if $\mu_{N,s}$ represents 
the associated positive measure, then we have
\begin{equation}
\label{eq:masse-first-moment-mu-s}
\mu_{N,s}(\mathbb{R}^{+}) = \frac{1}{M} \tr R_N, \; \int_{\mathbb{R}^{+}} \lambda \, d \mu_{N,s}(\lambda) =
c_N \, \frac{1}{M} \tr R_N \, \frac{1}{M} \tr R_N^{2}
\end{equation}
(\ref{eq:masse-first-moment-mu-s}) can be proved using the arguments of the proof of Proposition
\ref{prop:existence-uniqueness}.

As $ \frac{1}{M}\tr[S_N(z)-T_N(z)]$ is given by (\ref{eq:(S-T)F}) for $F=I$, (\ref{eq:evaluation-rate-S-T}) appears equivalent to the property
\begin{equation}
  \label{eq:evaluation-rate-R(S-T)}
\left| \frac{1}{M}\tr[R_N(S_N(z)-T_N(z))] \right|  = |s_N(z) - t_N(z)| \leq   \frac{1}{N^2}P_1(|z|)P_2(\im ^{-1} z)
  \end{equation}
In order to prove (\ref{eq:evaluation-rate-R(S-T)}), we define the following functions that appear formally similar to functions
$u(z)$ and $v(z)$ defined by (\ref{eq:def-u}) and (\ref{eq:def-v}):
\begin{align}
&u_{\alpha}(z)=c\dfrac{|cz\alpha(z)|^2\frac{1}{M}\tr(RS(z)S^*(z)R)}{|1-z(c\alpha(z))^2|^2}\\
&v_{\alpha}(z)=c\dfrac{\frac{1}{M}\tr(RS(z)S^*(z)R)}{|1-z(c\alpha(z))^2|^2}\\
&u_{t,\alpha}(z)=c\dfrac{|cz|^2t(z)\alpha(z)\frac{1}{M}\tr(RS(z)T(z)R)}{(1-z(c\alpha(z))^2)(1-z(ct(z))^2)}\label{eq:def_u_t_alpha}\\
&v_{t,\alpha}(z)=c\dfrac{\frac{1}{M}\tr(RS(z)T(z)R)}{(1-z(c\alpha(z))^2)(1-z(ct(z))^2)}\label{eq:def_v_t_alpha}
\end{align}
Using equation $t(z)=\frac{1}{M}\tr RT(z)$ and the definition of $s(z)$ and $S(z)$, we obtain easily that 
\begin{align}
\begin{pmatrix}
(s(z)-t(z))\\
z(s(z)-t(z))
\end{pmatrix}=\mathbf{D}_{t,\alpha}(z)\begin{pmatrix}
(s(z)-t(z))\\
z(s(z)-t(z))
\end{pmatrix}+\begin{pmatrix}
\epsilon_1(z)\\
\epsilon_2(z)
\end{pmatrix}
\end{align}
holds, where
\begin{align}
&\epsilon_1(z)=(\alpha(z)-s(z))(zv_{t,\alpha}(z)+u_{t,\alpha}(z))\\
&\epsilon_2(z)=z(\alpha(z)-s(z))(zv_{t,\alpha}(z)+u_{t,\alpha}(z))\\
&\mathbf{D}_{t,\alpha}(z)=\begin{pmatrix}
u_{t,\alpha}(z)& v_{t,\alpha}(z)\\
z^2v_{t,\alpha}(z)&u_{t,\alpha}(z)
\end{pmatrix}
\end{align} 
This can also be written as
\begin{align}
\label{eq:system-s-t}
(\mathbf{I}-\mathbf{D}_{t,\alpha}(z))\begin{pmatrix}
(s(z)-t(z))\\
z(s(z)-t(z))
\end{pmatrix}=\begin{pmatrix}
\epsilon_1(z)\\
\epsilon_2(z)
\end{pmatrix}
\end{align}
(\ref{eq:convergence-E(Q)-S}) leads to $\alpha(z)-s(z)=\mathcal{O}_z(N^{-2})$. In order to verify that $(\epsilon_i(z))_{i=1,2}$ are
$\mathcal{O}_z(N^{-2})$ as well, we have to control $u_{t,\alpha}$ and $v_{t,\alpha}$. As $t(z), \alpha(z), \|T(z)\|$ and $\|S(z)\|$ are
$\mathcal{O}_z(1)$ terms, it is sufficient to evaluate the denominator of the right handside of (\ref{eq:def_u_t_alpha}). As the mass and the first moment of $\mu$ and $\overline{\mu}$ (the measure 
associated to $\alpha(z)$) both verify the conditions of Lemma \ref{le:gestion-(1-z(cbeta)2}, 
this Lemma implies that $(1-z(ct(z))^2)^{-1} = \mathcal{O}_z(1)$ and  $(1-z(c\alpha(z))^2)^{-1} = \mathcal{O}_z(1)$. Therefore, we have checked that $(\epsilon_i(z))_{i=1,2}$ are $\mathcal{O}_z(N^{-2})$ terms.  $\blacksquare$ \\

In order to evaluate $s(z) - t(z)$, it is of course necessary to show that matrix $I-\mathbf{D}_{t,\alpha}(z)$ is invertible on
$\mathbb{C}^{+}$, and to control the action of its inverse on the vector $(\epsilon_1(z), \epsilon_1(z))^{T}$. 
We define matrix $\mathbf{D}_{\alpha}$ by 
\begin{align}
\mathbf{D}_{\alpha}(z)=\begin{pmatrix}
u_{\alpha}(z)& v_{\alpha}(z)\\
z^2v_{\alpha}(z)&u_{\alpha}(z)
\end{pmatrix}
\end{align}
and establish the following  result. 
\begin{lemma}\label{le:eval-det-D-alpha}
  For each $z \in \mathbb{C}^{+}$, it exist nice constants $\kappa$ and $\beta$ such that
  \begin{equation}
    \label{eq:lower-bound-det(I-D)}
    \det(I - \mathbf{D}(z)) \geq \dfrac{\kappa \, \left(\im  z\right)^{8}}{(|\beta|^2+|z|^2)^4}
  \end{equation}
  Moreover, it exist 2 nice polynomials $P_1$ and $P_2$ for which
  \begin{equation}
    \label{eq:1-ualpha>0}
    1 - u_{\alpha}(z) > 0
  \end{equation}
  and 
 \begin{equation}
    \label{eq:lower-bound-det(I-Dalpha)}
    \det(I - \mathbf{D}_{\alpha}(z)) \geq \dfrac{\kappa \, \left(\im  z \right)^{8}}{(|\beta|^2+|z|^2)^4}
  \end{equation}
for each $z \in B_N$, where $B_N$ is defined as 
\begin{align}
\label{eq:def-BN}
\mathcal{B}_N=\left\{z\in\mathbb{C^+},\;\dfrac{1}{MN}P_1(|z|)P_2
\left(\dfrac{1}{\im  z}\right)\le 1 \right\}
\end{align}
Finally, for each $z \in \mathcal{B}_N$, it holds that
\begin{equation}
    \label{eq:lower-bound-det(I-Dtalpha)}
    \det(I - \mathbf{D}_{t,\alpha}(z)) \geq \dfrac{\kappa \, \left(\im  z \right)^{8}}{(|\beta|^2+|z|^2)^4}
  \end{equation}
\end{lemma}
{\bf Proof.} 
To evaluate $\det(I-\mathbf{D}(z))$, we use the calculations of the proof of Lemma~\ref{le:properties-u-v}. In particular, we have
\begin{align}\label{eq:final-system-imt-imzt}
(I-\mathbf{D}(z))\begin{pmatrix}
\im  t(z)\\
\im  zt(z) 
\end{pmatrix}=\im  z\begin{pmatrix}
\frac{1}{M}\tr RT(z)T^*(z)\\
0
\end{pmatrix}
\end{align}
This implies that 
\begin{align}
1-u(z)=\dfrac{\im  z}{\im  t(z)}\cdot\dfrac{1}{M}\tr RT(z)T^*(z) +\dfrac{\im  zt(z)}{\im  t(z)}v(z)\ge \dfrac{\im  z}{\im  t(z)}\cdot\dfrac{1}{M}\tr RT(z)T^*(z) 
\end{align}
By applying Cramer's rule to (\ref{eq:final-system-imt-imzt}), we obtain that 
\begin{align}\label{ineq:evaluation_det_D_t}
\det(I-\mathbf{D}(z))=\dfrac{\im  z}{\im  t(z)}\cdot\dfrac{1}{M}\tr RT(z)T^*(z)(1-u(z))\ge \left(\dfrac{\im  z}{\im  t(z)}\cdot\dfrac{1}{M}\tr RT(z)T^*(z)\right)^2
\end{align}
It is clear that $\im  t(z) \leq |t(z)| \leq \frac{1}{M} \tr R \, \left( \im  z \right)^{-1} \leq b \left( \im  z \right)^{-1}$. Therefore, it holds that
$\dfrac{\im  z}{\im  t(z)} \geq \frac{1}{b} \left( \im  t(z) \right)^{2}$. We now evaluate $\frac{1}{M}\tr RT(z)T^*(z)$. For this, we remark that
\begin{equation}
  \label{eq:control-tr-RTT}
  \frac{1}{M}\tr RT(z)T^*(z) = \frac{1}{M}\tr RT(z)T^*(z)RR^{-1} \geq \frac{1}{b} \frac{1}{M} \tr( RT(z)T^*(z)R)
\end{equation}
Jensen's inequality implies that $\frac{1}{M} \tr( RT(z)T^*(z)R) \geq \left| \frac{1}{M} \tr RT(z) \right|^{2} = |t(z)|^{2} \geq \left( \im  \, t(z) \right)^{2}$. Therefore, the application of Lemma \ref{le:gestion-(1-z(cbeta)2} to $\beta(z) = t(z)$ implies that 
$$
\left(\dfrac{\im  z}{\im  t(z)}\cdot\dfrac{1}{M}\tr RT(z)T^*(z)\right)^2 \geq \dfrac{\kappa \, \left(\im  z\right)^{8}}{(|\beta|^2+|z|^2)^4}
$$
for some nice constants $\kappa$ and $\beta$. (\ref{eq:lower-bound-det(I-D)}) thus follows from (\ref{ineq:evaluation_det_D_t}). \\

We now establish (\ref{eq:1-ualpha>0}) and (\ref{eq:lower-bound-det(I-Dalpha)}), and denote by $\epsilon(z)$ the function $\epsilon(z) = \alpha(z) - s(z)$. Using
the equation $s(z) = \frac{1}{M} \tr R S(z)$, and calculating $\im  \, s(z)$ and $\im  \, z s(z)$, we obtain immediately that 
\begin{align}
  \label{eq:equation-ims-imzs}
(\mathbf{I}-\mathbf{D}_{\alpha}(z))\begin{pmatrix}
\im \alpha (z)\\ 
\im  z\alpha (z)
\end{pmatrix}=\im  z\begin{pmatrix}
\dfrac{1}{M}\tr RS(z)S^*(z)\\
0
\end{pmatrix}+\begin{pmatrix}
\im  \epsilon(z)\\
\im  z \epsilon(z)
\end{pmatrix}.
\end{align}
The first component of (\ref{eq:equation-ims-imzs}) leads to 
\begin{align}\label{ineq:1-u_alpha}
1-u_{\alpha}
=\dfrac{\im  z}{\im \alpha}\cdot\dfrac{1}{M}\tr RSS^*+\dfrac{\im  \epsilon}{\im \alpha}+\dfrac{\im  z\alpha}{\im \alpha}v_{\alpha}
\ge \dfrac{\im  z}{\im \alpha}\cdot\dfrac{1}{M}\tr RSS^*+\dfrac{\im  \epsilon}{\im \alpha}
\end{align}
Using the same arguments as above, we obtain that $\dfrac{1}{M}\tr RSS^* \geq \frac{1}{b} \, |s(z)|^{2} 
\geq \frac{1}{b} \, \left(\im  s(z) \right)^{2}$. As (\ref{eq:masse-first-moment-mu-s}) holds, 
we can apply Lemma \ref{le:gestion-(1-z(cbeta)2} to $\beta(z) = s(z)$ and obtain as above that
$$
\dfrac{\im  z}{\im  s(z)}\cdot\dfrac{1}{M}\tr RS(z)S^*(z) \geq \dfrac{\kappa \, \left(\im  z\right)^{4}}{(|\beta|^2+|z|^2)^2}
$$
for some nice constants $\beta$ and $\kappa$. We remark that $\dfrac{\im  \epsilon}{\im \alpha} \geq - \frac{|\epsilon|}{\im  \alpha}$. Therefore, by Lemma \ref{le:gestion-(1-z(cbeta)2} applied to $\beta(z) = \alpha(z)$, it holds 
that $\dfrac{\im  \epsilon}{\im \alpha} \geq - \kappa_1 |\epsilon| \frac{\beta_1^{2} + |z|^{2}}{\im  z}$
for some nice constants $\kappa_1$ and $\beta_1$.  
As $|\epsilon(z)| \leq \frac{1}{N^{2}} Q_1(|z|) Q_2(\frac{1}{\im  z})$ for some nice polynomials 
$Q_1$ and $Q_2$,we obtain that 
\begin{equation}
\label{eq:lower-bound-1-ualpha}
1 - u_{\alpha} \geq  \dfrac{\im  z}{\im \alpha}\cdot\dfrac{1}{M}\tr RSS^*+\dfrac{\im  \epsilon}{\im \alpha} \geq \dfrac{\im  z}{\im \alpha}\cdot\dfrac{1}{M}\tr RSS^* - \dfrac{|\epsilon|}{\im \alpha} \geq 
\frac{1}{2} \, \dfrac{\kappa \, \left(\im  z\right)^{4}}{(|\beta|^2+|z|^2)^2}
\end{equation}
if $z$ belongs to the set $\mathcal{B}_{1,N}$ defined by 
$$
 \dfrac{\kappa \, \left(\im  z\right)^{4}}{(|\beta|^2+|z|^2)^2} - \frac{1}{N^{2}} Q_1(|z|) Q_2(\frac{1}{\im  z}) \,  \kappa_1 \, \frac{\beta_1^{2} + |z|^{2}}{\im  z} \geq \frac{1}{2} \, \dfrac{\kappa \, \left(\im  z\right)^{4}}{(|\beta|^2+|z|^2)^2}
$$
The set $\mathcal{B}_{1,N}$ is clearly defined in the same way than $\mathcal{B}_N$, 
but from 2 other nice polynomials $P_{1,1}$ and $P_{2,1}$. 

Using the Cramer rule, we obtain that $\det(\mathbf{I}-\mathbf{D}_{\alpha})$ can be written as
\begin{align}
\det(\mathbf{I}-\mathbf{D}_{\alpha})
=\left(\dfrac{\im  z}{\im\alpha}\cdot\dfrac{1}{M}\tr RSS^*+\dfrac{\im \epsilon}{\im\alpha}\right)(1-u_{\alpha})+\dfrac{\im z \epsilon}{\im\alpha}v_{\alpha}
\end{align}
Plugging (\ref{eq:lower-bound-1-ualpha}) in the last equation, we get that 
the inequality 
\begin{align}
\det(\mathbf{I}-\mathbf{D}_{\alpha})\ge \left(\frac{1}{2} \, \dfrac{\kappa \, \left(\im z\right)^{4}}{(|\beta|^2+|z|^2)^2} \right)^{2} - \frac{|z| \, |\epsilon|}{\im \alpha} \, v_{\alpha} 
\end{align}
holds for each $z \in \mathcal{B}_{1,N}$. As $v_{\alpha} = \mathcal{O}_z(1)$, we obtain 
that 
$$
 \left(\dfrac{\kappa \, \left(\im z\right)^{4}}{(|\beta|^2+|z|^2)^2} \right)^{2} - \frac{|z| \, |\epsilon|}{\im \alpha} \, v_{\alpha} \geq \left(\frac{1}{4} \, \dfrac{\kappa \, \left(\im z\right)^{4}}{(|\beta|^2+|z|^2)^2} \right)^{2} 
$$
for each $z \in \mathcal{B}_{2,N}$, where $\mathcal{B}_{2,N}$ is defined as $\mathcal{B}_N$ from 
2 nice polynomials $P_{1,2}$ and $P_{2,2}$. We put $P_1(|z|)=P_{1,1}(|z|) + P_{1,2}(|z|)$ and 
$P_2(1/\im z) = P_{2,1}(1/\im z) + P_{2,2}(1/\im z)$, and consider the set $\mathcal{B}_N$ defined by  
(\ref{eq:def-BN}). It is clear that $\mathcal{B}_N \subset \mathcal{B}_{1,N} \cap \mathcal{B}_{2,N}$,
and that (\ref{eq:1-ualpha>0}) and (\ref{eq:lower-bound-det(I-Dalpha)}) hold if $z \in \mathcal{B}_N$. \\

It remains to establish (\ref{eq:lower-bound-det(I-Dtalpha)}). For this, we remark that the inequalities
\begin{multline*}
|\det (\mathbf{I}-\mathbf{D}_{t,\alpha}(z))|\ge |1-u_{t,\alpha}(z)|^2-|z|^2|v_{t,\alpha}(z)|^2\ge (1-|u_{t,\alpha}(z)|)^2\\
-|z|v_{\alpha}(z)\cdot|z|v_t(z)
\ge (1-\sqrt{u(z)u_{\alpha}(z)})^2-|z|v_{\alpha}(z)\cdot|z|v(z)
\ge (1-u(z))(1-u_{\alpha}(z))\\
-|z|v_{\alpha}(z)\cdot|z|v(z)
\ge \sqrt{((1-u(z))^2-|z|^2v(z))((1-u_{\alpha}(z))^2-|z|^2v_{\alpha}(z))}\\
=\sqrt{\det(\mathbf{I}-\mathbf{D}(z))\det(\mathbf{I}-\mathbf{D}_{\alpha}(z))}
\end{multline*} 
hold for each $z \in \mathcal{B}_N$. Therefore, (\ref{eq:lower-bound-det(I-Dtalpha)}) follows 
from (\ref{eq:lower-bound-det(I-D)}) and (\ref{eq:lower-bound-det(I-Dalpha)}). This completes
the proof of Lemma \ref{le:eval-det-D-alpha}.  $\blacksquare$ \\

Solving (\ref{eq:system-s-t}), we obtain immediately that it exists 2 nice polynomials 
$Q_1$ and $Q_2$ such that,
$$
|s_N(z) - t_N(z)| \leq \frac{1}{MN} Q_1(|z|) Q_2(\frac{1}{\im z}) 
$$
holds for each $z \in \mathcal{B}_N$. If $z \in \mathcal{B}_N^{c}$, we use the argument in \cite{HT:05}. More precisely, 
if $z \in \mathcal{B}_N^{c}$, the inequality $1 < \frac{1}{MN} P_1(|z|) P_2(1/\im z)$ holds. 
As $|s_N(z) - t_N(z)| \leq 2 \, \frac{1}{M} \tr R_N \, \frac{1}{\im z}$ on $\mathbb{C}^{+}$, 
we deduce that 
$$
|s_N(z) - t_N(z)| \leq 2 b  \frac{1}{MN} P_1(|z|) \frac{P_2(1/\im z)}{\im z}
$$
for each $z \in  \mathcal{B}_N^{c}$. This, in turn, leads to the conclusion that 
$s_N(z) - t_N(z) = \mathcal{O}_z(\frac{1}{N^{2}})$ for each $z \in \mathbb{C}^{+}$. 
This establishes (\ref{eq:evaluation-rate-R(S-T)}) and $\frac{1}{M} \mathrm{Tr}(T_N(z) -S_N(z)) = \mathcal{O}_z(\frac{1}{N^{2}})$ as expected. This completes the proof of Proposition~\ref{prop:expQ-T}.  $\blacksquare$

We now follow \cite{capitaine-donati-martin-feral-2009} and \cite{HT:05} and use the following Lemma

\begin{lemma}\label{lemma:exp_of_tr_phi}
	Let $\phi$ be a compactly supported real valued smooth function defined on $\mathbb{R}^+$, i.\,e. $\phi\in\mathcal{C}^{\infty}_c(\mathbb{R}^+,\mathbb{R}^+)$. Then,
	\begin{align}
	\mathbb{E}\left\{\dfrac{1}{ML}\tr\phi(W_fW_p^*W_pW_f^*)\right\}-\int_{\mathcal{S}_N}\phi(\lambda)d\mu_N(\lambda)=\mathcal{O}\left(\dfrac{1}{N^2}\right)
	\end{align}
\end{lemma}
\textbf{Proof.} Due to Proposition~\ref{stil_tran} we can write
\begin{align}
\ex\left\{\dfrac{1}{ML}\tr\phi(W_fW^*_pW_pW^*_f)\right\}
=\dfrac{1}{\pi}\lim\limits_{y\downarrow0}\im\left\{\int_{\mathbb{R^+}}\phi(x)\ex\left\{\dfrac{1}{ML}\tr Q(x+iy)\right\}dx\right\}
\end{align}
as well as 
\begin{align}
\int_{\mathcal{S}_N}\phi(\lambda)d\mu_N(\lambda)
=\dfrac{1}{\pi}\lim\limits_{y\downarrow0}\im\left\{\int_{\mathbb{R^+}}\phi(x)\ex\left\{\dfrac{1}{ML}\tr T(x+iy)\right\}dx\right\}
\end{align}
Using Proposition~\ref{prop:expQ-T}, we obtain
\begin{multline}\label{eq:exp_of_tr_phi}
\ex\left\{\dfrac{1}{ML}\tr\phi(W_fW^*_pW_pW^*_f)\right\}
-\int_{\mathcal{S}_N}\phi(\lambda)d\mu_N(\lambda)\\
=\dfrac{1}{N^2}\dfrac{1}{\pi}\lim\limits_{y\downarrow 0}\im\left\{\int_{\mathbb{R^+}}\phi(x)r_N(x+iy)dx\right\}
\end{multline}
Since the function $r_N(z)=\mathcal{O}_z(1)$, we can use the result which was proved in \cite[Section~3.3]{capitaine-donati-martin-feral-2007} and obtain
\begin{align}
\limsup_{y\downarrow0}\left|\int_{\mathbb{R}+}\phi(x)r_N(x+iy)dx\right|\le \kappa,
\end{align}
for some nice constant $\kappa$. This and (\ref{eq:exp_of_tr_phi}) complete the proof.  $\blacksquare$ \\

In order to establish Theorem~\ref{th:no_eigen}, we introduce a function $\phi\in\mathcal{C}_c^{\infty}$ such that $0\le\phi(\lambda)\le 1$ and
\begin{align*}
\phi(\lambda)=\begin{cases}
1,\, \text{for } \lambda\in[\kappa_1, \kappa_2],\\
0,\, \text{for } \lambda\in\mathbb{R}-(\kappa_1-\epsilon,\kappa_2+\epsilon)
\end{cases}
\end{align*}
Since for $N$ large enough $(\kappa_1-\epsilon,\kappa_2+\epsilon)\cap \mathcal{S}_N=\varnothing$ then $\int_{\mathcal{S}_N}\phi(\lambda)d\mu_N(\lambda)=0$ and according to Lemma~\ref{lemma:exp_of_tr_phi}
\begin{align*}
\mathbb{E}\left\{\dfrac{1}{ML}\tr\phi(W_fW_p^*W_pW_f^*)\right\}=\mathcal{O}\left(\dfrac{1}{N^2}\right).
\end{align*}
Now we show that
\begin{align*}
\var\left\{\dfrac{1}{ML}\tr\phi(W_fW_p^*W_pW_f^*)\right\}=\mathcal{O}\left(\dfrac{1}{N^4}\right)
\end{align*}
For this we use again the Poincare-Nash inequality 
\begin{multline*}
\var\{\tr\phi(W_fW_p^*W_pW_f^*)\}
\le \sum \ex\Big\{\left(\dfrac{\partial\tr\phi(W_fW_p^*W_pW_f^*)}{\partial \overline{W}^{m_1}_{i_1,j_1}}\right)^*\ex\{W^{m_1}_{i_1,j_1}\overline{W}^{m_2}_{i_2,j_2}\}\\
\times\dfrac{\partial\tr\phi(W_fW_p^*W_pW_f^*)}{\partial \overline{W}^{m_2}_{i_2,j_2}}\Big\}
+\sum \ex\left\{\dfrac{\partial\tr\phi(WW^*)}{\partial W^{m_1}_{i_1,j_1}}\ex\{W^{m_1}_{i_1,j_1}\overline{W}^{m_2}_{i_2,j_2}\}\left(\dfrac{\partial\tr\phi(WW^*)}{\partial W^{m_2}_{i_2,j_2}}\right)^*\right\}
\end{multline*}
We only evaluate the first term of the r.h.s. of the inequality, denoted by $\psi$, because the second is similar. For this we write first
\begin{align*}
\dfrac{\partial\tr \phi(W_fW_p^*W_pW_f^*)}{\partial \overline{W}^{m_1}_{i_1,j_1}}
=\tr \left(\phi^\prime(W_fW_p^*W_pW_f^*)\dfrac{\partial W_fW_p^*W_pW_f^*}{\partial \overline{W}^{m_1}_{i_1,j_1}}\right)\\
= \begin{cases}
1\le i_1\le L,\, (W_pW_f^*\phi^\prime(W_fW_p^*W_pW_f^*)W_f)^{m_1}_{i_1j_1},\\
L+1\le i_1\le 2L,\, (\phi^\prime(W_fW_p^*W_pW_f^*)W_f^*W_fW_p)^{m_1}_{(i_1-L)j_1}
\end{cases}
\end{align*}
Plugging this into (\ref{var_phi}) we obtain
\begin{multline*}
\psi=\sum_{i_1i_2=1}^{L} \sum_{j_1,j_2,m_1,m_2}\Big( \dfrac{1}{N}\ex\Big\{\left(W_pW_f^*\phi^\prime(W_fW_p^*W_pW_f^*)W_f\right)^{*m_1}_{i_1j_1}
R_{m_1m_2}\delta_{i_1+j_1,i_2+j_2}\\
\times
\left(W_pW_f^*\phi^\prime(W_fW_p^*W_pW_f^*)W_f\right)^{m_2}_{i_2,j_2}\Big\}
+\dfrac{1}{N}\ex\Big\{\left(\phi^\prime(W_fW_p^*W_pW_f^*)W_fW_p^*W_p\right)^{*m_1}_{i_1j_1}\\
\times R_{m_1m_2}\delta_{i_1+j_1,i_2+j_2}
\left(\phi^\prime(W_fW_p^*W_pW_f^*)W_fW_p^*W_p\right)^{m_2}_{i_2,j_2}\Big\}\Big).
\end{multline*}
Following the proof of Lemma~\ref{moment}, we obtain
\begin{multline}
\var\{\tr\phi(W_fW_p^*W_pW_f^*)\}
\le\dfrac{C}{N}\ex\{\tr W_f^*\phi^\prime(W_fW_p^*W_pW_f^*)W_fW_p^*W_pW_f^*\label{ineq:var_tr_phi}\\
\times \phi^\prime(W_fW_p^*W_pW_f^*)W_f\}
+\dfrac{C}{N}\ex\{\tr W_fW_p^*W_pW_p^*W_pW_f^* \left(\phi^\prime(W_fW_p^*W_pW_f^*)\right)^2\}.
\end{multline}
To evaluate the first term ($\psi_1$) of the r.h.s of (\ref{ineq:var_tr_phi}) we denote $\eta(\lambda)=(\phi^\prime(\lambda))^2\lambda$ and write 
\begin{multline*}
\dfrac{L}{N}\ex\left\{\tr W_f^*\phi^\prime(W_fW_p^*W_pW_f^*)W_fW_p^*W_pW_f^*\phi^\prime(W_fW_p^*W_pW_f^*)W_f\right\}\\
\le\dfrac{L}{N}\ex\left\{\|W_f\|^2\tr (\eta(W_fW_p^*W_pW_f^*))\right\}.
\end{multline*}
We recall that (\ref{eq:def-Wiid}) implies that $\|W_f\|^2 \leq b \|W_{iid}\|^2$. Therefore, it holds that 
\begin{multline*}
\psi_1 \leq \dfrac{\kappa}{N}\ex\{\|W_{iid}\|^2\mathbf{1}_{\|W_{iid}\|\le(1+\sqrt{c_*})^2+\delta}\tr(\eta(W_fW_p^*W_pW_f^*))\}\\
+  \dfrac{\kappa}{N}\ex\{\|W_{iid}\|^2\mathbf{1}_{\|W_{iid}\|>(1+\sqrt{c_*})^2+\delta}\tr(\eta(W_fW_p^*W_pW_f^*))\}\\ \le\dfrac{\kappa}{N}\ex\{\tr(\eta(W_fW_p^*W_pW_f^*))\}
+\kappa \ex^{1/2}\{\|W_{iid}\|^4\mathbf{1}_{\|W_{iid}\|>(1+\sqrt{c_*})^2+\delta}\}\\
\times\ex^{1/2}\left\{\left(\dfrac{1}{N}\tr(\eta(W_fW_p^*W_pW_f^*))\right)^2\right\}
\end{multline*}
Lemma~\ref{lemma:exp_of_tr_phi} implies that $\frac{1}{N}\ex\{\tr(\eta(W_fW^*_pW_pW^*_f))\}=\mathcal{O}(N^{-2})$. Throughout the proof of Lemma~\ref{moment}, we get that $\ex{\|W_{iid}\|^4\mathbf{1}_{\|W_{iid}\|>(1+\sqrt{c_*})^2+\delta}}=\mathcal{O}(N^{-k})$ for all $k$.
Since function $\phi^\prime\in\mathcal{C}_c^{\infty}$, there exists a nice constant $\kappa$ such that $|\phi^\prime(\lambda)|<\kappa$ for all $\lambda$ and $\phi^\prime(\lambda)=0$ for all $\lambda>b+2\epsilon$. We deduce from this it exists a nice constant $\kappa$ such that $\|\eta(W_{f,N}W_{p,N}^*W_{p,N}W_{f,N}^*)\|<~\kappa$ for each $N$. From what about we conclude that $\psi_1=\mathcal{O}(N^{-2})$.

As for the second term ($\psi_2$) of the r.h.s of (\ref{ineq:var_tr_phi}), we write
\begin{align*}
\psi_2=\dfrac{\kappa }{N}\ex\left\{\tr W_p^*W_pW_p^*W_pW_f^* \left(\phi^\prime(W_fW_p^*W_pW_f^*)\right)^2W_f\right\}\\
\le \kappa \ex\left\{\| W_p\|^2\dfrac{1}{N}\tr\left(\phi^\prime(W_fW_p^*W_pW_f^*)\right)^2W_fW_p^*W_pW_f^*\right\}
\end{align*}
It is easy to see that $\psi_2$ can be evaluated as $\psi_1$, leading to the conclusion that 
$\psi_2=\mathcal{O}(N^{-2})$. Therefore, we have checked that 
\begin{align*}
\var\{\tr\phi(W_fW_p^*W_pW_f^*)\}=\mathcal{O}\left(\dfrac{1}{N^2}\right).
\end{align*} 

Now we can complete the proof of Theorem~\ref{th:no_eigen} as in \cite{capitaine-donati-martin-feral-2009}. For this we apply the classical Markov inequality and combine what above
\begin{multline*}
\mathbf{P}\left\{\dfrac{1}{ML}\tr\phi(W_fW_p^*W_pW_f^*)>\dfrac{1}{N^{4/3}}\right\}
\le N^{8/3}\mathbb{E}\left\{\left(\dfrac{1}{ML}\tr\phi(W_fW_p^*W_pW_f^*)\right)^2\right\}\\
=N^{8/3}\left(\var\left\{\dfrac{1}{ML}\tr\phi(W_fW_p^*W_pW_f^*)\right\}+\left(\mathbb{E}\left\{\dfrac{1}{ML}\tr\phi(W_fW_p^*W_pW_f^*)\right\}\right)^2\right)\\
=\mathcal{O}\left(\dfrac{1}{N^{4/3}}\right).
\end{multline*}
Applying Borel-Cantelli lemma, for $N$ large enough, we have with probability one 
\begin{align*}
\dfrac{1}{ML}\tr\phi(W_fW_p^*W_pW_f^**)\le\dfrac{1}{N^{4/3}}
\end{align*}
By the very definition of function $\phi$, the number of eigenvalues of matrix $W_fW_p^*W_pW_f^*$ lying in the interval $[\kappa_1,\kappa_2]$ is upper bounded by $\tr\phi(W_fW_p^*W_pW_f^*)\le\frac{1}{N^{1/3}}$. Since this number of eigenvalues is an integer, we conclude that with probability one there is no eigenvalues in the interval $[\kappa_1,\kappa_2]$ for each $N$ large enough.  $\blacksquare$ \\

\begin{figure}[ht!]
	\centering
	\par
        \includegraphics[scale=0.8]{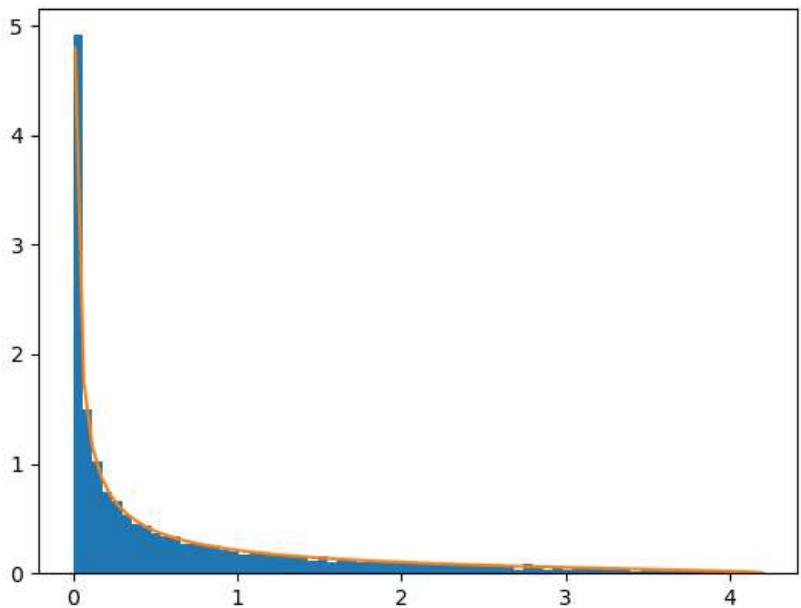}
	\caption{Histogram of the eigenvalues and graph of $g_N(x)$ for
            $M=500, N=1500, L=2$}
	\label{fig:histogram}
\end{figure}

We finally illustrate the above results by the following numerical experiment. 
$M,N,L$ are given by $M=500$, $N=1500$
and $L=2$ so that $c_N = 2/3$. The eigenvalues of matrix $R_N$ are defined by $\lambda_{k,N}  = 1/2 + \frac{\pi}{4} \cos\left(\frac{\pi (k-1)}{2M}\right)$
for $k=1, \ldots, M$. Matrix $R_N$ verifies $\frac{1}{M} \mathrm{Tr}(R_N) \simeq 1$. Fig. \ref{fig:histogram} 
represents the histogram of the eigenvalues of
a realization of  $W_{f,N} W_{p,N}^{*}  W_{p,N} W_{f,N}^{*}$ as well as the graph
of the density $g_N(x)$. We notice that the histogram and the graph of $g_N$ are in accordance, and
that, as expected, no eigenvalue of  $W_{f,N} W_{p,N}^{*} W_{p,N} W_{f,N}^{*}$  lies outside the 
support of $g_N$. 

	\section{Recovering the behaviour of the empirical eigenvalue distribution $\hat{\nu}_N$ using free probability tools}
\label{sec:free-probability}

        The purpose of this paragraph  is to show that it is possible to use free probability tools in order
        to characterize the limiting behaviour of the empirical eigenvalue distribution $\hat{\nu}_N$ of matrix
        $W_{f,N}W_{p,N}^{*}W_{p,N} W_{f,N}^{*}$. As the present paper is not focused on these kind of approach, we present briefly the following results and leave the details to the reader. \\

        The free probability approach is based on the following observations:
        \begin{itemize}
        \item Up to the zero eigenvalue, the eigenvalues of $W_{f,N}W_{p,N}^{*}W_{p,N} W_{f,N}^{*}$ coincide with the eigenvalues of $ W_{f,N}^{*}W_{f,N} W_{p,N}^{*}W_{p,N}$
        \item The matrices $W_{f,N}^{*}W_{f,N}$ and $W_{p,N}^{*}W_{p,N}$ are almost surely asymptotically free. Therefore, the eigenvalue distribution
          of $ W_{f,N}^{*}W_{f,N} W_{p,N}^{*}W_{p,N}$ converges towards the free multiplicative convolution product of the limit distributions of  $W_{f,N}^{*}W_{f,N}$ and $W_{p,N}^{*}W_{p,N}$. These two distributions appear to coincide both with the limit distribution of the well known random matrix model
          $\frac{1}{N} X_N^{*} (I_L \times R_N) X_N$ where $X_N$ is a $ML \times N$ complex Gaussian random matrix with unit variance i.i.d. entries.  
        \end{itemize}
        In the following, we follow the definitions of asymptotic freeness provided in \cite{hiai-petz-2000}
        (see in particular section 4.3) which need the existence of certain limit distributions. This is in contrast with
        the approach developed in the previous sections more focused on the behaviour of deterministic equivalents. We however mention that more recent free probability works (see e.g. \cite{mingo-speicher-book} and the references therein, \cite{belinschi-capitaine-2017}) allow to avoid the introduction of limit distributions, and would allow to recover the previous results on the deterministic equivalent $\nu_N$ of $\hat{\nu}_N$. \\

In order to be in accordance with \cite{hiai-petz-2000}, we thus formulate in this section the following 
assumption: 
\begin{assumption}
\label{as:limit-omegaN}
The empirical eigenvalue distribution $\omega_N = \frac{1}{M} \sum_{k=1}^{M} \delta_{\lambda_{k,N}}$
of matrix $R_N$ converges towards a limit distribution $\omega$
\end{assumption}
We remark that hypothesis \ref{eq:hypothesis-R-bis} implies that $\omega$ is compactly supported. Moreover, it can be shown that measures $(\mu_N)_{N \geq 1}$ and
$(\nu_N)_{N \geq 1}$ both converge weakly towards limits denoted $\mu$ and $\nu$ in this section. We also notice that
Lemma \ref{le:w+-x+-bounded} implies that $\mu$ and $\nu$ are compactly supported. It is also easily checked that the
Stieltjes transform $t(z)$ of $\mu$ verifies the equation
\begin{align}
\label{eq:equation-t}
t(z)=-\dfrac{1}{z}\int_{\mathbb{R}^{+}} \dfrac{\tau \, d\omega(\tau)}{1+\dfrac{c_*\tau t(z)}{1-zc_*^2t^2(z) }}
\end{align}
while the Stieltjes transform $t_{\nu}$ of $\nu$ is given by
\begin{equation}
  \label{eq:equation-tnu}
t_{\nu}(z) = -\frac{1}{z} - \frac{c_* t(z)^{2}}{1 - z (c_*t(z))^{2}}
\end{equation}
We recall that $c_*$ represents the limit of $c_N = \frac{ML}{N}$. In the following, we establish that
(\ref{eq:equation-t}) and (\ref{eq:equation-tnu}) can be obtained using free probability
technics. \\

Before going further, we first recall the main useful definitions introduced in \cite{hiai-petz-2000}.
        \begin{definition}
          \label{def:limit-joint-distribution}
          Consider a finite family of sequences of $N \times N$ possibly random matrices $((X_{i,N})_{N \geq 1})_{i=1, \ldots, r}$.
          Then $(X_{i,N})_{i=1, \ldots, r}$
          is said to have an almost sure joint limit if for each non commutative polynomial $P(x_1, \ldots,x_r)$
          in $r$ indeterminates, then $\frac{1}{N} \mathrm{Tr}P(X_{1,N}, \ldots,X_{r,N})$ converges almost surely towards $\mu(P)$
          where $\mu$ is a deterministic distribution defined on the set of all non commutative polynomials
          in $r$ indeterminates (i.e. $\mu$ is a linear form such that $\mu(1)=1$).
        \end{definition}
        We remark that if $r=1$ and $(X_{1,N})_{N \geq 1}$ are Hermitian matrices,
        the above condition is equivalent to the existence of a limit empirical eigenvalue distribution.
        \begin{definition}
          \label{def:asymptotic-freeness}
Consider $p$ families 
$(X^{(1)}_{i,N})_{i=1, \ldots, r_1}, \ldots, (X^{(p)}_{i,N})_{i=1, \ldots, r_p}$ of $N \times N$ 
possibly random matrices. Then, $X^{(1)}, \ldots, X^{(p)}$ are said to be almost surely asymptotically free
if the 2 following conditions hold:
\begin{itemize}
\item For each $q=1, \ldots, p$,  $(X^{(q)}_{i,N})_{i=1, \ldots, r_q}$ has an almost sure joint limit
\item $\forall m$, $i_1,\cdots, i_m\in \{1,2, \ldots, p \}$ with $i_1\ne i_2\ne\cdots\ne i_m$, and for each
non commutative polynomials $(P_j)_{j=1, \ldots, m}$ in $(r_{i_j})_{j=1, \ldots, m}$ indeterminates such that $\frac{1}{N}\tr(P_j(X^{i_j}_{1,N}, \ldots, X^{i_j}_{r_{i_j},N}))\rightarrow 0$ a.s. it holds that
	\begin{align}
	\dfrac{1}{N}\tr(P_1(X^{i_1}_{1,N}, \ldots, X^{i_1}_{r_{i_1},N}) \cdots P_m(X^{i_m}_{1,N}, \ldots, X^{i_m}_{r_{i_m},N}))\rightarrow 0 \quad a.s.
	\end{align}
\end{itemize} 
\end{definition}
We remark that when each family $X^{(q)}$ is reduced to a single sequence $(X^{(q)}_{N})_{N \geq 1}$ 
of $N \times N$ hermitian, or similar to hermitian matrices \footnote{in the sense that $X^{(q)}_N = U^{(q)}_N H^{(q)}_N (U^{(q)}_N)^{-1}$ for some
  $N \times N$  hermitian matrix $H^{(q)}_N$}, the almost sure freeness of $X^{(1)}, \ldots, X^{(p)}$ holds if
\begin{definition}
  \label{def:asymptotic-freeness-simplified}
\begin{itemize}
	\item For each $q=1, \ldots, p$, $(X^{(q)}_N)_{N \geq 1}$ has a limit eigenvalue distribution
	\item  $\forall m$, $i_1,\cdots, i_m\in \{1, 2, \ldots,p \}$ with $i_1\ne i_2\ne\cdots\ne i_m$, and for each 
1 variate polynomials $(P_j)_{j=1\cdots m}$ such that $\frac{1}{N}\tr(P_j(X^{i_j}_N))\rightarrow 0$ a.s. it holds that
	\begin{equation}
	\label{eq:free-matrices}
	\dfrac{1}{N}\tr(P_1(X^{(i_1)}_N)P_2(X^{(i_2)}_N) \cdots P_m(X^{(i_m)}_N))\rightarrow 0 \quad a.s.
	\end{equation}
\end{itemize}
\end{definition}
We also recall the definition of the $S$ transform of a probability measure, and recall that the 
$S$ transform of the free multiplicative convolution product of 2 probability measures is the product
of their $S$ transforms.          
\begin{definition}
\label{th:s-transf}
	Given a compactly supported probability measure $\mu$ carried by $\mathbb{R}^{+}$, we define $\psi_{\mu}(z)$ as the formal power series defined by
	\begin{align}\label{def:psi}
	\psi_{\mu}(z)=\sum_{k\ge1}z^k\int t^kd\mu(t)=\int\dfrac{zt}{1-zt}d\mu(t)
	\end{align}
Let $\chi_{\mu}$ be the unique function analytic in a neighbourhood of zero, satisfying
\begin{align}\label{def:chi}
\chi_{\mu}(\psi_{\mu}(z))=z
\end{align}
for $|z|$ small enough. The, we define the $S$ transform of $\mu$ as the function $S_{\mu}(z)$
defined in a neighbourhood of zero by 
\begin{align}\label{def:s-transform}
S_{\mu}(z)=\chi_{\mu}(z)\dfrac{1+z}{z}.
\end{align}
Moreover, if $\mu_1$ and $\mu_2$ are two compactly supported probability measures carried by $\mathbb{R}^{+}$, the S-transform $S_{\mu_1\boxtimes\mu_2}$ of $\mu_1\boxtimes\mu_2$ satisfies
\begin{align}\label{eq:s-transform_multip}
S_{\mu_1\boxtimes\mu_2}=S_{\mu_1}S_{\mu_2}.
\end{align}
\end{definition}
We are now in position to state the main result of this section. 
\begin{proposition}\label{th:asymp_free}
	Matrices $W_{f,N}^*W_{f,N}$ and $W_{p,N}^*W_{p,N}$ are almost surely asymptotically free.
\end{proposition}
\textbf{Proof.} We first notice that it possible to replace matrices $W_f$ and $W_p$ 
by finite rank perturbations because the very definition of almost sure asymptotic freeness is not affected 
by finite rank perturbations. We thus exchange $W_p$ and $W_f$ by $\tilde{W}_p = \frac{1}{\sqrt{N}} \tilde{Y}_p$ 
and $\tilde{W}_f = \frac{1}{\sqrt{N}} \tilde{Y}_f$ where $\tilde{Y}_p$ and $\tilde{Y}_f$ are defined by

\begin{equation}
  \tilde{Y}_p= \left( \begin{array}{cccccccc} y_1 & \ldots &\ldots & \ldots & \ldots & \ldots &  y_N \\
    y_2 & \ldots &\ldots & \ldots &\ldots &  y_N & y_1 \\
    y_3 & \ldots &\ldots &\ldots & y_{N} & y_1 & y_2 \\
    \vdots & \ldots & \ldots & \reflectbox{$\ddots$} & \reflectbox{$\ddots$} &  \reflectbox{$\ddots$} & \vdots \\
    y_L & \ldots & y_N & y_1 & y_2 & \ldots & y_{L-1} \end{array} \right)
  \end{equation}

\begin{equation}
  \tilde{Y}_f= \left( \begin{array}{cccccccccc} y_{L+1} & \ldots &\ldots  & \ldots &\ldots & \ldots  &  y_N & y_1 & \ldots & y_L \\
    y_{L+2} & \ldots  & \ldots & \ldots &\ldots  & y_{N} & y_1 &\ldots &  y_L & y_{L+1} \\
    y_{L+3} & \ldots &\ldots  & \ldots & y_N & y_1  &\ldots & y_{L} & y_{L+1} & y_{L+2} \\
    \vdots & \ldots  & \ldots & \reflectbox{$\ddots$} &\reflectbox{$\ddots$} &\ldots & \reflectbox{$\ddots$} & \reflectbox{$\ddots$} &  \reflectbox{$\ddots$} & \vdots \\
    y_{2L} & \ldots & y_N & y_1 & \ldots & y_{L} & y_{L+1} & y_{L+2} & \ldots & y_{2L-1} \end{array} \right)
\end{equation}

In other words, vectors $y_{N+1}, \ldots, y_{N+L-1}, \ldots, y_{N+2L-1}$ are replaced by vectors
$y_{1}, \ldots, y_{L-1}, \ldots, y_{2L-1}$. 
In order to simplify the notations, we still denote the above finite rank modifications by 
$Y_p, Y_f, W_p, W_f$. We define the $N\times N$ matrix $\Pi$ and $M\times N$ matrix $Y$ by
\begin{align}
\Pi=\begin{pmatrix}
0&\ldots&0&1\\
1&\ddots&&0\\
\vdots&\ddots&\ddots&\vdots\\
0&\ldots&1&0
\end{pmatrix}, \, \text{and }
Y=(y_1,y_2,\ldots,y_N)
\end{align} 
and rewrite $Y_p$ (and $Y_f$ respectively) as
\begin{align}
Y_p=\begin{pmatrix}
Y\\
Y\Pi\\
\vdots\\
Y\Pi^{L-1}
\end{pmatrix},\quad
Y_f=\begin{pmatrix}
Y\Pi^{L}\\
Y\Pi^{L+1}\\
\vdots\\
Y\Pi^{2L-1}
\end{pmatrix}
\end{align}
This allows us to obtain the useful expression for $W_p^*W_p$ and $W_f^*W_f$ 
\begin{eqnarray}
  \label{eq:expre-WpWp*}
  &W_p^*W_p=\sum_{k=0}^{L-1}\Pi^{*k}\left(\dfrac{Y^*Y}{N}\right)\Pi^k\\
   \label{eq:expre-WfWf*}
&W_f^*W_f=\sum_{k=L}^{2L-1}\Pi^{*k}\left(\dfrac{Y^*Y}{N}\right)\Pi^k
\end{eqnarray}
Since $N^{-1}Y^*Y$ can be written as $N^{-1}Y^*_{iid}R_NY_{iid}$, where $Y_{iid}$ has i.i.d. Gaussian entries, 
the hermitian matrix $N^{-1}Y^*Y$ is unitarily invariant. Moreover,  Assumption~\ref{as:limit-omegaN} implies that 
$N^{-1}Y^*Y$ has a limit distribution while it is easily checked that the family 
 $\{I,\Pi^*,\Pi,\ldots,\Pi^{*2L-1},\Pi^{2L-1}\}$ has the same property.  This and Theorem 4.3.5 in \cite{hiai-petz-2000} leads to the conclusion that $Y^*Y/N$ and $\{I,\Pi^*,\Pi,\ldots,\Pi^{*2L-1},\Pi^{2L-1}\}$ are almost surely asymptotically free. Proposition \ref{th:asymp_free} thus appears to be an immediate consequence of the 
following Lemma adapted from Lemma 6 in \cite{Evans_Tse_2000}. In order to make the connections between Lemma \ref{evans-tse} and Lemma 6 in \cite{Evans_Tse_2000},
we use nearly the same notations than in \cite{Evans_Tse_2000} in the following statement. 

\begin{lemma}\label{evans-tse}
	We consider a sequence of $N \times N$ hermitian random matrices $(X^{N})_{N \geq 1}$ and $N \times N$ deterministic 
matrices $U_{1}^{N},W_{1}^N,\ldots, U_{m}^N,W_{m}^N$ such that $X_N$ and $\{ U_{1}^{N},W_{1}^N,\ldots, U_{m}^N,W_{m}^N \}$ are almost surely asymptotically free. Then, if $U_{1}^{N},W_{1}^N,\ldots, U_{m}^N,W_{m}^N$ satisfy 
	\begin{align}\label{eq:cond_U_W}
	U^N_iW^N_i=W^N_iU^N_i= I_N
	\end{align}
	for each $i=1, \ldots, m$ as well as  $\frac{1}{N}\tr (U^N_iW^N_j)=\delta_{i-j}$ for all $i,j=1\ldots m$, then the random matrices $U^N_1X^NW^N_1,\ldots, U^N_mX^NW^N_m$ are almost surely asymptotically free.
\end{lemma}
{\bf Proof.} We prove Lemma \ref{evans-tse} by following step by step the proof from \cite{Evans_Tse_2000}.
For simplicity we omit index $N$ below.  Due to (\ref{eq:cond_U_W}) we have $W_i = U_i^{-1}$ so that matrices $(U_i X W_i)_{i=1, \ldots, m}$
are similar to the hermitian matrix $X$. We have thus to verify the 2 items of Definition \ref{def:asymptotic-freeness-simplified}. The first
item is obvious. To check condition (\ref{eq:free-matrices}), we consider any $k$, indexes $i_1,\cdots,i_k$ with $i_1\ne\cdots\ne i_k$ and polynomials $P_j$ such that $\frac{1}{n}\tr(P_j(U_{i_j}XW_{i_j}))\rightarrow 0$ a.s. Using again (\ref{eq:cond_U_W}) it is clear that $P_j(U_{i_j}XW_{i_j})=U_{i_j}P_j(X)W_{i_j}$ and, as a
consequence, $\frac{1}{n}\tr(P_j(X))\rightarrow 0$ a.s. We define $\eta_N$ as 
\begin{multline*}
	\eta_N=\dfrac{1}{N}\tr(P_1(U_{i_1}XW_{i_1})P_2(U_{i_2}XW_{i_2})\cdots (U_{i_k}XW_{i_k}))=\\
	\dfrac{1}{N}\tr(U_{i_1}P_1(X)W_{i_1}U_{i_2}P_2(X)W_{i_2}\cdots U_{i_k}P_k(X)W_{i_k})
	=\dfrac{1}{N}\tr\left(\prod_{j=1}^{k}W_{i_{j-1}}U_{i_j}P_j(X)\right)
\end{multline*} 
where $i_0=i_k$. If $i_{1}\ne i_k$ then by assumption $\frac{1}{n}\tr(W_{i_{j-1}}U_{i_j})=0$ for $j=1, \ldots, m$. As we  also have $\frac{1}{n}\tr(P_j(X))\rightarrow 0$ a.s, the almost sure asymptotic freeness of $X$ and $\{U_1,W_1,\cdots,U_m,W_m\}$ leads to the conclusion that $\eta_N \rightarrow 0$ a.s. In the case when $i_1=i_k$ we have $W_{i_k}U_{i_1}=I_N$ and the same conclusion holds.
$\square$

By taking $X=\frac{YY^*}{N}$, $U_i=\Pi^{*i-1}$ and $W_i=\Pi^{i-1}$, Lemma~\ref{evans-tse} gives us immediately that \\  $\frac{Y^*Y}{N},\Pi^*(\frac{Y^*Y}{N})\Pi,\ldots,\Pi^{*2L-1}(\frac{Y^*Y}{N})\Pi^{2L-1}$ are almost surely asymptotically free. Using the expression (\ref{eq:expre-WpWp*}, \ref{eq:expre-WfWf*}) of  $W_p^*W_p$ and $W_f^*W_f$, we obtain that  $W_p^*W_p$ and $W_f^*W_f$  are almost surely asymptotically free. $\blacksquare$ \\

We also deduce that the limit distributions of  $W_p^*W_p$ and $W_f^*W_f$ both coincide with the additive free convolution product of $L$ copies of the well
known limit distribution of $\frac{Y^*Y}{N}$. It is easily seen that the Stieljes transform, denoted $t_{MP}(z)$ in the following, of this free addditive convolution product is solution of the familiar equation 
\begin{align}\label{eq:stieltjes_MP}
t_{MP}(z)=-\dfrac{1}{z-c_*\int\dfrac{\tau\omega (d\tau)}{1+\tau t_{MP}(z)}}
\end{align} 
In the following, we denote by $\mu_{MP}$ the corresponding probability measure. 
It is clear that (\ref{eq:stieltjes_MP}) coincides with the equation verified by the Stieltjes transform of the limit eigenvalue distribution
of the random matrix  $\frac{1}{N} X_N^{*} (I_L \times R_N) X_N$ where $X_N$ is a $ML \times N$ complex Gaussian random matrix with unit variance i.i.d. entries.
We note that this result could also be easily obtained using the Gaussian technics developed in \cite{L:15} in the case where $R_N$ is reduced to a multiple
of $I_M$. \\

According to Proposition~\ref{th:asymp_free}, the limit eigenvalue distribution of $W_{f,N}^*W_{f,N}W_{p,N}^*W_{p,N}$ is $\mu_{MP}\boxtimes\mu_{MP}$.
In the following, we denote by $\tilde{\nu}$ this measure and by $\tilde{f}(z)$ its Stieltjes transform. To find an equation satisfied by $\tilde{f}(z)$, we  use (\ref{eq:s-transform_multip}). 
(\ref{def:s-transform}) and (\ref{eq:s-transform_multip}) give us immediately
\begin{align*}
\chi_{\tilde{\nu}}(z)=\dfrac{1+z}{z}\chi_{MP}^2(z)
\end{align*}
By replacing here $z$ with $\psi_{\tilde{\nu}}(z)$ and taking into account (\ref{def:chi}) we obtain
\begin{align}\label{eq:psi_nu-chi_MP}
z=\dfrac{1+\psi_{\tilde{\nu}}(z)}{\psi_{\tilde{\nu}}(z)}\chi_{MP}^2(\psi_{\tilde{\nu}}(z))
\end{align}
We notice that by definition (\ref{def:psi}), we have
\begin{align}\label{eq:psy-stieltjes}
\psi_{\tilde{\nu}}(z)=\int\dfrac{zt}{1-zt}d\tilde{\nu}(t)=\int\dfrac{d\tilde{\nu}(t)}{1-zt}-1=-\dfrac{1}{z}\tilde{f}\left(\dfrac{1}{z}\right)-1
\end{align}
Putting this into (\ref{eq:psi_nu-chi_MP}) and replacing $z$ with $\frac{1}{z}$ give us
\begin{align*}
\dfrac{z^2\tilde{f}(z)}{1+z\tilde{f}(z)}\chi_{MP}^2\left(\psi_{\tilde{\nu}}\left(\dfrac{1}{z}\right)\right)=1
\end{align*}
From this, it is straightforward to obtain the expression of $\tilde{f}(z)$. For more convenience, we introduce the function $g(z)=\chi_{MP}^2(\psi_{\tilde{\nu}}(z^{-1}))$ which is analytic in the neighbourhood of infinity. It holds that 
\begin{equation}
  \label{eq:expre-tildef-g}
\tilde{f}(z)=\left(z^2g^2(z)-z\right)^{-1}
\end{equation}  
It remains to determine $g(z)$. For this we use (\ref{eq:psy-stieltjes}) for $\psi_{MP}$, $t_{MP}$ and replace $z$ with $\chi_{MP}(z)$. Then (\ref{def:chi}) gives 
\begin{align*}
z=-1-\dfrac{1}{\chi_{MP}(z)}t_{MP}\left(\dfrac{1}{\chi_{MP}(z)}\right)
\Rightarrow t_{MP}(\chi_{MP}^{-1}(z))=-(1+z)\chi_{MP}(z)
\end{align*}
To obtain the equation for $\chi_{MP}$ it is sufficient to use the above expression of $t_{MP}(\chi_{MP}^{-1}(z))$, and to plug it in (\ref{eq:stieltjes_MP}) with $z=\chi_{MP}^{-1}(z)$. Therefore, we obtain that 
\begin{align*}
(1+z)\chi_{MP}(z)=\dfrac{1}{\dfrac{1}{\chi_{MP}(z)}-c_*\int\dfrac{\tau d\omega_(\tau)}{1-\tau(1+z)\chi_{MP}(z)}}
\end{align*}
After simple algebra we get that
\begin{align*}
\dfrac{z}{(1+z)\chi_{MP}(z)}=c_*\int\dfrac{\tau d\omega_(\tau)}{1-\tau(1+z)\chi_{MP}(z)}
\end{align*} 
We finally replace $z$ by $\psi_{\tilde{\nu}}(z^{-1})$. With (\ref{eq:psi_nu-chi_MP}) it is easy to see that the l.h.s. equal to $zg(z)$. To treat the r.h.s. we use again (\ref{eq:psi_nu-chi_MP}) to obtain that  $\psi_{\tilde{\nu}}(z^{-1})=zg^2(z)(1-zg^2(z))$, and get that
\begin{align}
  \label{eq:equation-g}
g(z)=
\dfrac{1}{z}\int_{\mathbb{R}^{+}} \dfrac{c_*\tau \, d\omega(\tau)}{1-\dfrac{\tau g(z)}{1-zg^2(z)}}
\end{align}
Now we recall the equation obtained above for $t(z)$
\begin{align}\label{eq:t}
t(z)=-\dfrac{1}{z}\int\dfrac{\tau\omega(d\tau)}{1+\dfrac{c_*\tau t(z)}{1-zc_*^2t^2(z) }}
\end{align}
The equations (\ref{eq:equation-g}) and (\ref{eq:t}) are identical up to factor $-c_*$. Since it can be shown that Eq. (\ref{eq:t}) has a unique solution on the set of Stieltjes transforms, we obtain that $g(z)=-c_*t(z)$. Therefore, (\ref{eq:expre-tildef-g})
leads to the equation
$$
\tilde{f}(z) = -\frac{1}{z\left[1 - z(c_*t(z))^{2}\right]}
$$
The Stieltjes transform of the limit eigenvalue distribution of $W_{f} W_p^{*} W_p W_{f}^{*}$ is clearly equal
to $\frac{1}{c_*} \left( \tilde{f}(z) + \frac{1-c_*}{z} \right)$. Using the expression (\ref{eq:equation-tnu}) of $t_{\nu}(z)$, we obtain immediately that
$$
\frac{1}{c_*} \left( \tilde{f}(z) + \frac{1-c_*}{z} \right) = t_{\nu}(z)
$$
We have thus proved that the limit eigenvalue distribution of  $W_{f} W_p^{*} W_p W_{f}^{*}$ can be evaluated using free probability technics.

\section*{Acknowledgments}
The authors thank warmly Prof. Leonid Pastur for fruitful discussions.

\end{document}